%% file: Picklo_LSIACMRA.tex
\documentclass{siamart190516}

\usepackage{geometry}
\usepackage{chapterbib}
\usepackage[utf8]{inputenc}
\usepackage{tabularx}
\usepackage{tikz}
\usepackage{amsmath,amssymb,amsfonts}
\usetikzlibrary{patterns,patterns.meta}
\usetikzlibrary{arrows,arrows.meta}

\usepackage{algorithm}
\usepackage{algpseudocode}
\usepackage{mathtools}
\usepackage{multirow}
\usepackage{pgfplots}
\def\addlegendimage{\csname pgfplots@addlegendimage\endcsname}
\newenvironment{customlegend}[1][]{%
    \begingroup
    \csname pgfplots@init@cleared@structures\endcsname
    \pgfplotsset{#1}%
}{%
    \csname pgfplots@createlegend\endcsname
    \endgroup
}%

\usepackage{mathrsfs}
\usepackage{graphicx}
\usepackage{pstricks-add}
\usepackage{pgf}
\usepackage{lineno}
\usepackage[toc,page]{appendix}
\usepackage{bbm}
\usepackage{bbold}
\usepackage{lscape,rotating,pdflscape}
\usepackage{enumerate}
\usepackage{float,multirow}
\usepackage{array}
 \usepackage{relsize}
\usepackage{setspace}
\usepackage{colortbl}
\usepackage{mathdots}
\usepackage{xcolor}
 \usepackage{blkarray,stackengine}
 \usepackage{footnote}
 \usepackage[bottom]{footmisc}

 \newcommand{\norm}[1]{\left\lVert#1\right\rVert}
 
\usetikzlibrary{decorations.markings}
\definecolor{DarkRed}{rgb}{0.65,0.00,0.05}
\definecolor{structure}{rgb}{0.2,0.2,0.7}

\newcolumntype{M}[1]{>{\centering\arraybackslash}m{#1}}
\newcolumntype{N}{@{}m{0pt}@{}}

\title{Enhanced Multi-Resolution Analysis for Multi-Dimensional Data Utilizing Line Filtering Techniques}

\author{Matthew J. Picklo{\footnotemark[1] \footnotemark[2] \footnotemark[3]} \and
Jennifer K. Ryan{\footnotemark[1] \footnotemark[2] \footnotemark[4]}}

\footnotetext[1]{Applied Mathematics and Statistics, Colorado School of Mines, Golden, CO 80401, USA.}
\footnotetext[2]{Research supported by Air Force Office of Scientific Research (AFOSR), Computational Mathematics Program (Program Manager: Dr. Fariba Fahroo), under {grant number FA9550-20-1-0166.}}

\footnotetext[3]{Email: {mpicklo@mines.edu}}
\footnotetext[4]{{\it Corresponding Author.} {Email:} jkryan@mines.edu}

\date{\today}

\begin{document}

\maketitle

\begin{abstract}
 In this article we introduce Line Smoothness-Increasing Accuracy-Conserving Multi-Resolution Analysis\linebreak (LSIAC-MRA). This is a procedure for exploiting convolution kernel post-processors for obtaining more accurate {multi-dimensional multi-resolution analysis (MRA) in terms of error reduction.} This filtering-projection tool allows for the transition of data between different resolutions while simultaneously decreasing errors in the fine grid approximation.  It specifically allows for defining detail multi-wavelet coefficients when translating coarse data onto finer meshes.  These coefficients are usually not defined in such cases.  We show how to analytically evaluate the resulting convolutions and express the {filtered} approximation in a new basis. This {is done by} combining the filtering procedure with projection operators that allow for {computational} implementation of this scale transition procedure. Further, this procedure can be applied to piecewise constant approximations to functions, {as it provides error reduction.} We demonstrate the effectiveness of this technique in {two} and {three} dimensions.
\end{abstract}

\begin{keywords}
Discontinuous Galerkin, post-processing,
SIAC filtering, Line filtering, accuracy enhancement, error reduction, Multi-Resolution Analysis
\end{keywords}

\begin{AMS}{65M60}\end{AMS} 


\section{Introduction}

\input{Sections_LSIAC/Introduction}

\section{Background}

In the following sections we {present} the {ideas} necessary for {understanding} enhanced multi-resolution analysis for multi-dimensions.  To better illustrate the ideas, we begin with a discussion of the one-dimensional ideas {covered} in \cite{SIACMRA}.

\input{Sections_LSIAC/InputDataFormat}

\input{Sections_LSIAC/MRA}

\input{Sections_LSIAC/Filter_HD}

\section{Multi-dimensional Line SIAC-MRA}

\input{Sections_LSIAC/SIAC_MRA_HD}

\input{Sections_LSIAC/Implementation}

\input{Sections_LSIAC/Numerical_Results}

\input{Sections_LSIAC/Conclusions}
\input{Sections_LSIAC/Acknowledgements}

\section*{Statement of Contribution}
The ideas in this article were conceived {by} Jennifer Ryan, who is also responsible for the presentation of information, writing the introduction and conclusion sections, and final edits.  Matt Picklo performed the numerical experiments and wrote the majority of the article.

\bibliographystyle{siamplain.bst}
\bibliography{Picklo_LSIACMRA}
\newpage
\end{document}


\maketitle

\section{SIAC-MRA: Filtered Basis Functions \& Coefficients:}
In these supplementary materials we detail the $\chi-$functions and associated coefficients appearing in the simplified filtered approximations. In addition, figures are provided delineating the different domains of definition for these functions. Denote by $\zeta$ the mapping of the filtering point to the reference element $[-1,1]$, i.e. $\zeta_x=\frac{2}{h}(x-x_i)$, $\zeta_y=\frac{2}{h}(y-y_j)$, $\zeta_z=\frac{2}{h}(z-z_k)$.  

\subsection{One-dimensional SIAC-MRA:} We first recall that the filtered approximation given in Ryan, (2021) is of the form 
\[
u_h^{\star}(x)=\begin{cases}
\sum_{m=0}^{p+1}d_m^j\chi^m_L(\zeta_x),\hspace{1cm} &\zeta_x\in(-1,0)\\
\sum_{m=0}^{p+1}\tilde{d}_m^j\chi^m_R(\zeta_x),\hspace{1cm} &\zeta_x\in(0,1)
\end{cases},
\]
with the $\chi$ functions and coefficients given in Table S\hspace{-0.05em}\ref{tab:1Dchi}.

\input{Supplementary_Material/1D_Basis_Table}

\subsection{Two-dimensional LSIAC-MRA}
There are seven different forms for the coefficients included in the two-dimensional L-SIAC case with the basis given in Table S\hspace{-0.05em}\ref{tab:2D_LSIAC_basis}. These coefficients are of the form
\[
a^{S,P}_{m,r}=a_{m,r}(q_x,q_y)=\sum_{\gamma=-p}^pc_{\gamma}u^{i+\gamma+qx,j+\gamma+qy}_{m,r},
\]
 where the arguments $(q_x,q_y)$ are determined by the region $S$ and position $P$ of the filter point. A table of these dependencies by region and position is provided in Table S\hspace{-0.05em}\ref{tab:2D_LSIAC_coefficient}. The regions themselves are depicted in Figure S\hspace{-0.05em}\ref{fig:2D_basis_function_domain}. An algorithm for the LSIAC-MRA procedure is included in Algorithm S\hspace{-0.05em}\ref{alg:LSIAC}.

\input{Supplementary_Material/2D_L-SIAC_Basis_Table}

\input{Supplementary_Material/2D_L-SIAC_coefficients}
\input{Supplementary_Material/2D_LSIAC_algorithm}

\newpage
\subsection{Three-dimensional LSIAC-MRA}
\begin{figure}
    \centering
    \begin{tabular}{|c|}
             \hline \includegraphics[width=0.9\linewidth]{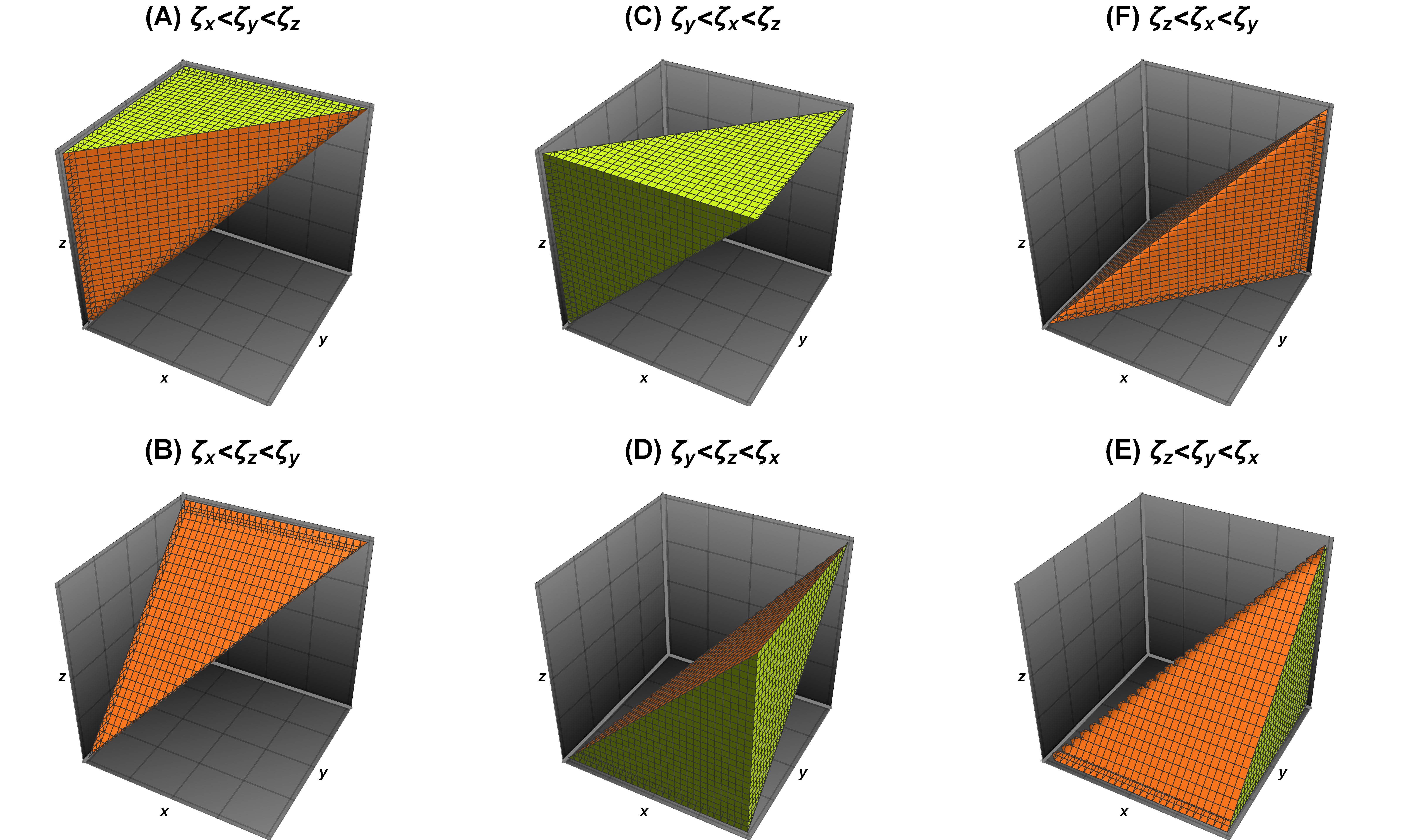}
 \\ \hline
              \includegraphics[width=0.9\linewidth]{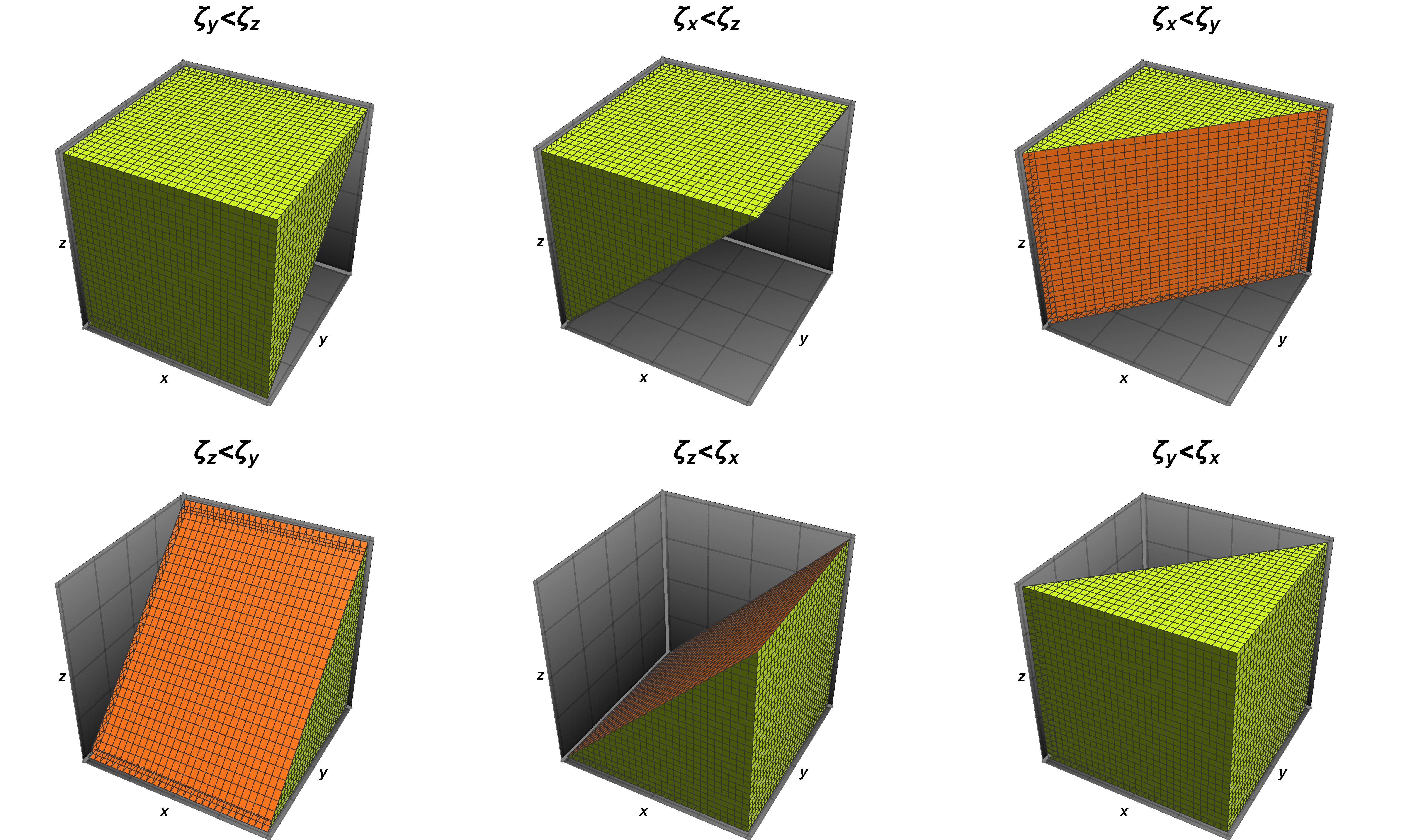}
\\ \hline
   \end{tabular}
    \caption{    \label{fig:3D_Regions} Depiction of how the domain of dependence of the filtered approximation form varies from octant to octant. Octants 1 and 8 have six regions (Top), while octants 2 through 7 have only a pair of regions each (Bottom). The form of these regions is a consequence of the filtering coordinate inequalities given in Table S\hspace{-0.05em}\ref{tab:3D_Classification} coupled with the relative location of each octant in the reference element.}
\end{figure}

In three-dimensions, the tensor product expansion functions are of the form \[\chi^{m,r,\ell}_{S,P}=C\int_{-1}^1\phi_m(\text{arg}_x)\phi_r(\text{arg}_y)\phi_{\ell}(\text{arg}_{z})\]
where the row of the table determines the Region $S$ and Position $P=\{L,LM,RM,R\}$. The regions are described in Tables S\hspace{-0.05em}\ref{tab:3D_Regions_Present} and S\hspace{-0.05em}\ref{tab:3D_Classification} with visual depictions in Figures S\hspace{-0.05em}\ref{fig:2D_basis_function_domain} and S\hspace{-0.05em}\ref{fig:3D_Regions}. The coefficient $C$ is given in the coefficient column, and the arguments of the Legendre polynomials are given in the appropriate argument column. These can be found in Table S\hspace{-0.05em}\ref{tab:3D_LSIAC_Basis}. There are fifteen different forms for the coefficients included in the 3D L-SIAC case. These coefficients are of the form
\[
a^{S,P}_{m,r,\ell}=a_{m,r,\ell}(q_x,q_y,q_z)=\sum_{\gamma=-p}^pc_{\gamma}u^{i+\gamma+qx,j+\gamma+qy,k+\gamma+qz}_{m,r,\ell},
\]
 where the coefficients $(q_x,q_y,q_z)$ are determined by the region $S$ and position $P$ of the filter point. A table of these dependencies by region and position is provided in Table S\hspace{-0.05em}\ref{tab:3D_LSIAC_coefficient}. 
  \input{Supplementary_Material/3D_L-SIAC_Classification_Scheme}

\input{Supplementary_Material/3D_L-SIAC_Basis_Table}

\input{Supplementary_Material/3D_L-SIAC_coefficients}

%% file: Sections_LSIAC/Introduction.tex
Approximating multi-scale phenomena accurately and efficiently has its significance illustrated in such areas as turbulence modeling and kinetics among many other applications.  Of considerable importance is the ability to match data on a coarse grid with that on a fine grid.  Multi-resolution analysis (MRA) is one technique that is useful in accomplishing such a task {\cite{Alp93,Hov10MS,Harten,M_ller_2003}}.  MRA expresses the approximation in terms of averages and differences, with the differences being the details between a fine grid approximation and coarse grid approximation.  The MRA approximation usually starts at the finest grid and moves to a coarser grid.  The reverse is a challenging procedure as the difference coefficients {that allow for higher levels of resolution do not exist}.  In this paper, we introduce a procedure which allows for approximating the difference coefficients {as well as error reduction when moving to a finer grid}.  To {accomplish this}, we introduce Line Smoothness-Increasing Accuracy-Conserving multi-resolution analysis (LSIAC-MRA).  As the name suggests, it relies on the Line SIAC filter found {in} \cite{LSIAC} that is an improvement over the tensor-product SIAC filters {\cite{CLSS2,ICERM,RSA}}. {Previously, SIAC filters were used for one-dimensional MRA in \cite{SIACMRA}}.  We present the necessary operators and algorithms to accomplish this task as well as demonstrating its effectiveness on two- and three-dimensional functions.  Aside from demonstrating that LSIAC-MRA can reduce the errors with mesh refinement, we show that it is also effective on piecewise constant approximations. 

Multi-resolution Analysis (MRA) is a useful technique for moving data between successive mesh \linebreak refinements. {Initially introduced in the 1980s by Meyer and Mallat \cite{ Mallat,Meyer}, the MRA framework was extended to discrete data representations by Harten \cite{Harten} and expanded with the addition of multi-wavelets by Alpert \cite{Alp93}. MRA-based wavelet techniques have found application in a variety of numerical methods for partial differential equations \cite{Alp02BGV,Wavelets2,Hov10MS,M_ller_2003,Dahmen1996}. A survey of wavelet-based methods in fluid mechanics can be found in \cite{Schneider}. } {MRA} is based on the idea that an approximation can be represented as {a linear combination of} scaling functions, which are typically the approximation basis, and wavelets. The scaling functions give the averages of the approximation over a coarse grid and the wavelet information contains the information necessary to move to a finer grid.  Mathematically, the expression is the following:
\begin{equation}
\label{eq:DGMRA}
\underbrace{u_h^f(x,t)}_{\text{fine grid approximation}} = \underbrace{u_h^c(x,t)}_{\text{coarse grid information}} + \underbrace{\sum_{j=1}^N\, \sum_{k=0}^p\, d_{k,j}^N(t)\psi^c_{k,j}(x)}_{\text{details necessary to move to a finer grid}}.
\end{equation}
{In this article, we use the basis functions for the multi-wavelet space given by Alpert, $\psi^c_{k,j}(x)$}. {These} multi-wavelets {are associated with} the piecewise orthonormal Legendre approximation {given} in \cite{Alp93}.  In this\newpage \noindent expression, both the scaling functions and wavelets are piecewise polynomials of degree $p$.   {We expect that this procedure can be extended to generalized wavelets.}

In creating this variant of Nystr\"{o}m reconstruction \cite{Nystrom}, we utilize the LSIAC filter \cite{LSIAC} for multi-dimensional data. Nystr\"{o}m reconstruction also uses convolution for accurate reconstructions.  For LSIAC-MRA, it allows for reducing the computational cost of the multi-dimensional filtering operation by introducing a one-dimensional support for the convolution kernel.

Pairing Smoothness-Increasing Accuracy-Conserving (SIAC) post-processing with the multi-resolution analysis presented in \cite{Alp93,Hov10MS} is natural. Although SIAC has broad applicability, it has mostly been developed for discontinuous Galerkin (DG) methods, which rely on a piecewise polynomial approximation basis. {An introduction} to this pairing in one-dimension can be found in \cite{SIACMRA}. However, the {multi}-dimensional equivalent is challenging as discussed in Section \ref{sec:Implementation}. Therefore in order to get to this stage we discuss the necessary background in multi-resolution analysis (Section \ref{sec:mra}), {followed by Line SIAC filters} (Section \ref{sec:SIACmultid}) before introducing LSIAC-MRA as a means to approximate difference coefficients.  Demonstration of the effectiveness of this technique can be found in Section \ref{sec:NumericalResults}.

%% file: Sections_LSIAC/InputDataFormat.tex
\subsection{Input Data Format}
For ease of discussion, we will assume our data is given in a modal format.  This generalizes to the nodal format by simply assuming that the approximation basis is the Lagrange polynomial basis. Therefore, we describe the procedure for constructing projected approximations for higher dimensional quadrilateral and hexahedral meshes under consideration. In the following, we {discuss the definitions necessary to implement LSIAC-MRA.}

 For a given function $u(\mathbf{x}),$  we wish to construct a piecewise-polynomial approximation $u_h(\mathbf{x})$. To do so, we utilize the following notation: A multi-index $\alpha=(\alpha_1,\hdots,\alpha_n)$ is an $n$-tuple of non-negative integers. Define $|\alpha|=\sum_{i=1}^n \alpha_i$, and  $\partial^{\alpha}=\partial^{\alpha_1}\hdots \partial^{\alpha_n}$, where $\partial^{\alpha_i}=\frac{\partial^{\alpha_i} }{\partial x^{\alpha_i}}$. Let $\mathbf{x}=(x_1,\hdots,x_d)$. Consider a $d$-dimensional rectangular domain and its partition into $N^d$ $d$-dimensional rectangular elements: $\Omega=\cup_{\beta \in \mathbb{B}} {\mathcal{T}}_{\beta}$, where the indexing set is given by
\[\mathbb{B}=\{\beta\in \mathbb{N}^d: \norm{\beta}_{\ell^{\infty}}\leq N\}.
\]
{Utilizing the form of} a discontinuous Galerkin type approximation, we define an approximation space over our partition that will contain piecewise {polynomials} in $d$ variables up to a given degree $p$:
\[V_h^p=\{v \in L^2\;:\; v\in\mathbb{P}^p({\mathcal{T}}_{\beta}),\; \beta\in\mathbb{B}\}, \]
{where $h$ is typically associated with the mesh size.}  We {approximation $u(\mathbf{x})$ via $L^2-$projection, that is on a given element ${\mathcal{T}}_{\beta}$, we write }
\[
{u_h(\mathbf{x})\Big|_{{\mathcal{T}}_{\beta}}=\sum_{\norm{\alpha}_{\ell^{\infty}}\leq p}u^{\beta}_{\alpha}\phi^{\beta}_{\alpha}(\mathbf{x}),\qquad \langle u,  \phi^{\beta}_{\alpha} \rangle_{{\mathcal{T}}_{\beta}}=\langle u_h,\phi^{\beta}_{\alpha} \rangle_{{\mathcal{T}}_{\beta}},}
\]
where $\{\phi^{\beta}_{\alpha}\}_{\norm{\alpha}_{\ell^{\infty}}\leq p}$ is a basis for our approximation space on ${\mathcal{T}}_{\beta}$. {Note the standard approximation space notation above conflicts with the Multi-resolution analysis notation to be used below. For our purposes we equate $V_h^p$ and $V_n^p$ when $N=2^n$.}

%% file: Sections_LSIAC/MRA.tex
\subsection{Multi-Resolution Analysis}
\label{sec:mra}
Multi-resolution analysis (MRA) {introduced in \cite{Mallat,Meyer} provides} a framework in which to analyze modal approximations under scale transition. Specifically, MRA views {approximations} as belonging to a hierarchy of nested approximation spaces, where transition between scales is simply the addition or removal of finer-detail basis functions called multi-wavelets. To better express the utility of this setting as it relates to our modal projections, we provide a brief summary of MRA as detailed in \cite{SIACMRA}. More in depth introductions to multi-resolution analyses generated by scaling functions can be found in {\cite{Alp93,Alp02BGV,wavelets,Wavelets2}}.  For ease of presentation, the one-dimensional case is given.

Consider a nested sequence of approximation spaces
\[ V^p_{0}\subset V^p_{1}\subset\hdots\subset V^p_n \subset \hdots, \]
where each approximation space is given by
\[
V^p_n=\{v: v\in \mathbb{P}^p(I_j^n),\;j=0,\hdots,2^n-1\},
\]
with 
\[I^n_j=(-1+2^{-n+1}j,-1+2^{-n+1}(j+1)),\quad j=0,\dots,2^n-1.
\]

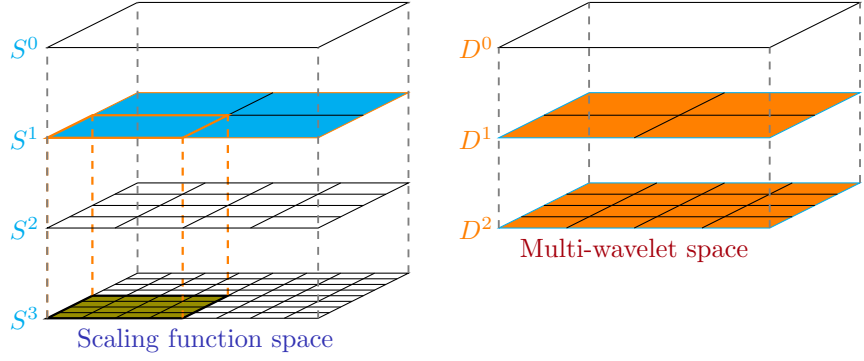
\begin{figure}[ht!]
\begin{center}

\begin{tikzpicture}[scale=0.6]
\draw (-8,3) -- (-2,3) --(-4,2) -- (-10,2) -- (-8,3);

\draw [color=orange, fill= cyan] (-8,1) -- (-2,1) --(-4,0) -- (-10,0) -- (-8,1);
\draw (-9,0.5) -- (-3,0.5);
\draw (-5,1) --(-7,0);

\draw[line width =0.3mm, fill=olive](-7,-4)--(-10,-4)--(-9,-3.5)--(-6,-3.5)--(-7,-4);
\draw[line width =0.3mm, color=orange](-7,0)--(-10,0)--(-9,0.5)--(-6,0.5)--(-7,0);
\draw[line width =0.3mm, color=orange, dashed] (-10,-4)--(-10,0);
\draw[line width =0.3mm, color=orange, dashed] (-7,-4)--(-7,0);
\draw[line width =0.3mm, color=orange, dashed] (-9,-3.5)--(-9,0.5);
\draw[line width =0.3mm, color=orange, dashed] (-6,-3.5)--(-6,0.5);

\draw (-8,-1) -- (-2,-1) --(-4,-2) -- (-10,-2) -- (-8,-1);
\draw (-9,-1.5) -- (-3,-1.5);
\draw (-5,-1) --(-7,-2);
\draw (-8.5,-1.25) --(-2.5,-1.25);
\draw (-9.5,-1.75) --(-3.5,-1.75);
\draw (-3.5,-1) -- (-5.5,-2);
\draw (-6.5,-1) -- (-8.5,-2);

\draw  (-8,-3) -- (-2,-3) --(-4,-4) -- (-10,-4) -- (-8,-3);
\draw (-9,-3.5) -- (-3,-3.5);
\draw (-5,-3) --(-7,-4);
\draw (-8.25,-3.125) --(-2.25,-3.125);
\draw (-8.5,-3.25) --(-2.5,-3.25);
\draw (-8.75,-3.375) --(-2.75,-3.375);
\draw (-9.25,-3.625)--(-3.25,-3.625);
\draw (-9.5,-3.75) --(-3.5,-3.75);
\draw (-9.75,-3.875)--(-3.75,-3.875);
\draw (-3.5,-3) -- (-5.5,-4);
\draw (-4.25,-3) -- (-6.25,-4);
\draw (-6.5,-3) -- (-8.5,-4);
\draw (-7.25,-3) -- (-9.25,-4);
\draw(-5.75,-3) --(-7.75,-4);
\draw(-2.75,-3) -- (-4.75,-4);

\draw[line width =0.25mm, color=gray, dashed] (-8,-3)--(-8,3);
\draw[line width =0.25mm, color=gray, dashed] (-2,-3)--(-2,3);
\draw[line width =0.25mm, color=gray, dashed] (-4,-4)--(-4,2);
\draw[line width =0.25mm, color=gray, dashed] (-10,-4)--(-10,2);
\draw (-6.5,-4.5) node{\color{structure} Scaling function space};
\draw (-10.5,-4) node {\color{cyan} $S^3$};
\draw (-10.5,-2) node {\color{cyan} $S^2$};
\draw (-10.5,0) node {\color{cyan} $S^1$};
\draw (-10.5,2) node {\color{cyan} $S^0$};

\draw (2,3) -- (8,3) --(6,2) -- (0,2) -- (2,3);

\draw [color=cyan, fill= orange] (2,1) -- (8,1) --(6,0) -- (0,0) -- (2,1);
\draw(1,0.5) -- (7,0.5);
\draw(3,0) -- (5,1);

\draw [color=cyan, fill= orange] (2,-1) -- (8,-1) --(6,-2) -- (0,-2) -- (2,-1);
\draw(1,-1.5) -- (7,-1.5);
\draw (5,-1)--(3,-2) ;
\draw(3.5,-1) -- (1.5,-2);
\draw(6.5,-1) --(4.5,-2);
\draw (1.5,-1.25) -- (7.5,-1.25);
\draw (0.5,-1.75) -- (6.5,-1.75);


\draw[line width =0.25mm, color=gray, dashed] (2,-1)--(2,3);
\draw[line width =0.25mm, color=gray, dashed] (8,-1)--(8,3);
\draw[line width =0.25mm, color=gray, dashed] (6,-2)--(6,2);
\draw[line width =0.25mm, color=gray, dashed] (0,-2)--(0,2);
\draw (3,-2.5) node{\color{DarkRed} Multi-wavelet space};
\draw (-0.5,-2) node{\color{orange} $D^2$};
\draw (-0.5,0) node{\color{orange} $D^1$};
\draw (-0.5,2) node{\color{orange} $D^0$};

\end{tikzpicture}
\end{center}
\caption{\label{fig:2dmw} Illustration of the two-dimensional multi-wavelet idea. {The approximation} at the finest level, $S^3$, can be represented as the approximation from the $S^1$ scaling function space plus the multi-wavelet coefficients from wavelet spaces $D^1$ and $D^2$.}
\end{figure}

{Denote} by $\{\phi_k(x)\}_{k=0}^p$ the orthonormal basis on the coarsest mesh consisting of one element, $n=0$. These basis functions are called \textit{scaling functions}, owing to fact that the systematic manner in which our approximation space and the underlying mesh is refined allows for construction of bases for any $V^p_n$ simply by scaling and translating the coarse mesh basis functions. For example, the basis functions for $V^p_n$ can be chosen as 
\[
\phi^n_{k,j}=2^{n/2}\phi_{k}\big(2^n(x+1)-2j-1\big),\;\;\;k=0,\hdots,p,\;\;\;j=0,\hdots,2^n-1.
\] 
We can then express our global modal approximation on a mesh of $N=2^n$ elements as
\begin{equation}
\label{eq:coarsefine}
u_h^{N}(x)=2^{-n/2}\sum_{j=0}^{2^n-1}\sum_{k=0}^p u_j^{k}\phi^n_{k,j}(x).
\end{equation}
The intermediary information which is required to advance from a coarse approximation space $V_{n}^p$ to a finer approximation space $V_{n+1}^p$ is contained within the wavelet space $W_{n}^p$ which is defined as the orthogonal compliment of $V_{n}^p$ in $V_{n+1}^p$:
\[
V_{n+1}^p=V_{n}^p\oplus W_{n}^p,\hspace{1cm} W_n^p\subset V^p_{n+1},\hspace{1cm} W^p_n\perp V_n^p.
\]
Inductively, this hierarchy allows the expression of the fine-mesh approximation space as the direct sum of the coarsest mesh approximation spaces and the intermediary wavelet spaces:
\[
     V_n^{p} = V_0^{p} \oplus W_0^{p} \oplus W_1^{p} \oplus \cdots \oplus W_{n-1}^{p}.
\]
Alternatively, this means our {fine-mesh} approximations are simply coarse mesh approximations plus details contained within these wavelet spaces. We denote the basis functions for $W^p_n$ by $\psi^n_{k,j},$ $k=0,\hdots,p$, $j=1,\hdots,2^n-1$. These functions are known as multi-wavelets. Decomposing our approximation in an analogous manner to our approximation spaces, {we have that} our {fine-mesh} approximation is just the sum of a coarse mesh approximation (on a mesh consisting of $N$ elements) and the details provided by the multi-wavelets:
\[
u_h^f(x)=u_h^{N}(x)+\sum^{2^{n-1}-1}_{j=0}\sum_{k=0}^p d^{n}_{k,j}\psi^{n}_{k,j}(x).
\]
An illustration of this idea in two-dimensions is given in Figure \ref{fig:2dmw}.  As a result, we see that one manner of obtaining fine-mesh approximations is to take a coarse-mesh approximation and add a linear combination of multi-wavelets, where the detail coefficient $d^n_{k,j}$ serve as the weights. {Note that the magnitude of the detail coefficients $d^n_{k,j}$ is determined by the local regularity of the function being approximated \cite{Wavelets2}.} The difficulty is that {when approximation information is only given} on a coarse grid, we do not have these finer details.  SIAC-MRA \cite{SIACMRA} enables a way of approximating detail coefficients $\hat{d}^{n}_{k,j}$ when only coarse mesh information is given.  Or, as an alternative for moving approximations from a coarse mesh to a fine mesh. 

%% file: Sections_LSIAC/Filter_HD.tex
\subsection{Multi-dimensional SIAC Filter}
\label{sec:SIACmultid}

The Smoothness-Increasing Accuracy-Conserving (SIAC)\linebreak post-processor is a convolution kernel originally designed to enable superconvergence in {finite element methods} \cite{CLSS2,RSA}.  They have since been extended for derivatives, boundaries and nonlinear hyperbolic equations {\cite{RC,XLiOne,SIACLinf,JiXuRyan2,SRV}} as well as mesh adaptivity \cite{Adaptivity}. The SIAC filter has been generalized for application to higher dimensional data using two methods {\cite{LSIAC,RSA}}: the first is a tensor product SIAC filter, whereby a convolution kernel is constructed by taking the tensor product of multiple one-dimensional SIAC kernels; the second is to construct a Line SIAC filter (L-SIAC), which is simply a one-dimensional SIAC kernel that has been rotated to align with a non-Cartesian axis.  Both of these techniques rely on understanding the one-dimensional kernel presented in this section.

 In one dimension, the SIAC filtered approximation is defined by 
\begin{equation}
    \label{eq:1Dconv}
    u_h^{\star}(x)=K_H\star u_h=\int_{-\infty}^{\infty}K^{(r+1,\ell)}_H(x-y)u_h(y)\;dy,
\end{equation}  
where the convolution kernel $K_H(x)=\frac{1}{H}K\big(\frac{x}{H}\big)$ is a scaled linear combination of $r+1$ central B-splines of order $\ell$, {that is normalized in $L^1$:}
\[K^{(r+1,\ell)}(t)=\sum_{\gamma=-r/2}^{r/2}c_{\gamma}B^{(\ell)}\Big(t+\gamma\Big). \]
The coefficients $c_{\gamma}$ are chosen so that the kernel reproduces polynomials {(consistency and moments)} up to degree $r$: 
\[K^{(r+1,\ell)}\star x^k=x^k \;\;\text{for}\;\; k=0,\hdots,r,\]
while the
central B-splines are piecewise polynomials defined recursively via the relation
\begin{align*}
   B^{(1)}&=\chi_{[-1/2,1/2)}\\
    B^{(n+1)}&=B^{(n)}\star B^{(1)}\\
    &=\frac{1}{n}\Big[(n/2+t)B^{(n)}(t+1/2)+(n/2-t)B^{(n)}(t-1/2)\Big].
\end{align*}
Central B-splines are used in the kernel construction as they provide finite support to the kernel, while allowing derivatives of the kernel to be expressed as divided differences. This later property and the {preservation of moments condition are crucial} in deriving {error estimates} \cite{CLSS2,SIACLinf,MengRyan1}. {An introduction to B-splines and their properties can be found in \cite{deBoor}.}
The kernel is scaled by a parameter $H$ typically set as the uniform element width $h$. To enable superconvergence in DG approximations{,} we require $r=2p$ and $\ell=p+1$. This allows for our filtered approximation to {obtain} $\mathcal{O}(h^{2p+1})$ errors in the $L^2-$norm and $L^{\infty}$-norm for linear hyperbolic equations.

{In two dimensions,} the Line SIAC filter {is} given by
\begin{equation}
    \label{eq:LSIAC}
u_h^{\star}(\mathbf{x})=\int_{\Gamma}K^{(r+1,\ell)}_{\Gamma,H}(t)u_h(t)\;dt, 
\end{equation}
where
\[\Gamma(t)={\mathbf{x}}+t(\cos(\theta),\sin(\theta)),\]
and
\begin{align}
\label{eq:kernel}
K^{(r+1,\ell)}_{\Gamma,H}(t)&=\sum_{\gamma=-r/2}^{r/2}c_{\gamma}B^{(\ell)}_H(t-\gamma)
    =\frac{1}{H}\sum_{\gamma=-r/2}^{r/2}c_{\gamma}B^{(\ell)}\Big(\frac{t}{H}-\gamma\Big).
    \intertext{The angle of rotation $\theta$ and the scaling parameter $H$ are selected to be}
    \theta&=\tan^{-1}\Bigg(\frac{h_{x_2}}{h_{x_1}}\Bigg),\quad
    H=h_{x_1}\cos(\theta)+h_{x_2}\sin(\theta).
    \end{align}
    In this article we choose the uniform mesh diagonal{,} $\theta=\pi/4$, and subsequently $H=\sqrt{2}h_{x_1}=\sqrt{2}h_{x_2}$.
    In higher dimensions we simply need to choose a new orientation and scaling factor. In three-dimensional space, the filtered approximation is given similarly where 
  the line $\Gamma$ parameterized by $t$ is given by
   \begin{equation}
  \Gamma(t)=\mathbf{x}+t\big(\cos(\theta)\sin(\phi),\sin(\theta)\sin(\phi),\cos(\phi))\big)=\mathbf{x}+t\mathbf{v},
   \end{equation}
  with $\theta,\phi$ being the orientation of the line in spherical coordinates. Because our approximation lies on a uniform three-dimensional mesh, the line filter will be aligned with one of the four diagonals connecting antipodal points of the cubic element. Hence, we choose $H=\sqrt{3}h$, where $h$ is the uniform mesh spacing, and $\theta=\pi/4$, $\phi=\tan^{-1}\big({\sqrt{2}}\big)$. This results in the direction vector $\mathbf{v}=\frac{1}{\sqrt{3}}(1,1,1)$.

     The utility of the L-SIAC filter is its ability to filter higher-dimensional data without increasing the dimensionality of the kernel support. This makes it an {alternative} to the tensor product filter when evaluating the convolution via quadrature sums. For example, in 2D for a quadratic approximation it can reduce the number of quadrature sums from 4000 to {ten \cite{LSIAC}}.  Additionally, the {constant} in the error estimate is reduced for the filtered solution. We consider the L-SIAC filter specifically because of its reduced and rotated kernel support in this context. The ability to manipulate the kernel orientation and dimensionality is anticipated to {be} useful in generalizing to triangulations due to the nature of the multi-variate basis functions.

%% file: Sections_LSIAC/SIAC_MRA_HD.tex
In this section, we introduce Line SIAC for more accurate multi-resolution analysis. This procedure allows us to obtain detail coefficients for moving information to successively finer grids. \textit{{These} detail coefficients are obtained for a uniform mesh by setting ${r}=2p$ and $\ell=1$ in the LSIAC filter} (Equation \eqref{eq:kernel}). {This} ensures enough moments in the underlying data are respected as well as ensuring {that} over-smoothing does not occur. For this choice of parameters, the filter coefficients are given by the equation
\begin{equation}
\sum_{\gamma = -p}^p\, c_\gamma \frac{1}{m+1}\left[\left(\frac{1}{2}-\gamma\right)^{m+1}+\left(\frac{1}{2}+\gamma\right)^{m+1}\right] = 
\begin{cases}
1, & m=0\\
0, & m=1,\dots,2p
\end{cases}.
\end{equation}

As presented in \cite{SIACMRA}, the one-dimensional convolution, \eqref{eq:1Dconv}, can be analytically calculated and thereby allows for expressing our filtered approximation in a new basis. As a result of the ordering of mesh and kernel breaks, we have a piecewise definition for the filtered approximation itself:
\[
u_h^{\star}(x)=\begin{cases}
\sum_{m=0}^{p+1}d_m^j\chi^m_L(\zeta),\hspace{1cm} &\zeta\in(-1,0)\\
\sum_{m=0}^{p+1}\tilde{d}_m^j\chi^m_R(\zeta),\hspace{1cm} &\zeta\in(0,1)
\end{cases},
\]
for each element. {Here} $\zeta=\xi_j(x)$ represents a local element mapping. The {filtered basis functions $\chi^m_L$ and $\chi^m_R$} used in this article and {their} associated coefficients are given in the supplementary materials.
{Assuming that our initial approximation is in $V_{n}^p$ and letting $N$ be the number of elements of this approximation space, denote the fine-mesh approximation by $u_h^f=\mathcal{P}_{n+1}u_h^{\star}$, where $\mathcal{P}_{n+1}$ represents the projection operator onto the approximation space $V^p_{n+1}$ consisting of $2N$ elements.} To obtain the fine-mesh modal coefficients in the ansatz $u_h^f(x)\big|_{I_n}=\sum_{k=0}^pu^n_k\phi^{{n}}_k(x),$ $n=1,\hdots, 2N$ we require
\[\langle u^f_{h},\phi^j_m\rangle_{I_{n}}=\langle u^{\star}_{h},\phi^j_m\rangle_{I_{n}},\hspace{1cm} m=0,\hdots,p.
\]

{Having} determined our fine-mesh approximation $u_h^f={\mathcal{P}_{n+1}}u_h^{\star}$, we can isolate the multi-wavelet \linebreak component simply by subtraction of the coarse-approximation in Equation \eqref{eq:coarsefine}:

\begin{equation}
{\sum^{2^{n-1}-1}_{j=0}\sum_{k=0}^p \hat{d}^{n}_{k,j}\psi^{n}_{k,j}(x)=\mathcal{P}_{n+1}u^{\star}_h-u_h,}
\end{equation}
{where $\psi^{n}_{k,j}$ are the basis functions of the multi-wavelet space.}
{These alternative multi-wavelet coefficients provide an improvement to our approximation, and the SIAC-MRA procedure provides a definition for them.} This allows for effectively transitioning the approximation onto a finer mesh. Additionally, it only requires coarse mesh data to be constructed. This differs significantly from traditional MRA where multi-wavelet coefficients are intrinsically dependent upon {the} initial function $u$, {which means that} we cannot construct multi-wavelet coefficients for improving our approximation if only coarse modal data is available. 


{Knowing} the form of the filtered approximation, we can map {the} coefficients to the fine-mesh modes simply by application of a projection operator. Furthermore, knowing the forms of the filtered approximation allows us to discretize this procedure and construct projection {matrices} to perform this scale transition. This will be discussed further in the implementation Section \ref{sec:Implementation}.

As demonstrated in \cite{SIACMRA}, repeated application of the filtering-projection procedure results in increasingly better approximations. It should be noted that SIAC-MRA is mesh-dependent in that our ability to simplify the convolution is intrinsically linked to patterns in the mesh construction. In this paper we consider uniform quadrilateral and hexahedral meshes, to generalize this procedure to higher dimensional spaces and introduce LSIAC-MRA.

In a manner analogous to the one-dimensional case \cite{SIACMRA}, we can take advantage of the simplicity of the kernel's B-spline components and the predictability of the element boundaries to analytically evaluate the convolutions given by Equation \eqref{eq:LSIAC}. Following this procedure, and letting $\zeta_x=\xi_i(x)$, $\zeta_y=\xi_j(y)$, and $\zeta_z=\xi_k(z)$ denote the mapping of our filtering point to the reference element, we obtain the following expressions for our filtered approximations:

%
%

For the two-dimensional Line SIAC, the filtered approximation can be written as
\begin{equation}
\label{eq:2DLSIAC}
u_h^\star(x,y)=\sum_{P=\{L,M,R\}}\sum_{m=0}^{p}\sum_{r=0}^{p} a^{S,P,i,j}_{m,r}\chi^{m,r}_{S,P}(\zeta_x,\zeta_y),
 \end{equation}
where $(x,y)$ belongs to region $S$ of element $(i,j)$. These regions are delineated in { Figure \ref{fig:L-SIAC_sppt}}. {The $P$ index results from splitting up the convolution into ``left",``middle", and ``right" components owing to element breaks.} 
Analogously, the three-dimensional LSIAC filtered approximation can be written as
\begin{equation}
\label{eq:3DLSIAC}
u_h^\star(x,y,z)=\sum_{P=\{L,LM,RM,R\}}\sum_{m=0}^{p}\sum_{r=0}^{p}\sum_{\ell=0}^{p} a^{S,P,i,j,k}_{m,r,\ell}\chi^{m,r,\ell}_{S,P}(\zeta_x,\zeta_y,\zeta_z), \end{equation}
where $(x,y,z)$ belongs to region $S$ of element $(i,j,k)${, where again the $P$ index denotes an ordering with respect to element breaks}. A description of the aforementioned regions and an associated classification scheme are provided in the supplementary materials.

  In two dimensions, {the} LSIAC implementation requires 6 quadrature domains per element. In three {dimensions} there are 24 quadrature domains per element. Note that for the L-SIAC filters, more complicated integrals of shifted Legendre Polynomials are produced. Though the reduced dimension of the L-SIAC kernel support means less information is needed to compute the $a$ coefficients, the integrand is no longer separable which leads to more complicated functions in the expansion. All the $\chi$ functions and their associated coefficients are detailed in the supplementary materials. Much the same as in one dimension, we can project these filtered approximations onto the finer mesh, and thereby construct fine-mesh approximations. The only difference is that we are now projecting onto either the quadrant or octants of reference element in two and three dimensions respectively. In the L-SIAC case this requires us to account for the piecewise definition of the filtered approximations and split up projections according to the associated domains of definition.

 {In what follows,} we demonstrate the performance of {LSIAC-MRA} . Because $\psi^{(1)}(x)=\chi_{[-1/2,1/2)}(x)$,  we {are} able to analytically evaluate the convolution for the L-SIAC kernel by easily accounting for mesh and kernel breaks. This {enables} us to obtain a closed form expression for the filtered approximation expressed in a new basis {and allows for our improved multi-resolution analysis technique.}

%% file: Sections_LSIAC/Implementation.tex
\section{Implementation} \label{sec:Implementation}
{In this section we describe a brief overview of the implementation of 2D LSIAC-MRA as the same ideas extend to higher dimensions. Additional details on the filtered approximation basis functions and a computational algorithm are included in the supplementary materials.}
\subsection{Transition Operators}
\input{Pictures/Region_Figure}

{Given {coarse}-mesh modal information, we wish to construct transition matrices for mapping the coarse-mesh modal coefficients to their refined fine-mesh counterparts. In our discrete setting these transition matrices are the composition of the discrete filtering and projection matrices $T={\mathcal{P}_{n+1}}K_H$. In two dimensions we will have a separate transition matrix for each quadrant. It is important to note that the sequencing of mesh and kernel breaks, and the elements within the kernel's support varies from region to region (see {Figure \ref{fig:L-SIAC_sppt}}). As can be seen when computing the refined approximation, we must project our filtered approximation onto the finer mesh. In the {2D} L-SIAC case we must sum the projections over the two triangular regions composing quadrants I and III. This requires different transition matrix construction from quadrant-to-quadrant.} {In the following we outline the construction procedure and refer the reader to Figure \ref{fig:matrix_Description} for a depiction of the effects of each matrix in the composition for the case $p=0$.} 

\input{Sections_LSIAC/2D_L-SIAC_transition_matrices}

%% file: Pictures/Region_Figure.tex
 \begin{figure}[hbt!]
    \centering

  \scalebox{0.9}{    
\begin{tabular}{c c c }

\multicolumn{3}{c}{
\begin{tikzpicture}

    \begin{customlegend}[legend columns=-1,legend entries={\phantom{aa}Element\phantom{aaaaaaa},\phantom{aa}Filter Region\phantom{aaaaaaa}, \phantom{aa}Kernel Support}]
    \addlegendimage{black,fill=gray!20,area legend} 
    \addlegendimage{ color =orange ,pattern=north west lines, pattern color=orange,area legend}
    \addlegendimage{gray,fill=gray,dashed,ultra thick,sharp plot}
    legend style={at={(0.5,0.5)},anchor=center}
    \end{customlegend}
\end{tikzpicture}
}

\\

\resizebox{0.23\linewidth}{!}{
\begin{tikzpicture}[every node/.style={minimum size=.5cm-\pgflinewidth, outer sep=0pt}]
        \fill[color = gray!20] (0.505,0.505) rectangle (0.995,0.995);
\draw[step=0.5cm,color=black] (0,0) grid (1.5,1.5);
    
    \node at (0.75,1.75){\tiny Q1 T};

    \draw[ color =orange ,pattern=north west lines, pattern color=orange] (0.75, 0.75) -- (1.0, 1.0)--(0.75,1.0)--(0.75,0.75);
    
    \draw[dashed, thin] (0.75,0.5)--(0.75,1.0);
\draw[dashed, thin] (0.5,0.75)--(1.0,0.75);
     \draw[dashed, thin] (0.5,0.5)--(1.0,1.0);
     
\draw[dashed,thick, color=gray] (0.5,0.5)--(1.25,1.25)--(1.0,1.25)--(0.5,0.75)--(0.5,0.5);
\end{tikzpicture}
}

&

\resizebox{0.23\linewidth}{!}{
\begin{tikzpicture}[every node/.style={minimum size=.5cm-\pgflinewidth, outer sep=0pt}]
       \fill[color = gray!20] (0.505,0.505) rectangle (0.995,0.995);
 \draw[step=0.5cm,color=black] (0,0) grid (1.5,1.5);
        \node at (0.75,1.75){\tiny Q2 T};

    \draw[ color =orange ,pattern=north west lines, pattern color=orange] (0.5, 0.75) rectangle (0.75, 1.0);
    
    \draw[dashed, thin] (0.75,0.5)--(0.75,1.0);
\draw[dashed, thin] (0.5,0.75)--(1.0,0.75);
     \draw[dashed, thin] (0.5,0.5)--(1.0,1.0);
     
\draw[dashed,thick, color=gray] (0.25,0.5)--(0.5,0.5)--(1.0,1.0)--(1.0,1.25)--(0.75,1.25)--(0.25,0.75)--(0.25,0.5);
\end{tikzpicture}
}
&

\resizebox{0.23\linewidth}{!}{
\begin{tikzpicture}[every node/.style={minimum size=.5cm-\pgflinewidth, outer sep=0pt}]
       \fill[color = gray!20] (0.505,0.505) rectangle (0.995,0.995);
 \draw[step=0.5cm,color=black] (0,0) grid (1.5,1.5);
        \node at (0.75,1.75){\tiny Q3 T};

    \draw[ color =orange ,pattern=north west lines, pattern color=orange] (0.5, 0.5) -- (0.75, 0.75)--(0.5,0.75)--(0.5,0.5);
    
    \draw[dashed, thin] (0.75,0.5)--(0.75,1.0);
\draw[dashed, thin] (0.5,0.75)--(1.0,0.75);
     \draw[dashed, thin] (0.5,0.5)--(1.0,1.0);
     
\draw[dashed,thick, color=gray] (0.25,0.25)--(1.0,1.0)--(0.75,1.0)--(0.25,0.5)--(0.25,0.25);
\end{tikzpicture}
}

\\

\resizebox{0.23\linewidth}{!}{
\begin{tikzpicture}[every node/.style={minimum size=.5cm-\pgflinewidth, outer sep=0pt}]
    \fill[color = gray!20] (0.505,0.505) rectangle (0.995,0.995);

    \draw[step=0.5cm,color=black] (0,0) grid (1.5,1.5);
        \node at (0.75,1.75){\tiny Q1 B};

    \draw[ color =orange ,pattern=north west lines, pattern color=orange] (0.75, 0.75) -- (1.0,0.75) --(1.0,1.0)--(0.75,0.75);
    
    \draw[dashed, thin] (0.75,0.5)--(0.75,1.0);
\draw[dashed, thin] (0.5,0.75)--(1.0,0.75);
     \draw[dashed, thin] (0.5,0.5)--(1.0,1.0);
     
\draw[dashed,thick, color=gray] (0.5,0.5)--(0.75,0.5)--(1.25,1.0)--(1.25,1.25)--(0.5,0.5);
\end{tikzpicture}
}
&
\resizebox{0.23\linewidth}{!}{
\begin{tikzpicture}[every node/.style={minimum size=.5cm-\pgflinewidth, outer sep=0pt}]
       \fill[color = gray!20] (0.505,0.505) rectangle (0.995,0.995);
 \draw[step=0.5cm,color=black] (0,0) grid (1.5,1.5);
    
    \draw[ color =orange ,pattern=north west lines, pattern color=orange] (0.75, 0.5) rectangle (1.0,0.75);
        \node at (0.75,1.75){\tiny Q4 B};

    \draw[dashed, thin] (0.75,0.5)--(0.75,1.0);
\draw[dashed, thin] (0.5,0.75)--(1.0,0.75);
     \draw[dashed, thin] (0.5,0.5)--(1.0,1.0);
     
\draw[dashed,thick, color=gray] (0.5,0.25)--(0.75,0.25)--(1.25,0.75)--(1.25,1.0)--(01.0,1.0)--(0.5,0.5)--(0.5,0.25);
\end{tikzpicture}
}
&

\resizebox{0.23\linewidth}{!}{
\begin{tikzpicture}[every node/.style={minimum size=.5cm-\pgflinewidth, outer sep=0pt}]
      \fill[color = gray!20] (0.505,0.505) rectangle (0.995,0.995);
  \draw[step=0.5cm,color=black] (0,0) grid (1.5,1.5);
        \node at (0.75,1.75){\tiny Q3 B};

    \draw[ color =orange ,pattern=north west lines, pattern color=orange] (0.5, 0.5) -- (0.75,0.5) --(0.75,0.75)--(0.5,0.5);
    
    \draw[dashed, thin] (0.75,0.5)--(0.75,1.0);
\draw[dashed, thin] (0.5,0.75)--(1.0,0.75);
     \draw[dashed, thin] (0.5,0.5)--(1.0,1.0);
     
\draw[dashed,thick, color=gray] (0.25,0.25)--(0.5,0.25)--(1.0,0.75)--(1.0,1.0)--(0.25,0.25);
\end{tikzpicture}
}

\end{tabular}
}
         \caption{{2D depiction of the varying element overlap of the L-SIAC kernel resulting in varying filtered approximation forms for $p=0$. The element is shaded in gray, filter region has diagonal lines, and the kernel support is encapsulated by dashed lines.  The filter region represents the location of possible filtering points, and the kernel support represents the overlap of the kernel support for those filtering points.}}
    \label{fig:L-SIAC_sppt}
\end{figure}
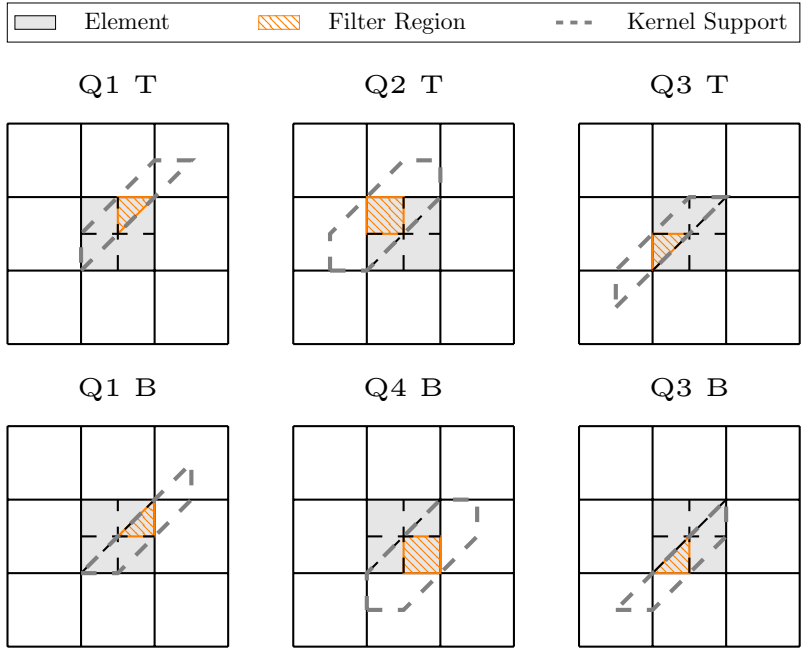

%% file: Sections_LSIAC/2D_L-SIAC_transition_matrices.tex
{For the purposes of performing the enhancement procedure via matrix operations it is important to have a standard ordering for the modal information. As constructing the transition matrices is easier if modal coefficients are grouped by element, we denote} the $(N^2(p+1))^2\times 1$ vector of modal coefficients by
\[\Vec{u}=\begin{bmatrix}\vec{u}^1 &|& \hdots &|&\vec{u}^N\end{bmatrix}^T, \]
where
\[{\vec{u}^j}=\begin{bmatrix}\vec{u}^{1,j} \hdots\vec{u}^{N,j} \end{bmatrix}, \]
\[{\vec{u}^{i,j}}=\begin{bmatrix} \vec{u}^{i,j}_{0} &|& \hdots &|& \vec{u}^{i,j}_{p}  \end{bmatrix}, \]
and
\[{\vec{u}^{i,j}_{ky}}=\begin{bmatrix} u^{i,j}_{0,ky} &|& \hdots &|& u^{i,j}_{p,ky}  \end{bmatrix}. \]

\subsubsection{Sifting Matrices}
All of {the coefficients in \eqref{eq:2DLSIAC}} fit the form {\begin{equation}
    \label{eq:coeff}
a^{i,j}_{m,r}(q_x,q_y)=\sum_{\gamma=-p}^{p}\, c_\gamma u_{kx,ky}^{i+\gamma+q_x,j+\gamma+q_y},
\end{equation}}{where $q_x,q_y\in\{-1,0,1\}$ are uniquely determined by the region $S$ and position $P$ as detailed in the supplementary materials}. To compute these coefficients, {a} selection matrix, denoted $\mathscr{S}^{i,j}$, is first constructed. Assuming periodic boundary conditions, this matrix will select all {modal values from the $(2p+3)^2$ elements about element $(i,j)$. We will later further pair this information to only those elements under the line filter's support.} Define the $(2p+3)N(p+1)^2\times N^2(p+1)^2$ matrix $S^j$ by
\[S^j=\begin{bmatrix}\ddots& \vdots& \iddots\\
\hdots&\omega_j(m,r)\mathbb{I}_{N(p+1)^2\times N(p+1)^2}&\hdots \\
\iddots &\vdots& \ddots 
\end{bmatrix} \]
where
\[\omega_j(m,r)=\begin{cases} 1,\;\; j\neq N\;\&\; p+2-m+r=j\;\mod(N)\\
 1,\;\; j= N\;\& \;p+2-m+r=j \text{ or }0\\
 0,\;\; \text{ else }
 \end{cases},\]
 for $m=1:(2p+3)$ and $r=1:N$.
 This matrix selects the relevant $j$-coordinates:
\[S^{j}\vec{u}=[u^{j-(p+1)}\;\hdots\;u^{j+(p+1)}]^T. \]

Next, construct an $(2p+3)(p+1)^2 \times N(p+1)^2$ matrix $\tilde{S}^i$ to select the relevant $i-$coordinates from within each vector.
Define
\[\tilde{S}^i=\begin{bmatrix}\ddots& \vdots& \iddots\\
\hdots&{\omega}_i(m,r)\mathbb{I}_{(p+1)^2\times (p+1)^2}&\hdots \\
\iddots &\vdots& \ddots 
\end{bmatrix} \]
 for $m=1:(2p+3)$ and $r=1:N$.
 We have that
\[\tilde{S}^iu^j=[u^{i-(p+1),j} \;\hdots\; u^{i+(p+1),j}]^T. \]
 Hence, defining $\mathscr{S}^{i,j}$ to be the $(2p+3)^2(p+1)^2\times N^2(p+1)^2$ matrix given by
 \[\mathscr{S}^{i,j}=\begin{bmatrix}
 \tilde{S}^i& &0\\
 & \ddots & \\
 0& & \tilde{S}^i
 \end{bmatrix}S^j,\]
 where there are $(2p+3)$ $\tilde{S}^i$ matrices on the diagonal, we have
\[\mathscr{S}^{i,j}u=\Big[[u^{i-(p+1),j-(p+1)}\;\dots u^{i+(p+1),j-(p+1)}]\;\; \hdots\;\;[u^{i-(p+1),j+(p+1)}\;\dots u^{i+(p+1),j+(p+1)}] \Big]^T=\tilde{u}^{i,j}\]
which, {as stated previously, only} contains modal information from the $(2p+3)^2$ elements around element $(i,j)$.

{Now we construct an intermediary matrix ${\mathscr{C}^{q_x,q_y}}$ whose function is to select the elements from $\tilde{u}^{i,j}$ which fall within the support of the line integral as affected by the shifting arguments $(q_x,q_y)$. Define the $(2p+1)(2p+3)(p+1)^2\times(2p+3)^2(p+1)^2$ matrix $C^{q_y}$ by}
\[{C^{q_y}=\begin{bmatrix}
\delta_{-1,q_y}\mathbb{I}&\delta_{0,q_y}\mathbb{I}&\delta_{1,q_y}\mathbb{I}&  &0\\
&\ddots &\ddots &\ddots& &\\
0&&\delta_{-1,q_y}\mathbb{I}&\delta_{0,q_y}\mathbb{I}&\delta_{1,q_y}\mathbb{I}
\end{bmatrix},}\]
{where $\delta$ is the Kronecker delta and $\mathbb{I}$ is the $(2p+3)(p+1)^2\times(2p+3)(p+1)^2$ identity matrix. This matrix has the effect of further restricting the relevant $j-$indices based off the shifting argument $q_y$. Then define the $(2p+1)(p+1)^2\times (2p+3)(2p+1)(p+1)^2$ matrix $\mathscr{C}^{q_x}$ by}
\[{\mathscr{C}^{q_x}=\begin{bmatrix}
\delta_{-1,q_x}\mathbb{I}&\delta_{0,q_x}\mathbb{I}&\delta_{1,q_x}\mathbb{I}&  &0\\
& &R_2 & &\\
& &\vdots & &\\
& &R_{2p+1} & &
\end{bmatrix},}  \]
{where $\mathbb{I}$ is $(p+1)^2\times(p+1)^2$, and $R_j$ is the first block row defined above with the non-zero entries translated $(j-1)(2p+4)(p+1)^2$ places to the right. We then have the desired $(2p+1)(p+1)^2\times(2p+3)^2(p+1)^2$ intermediary matrix $\mathscr{C}^{q_x,q_y}$ given by the composition}
\[{\mathscr{C}^{q_x,q_y}=\mathscr{C}^{q_x}C^{q_y}. }\]
{Applying this matrix produces}:
\[\mathscr{C}^{q_x,q_y}\tilde{u}^{i,j}=\tilde{u}^{i,j}_{q_x,q_y}=[u^{i-(p-q_x),j-(p-q_y)},\hdots,u^{i+(p+q_x),j+(p+q_y)} ]^T{,}\]
{which are the modal coefficients appearing in \eqref{eq:coeff}.}

\input{Pictures/matrix_figure/total}

\subsubsection{{Coefficient Matrices}}

{We now describe the construction of a matrix for computing \eqref{eq:coeff} from the sifted modal coefficients. Define} the matrix $A$ to be the $(p+1)^2\times(2p+1)(p+1)^2$ matrix given by
\[A=\frac{1}{2}\begin{bmatrix}
c_{-p}\mathbb{I}&&&0&|&\hdots&|&c_{p}\mathbb{I}&&&0 \\
0&c_{-p}\mathbb{I}&&0&|&\hdots&|&0&c_{p}\mathbb{I}&&0&\\
&&\ddots&&|&\hdots&|&&\ddots&&\\
0&&&c_{-p}\mathbb{I}&|&\hdots&|&0&&&c_{p}\mathbb{I}&
\end{bmatrix}\]
where the identity matrices are ${(p+1)\times(p+1)}$. {This is the matrix used to perform the summation and leads to}
\[A\tilde{u}^{i,j}_{q_x,q_y}=\vec{a}^{i,j}(q_x,q_y)\]
where 
\[\vec{a}^{i,j}(q_x,q_y)=[a^{i,j}_{r=0}\;\hdots\; a^{i,j}_{r=p}]^T(q_x,q_y)\]
with 
\[{a^{i,j}_{r}}(q_x,q_y)=[a^{i,j}_{0,r}(q_x,q_y)\;\hdots\; a^{i,j}_{p,r}(q_x,q_y)]. \]
Hence, the seven coefficient matrices $K^{i,j}_{q_x,q_y}$ are defined by the composition
\[ K^{i,j}_{q_x,q_y}=A\mathscr{C}^{q_x,q_y}\mathscr{S}^{i,j}\]
and {give}
\[K^{i,j}_{q_x,q_y}u=\vec{a}^{i,j}(q_x,q_y). \]
Rather than {concatenating} these matrices at this point in the procedure, we will wait until after applying the projection {in order} to minimize the number of transformations performed. {It is important to note that each choice of $(q_x,q_y)$ computes a different set of coefficients pertaining to a different set of basis functions, each of which is required for obtaining the enhanced approximation.}

\subsubsection{Projection Matrices}

{We now discuss the procedure necessary to} perform the projections. Because of the different filtered approximation forms in each quadrant, we must perform each projection slightly differently. Quadrants II and IV are similar as are I and III. Here we detail a procedure for QII.

\section*{QII}
In QII we define the projection of this approximation by $u^f_h|_{K^{2i-1,{2j}}}(x,y)$ where the new modal coefficients are computed by
\begin{align}
    \langle u_h^{\star},\phi_{kx}\phi_{ky}\rangle_{QII}&=\langle u_h^f,\phi_{kx}\phi_{ky}\rangle_{QII}\\
    &=u_{kx,ky}^{2i-1,{2j}},\;\;\;\;kx,ky=0,\hdots,p,\;\;i,j=1,\hdots,N.
\end{align}
To account for the change in domain scaling define $\xi_x(\zeta_x)=\frac{\zeta_x-1}{2}$ and $\xi_y(\zeta_y)=\frac{\zeta_y+1}{2}$. We have
\begin{align*}
    u_{kx,ky}^{2i-1,{2j}}=\sum_{m=0}^p\sum_{r=0}^p\int_{-1}^1\int_{-1}^1\Big\{&a^{i,j}_{m,r}(-1,0)\chi^{m,r}_{2T,L}(\xi_x,\xi_y)\\
    +&a^{i,j}_{m,r}(0,0)\chi^{m,r}_{2T,M}(\xi_x,\xi_y)\\
    +&a^{i,j}_{m,r}(0,1)\chi^{m,r}_{2T,R}(\xi_x,\xi_y)\Big\}\phi_{kx}(\zeta_x)\phi_{ky}(\zeta_y)\;d\zeta_y\;d\zeta_x.
\end{align*}
This expression will be broken up into three separate integrals, one for each of the different $\chi_{(\_)}$ functions where $(\_)\in\{L,M,R\}$. For each $\chi_{(\_)}$, define \[P_{(\_)}(kx,ky,m,r)=\int_{-1}^1\int_{-1}^1\chi^{m,r}_{(2T,\_)}(\xi_x,\xi_y)\phi_{kx}(\zeta_x)\phi_{ky}(\zeta_y)\;d\zeta_y\;d\zeta_x.\]
Define the projection matrix $P^{2T}_{(\_)}$ by
\[P^{2T}_{(\_)}={\begin{bmatrix}\vec{P}(0,0)&
\hdots&
\vec{P}(p,0)&
&
\hdots&
&
\vec{P}(0,p)&
\hdots&
\vec{P}(p,p)\end{bmatrix}}^{{T}}, \]
where 
\[\vec{P}(kx,ky)=[P(kx,ky,0,0)\hdots P(kx,ky,p,0) \hdots\hdots P(kx,ky,0,p)\hdots P(kx,ky,p,p)]^{{T}}. \]
Hence, we have 
\[u^f_{2i-1,{2j}}=P_L^{2T}\vec{a}^{i,j}(-1,0)+P_M^{2T}\vec{a}^{i,j}(0,0)+P_R^{2T}\vec{a}^{i,j}(0,1).\]
Noting the necessary choice of $K^{i,j}_{q_x,q_y}$ to obtain each of the coefficient vectors, we have 
\[u^f_{2i-1,{2j}}=(P_L^{2T}K^{i,j}_{-1,0}+P_M^{2T}K^{i,j}_{0,0}+P_R^{2T}K^{i,j}_{0,1})\vec{u}.\]
Denoting this composition by $T^{{2}}_{i,j}$ we can concatenate these matrices together to obtain the global transition matrix $T^{{2}}$ defined by
\[T^{{ 2}}=\begin{bmatrix}
T^{{2}}_{1,1}&\hdots& T^{{2}}_{N,1} &|\hdots|&T^{{2}}_{1,N}&\hdots T^{{2}}_{N,N} 
\end{bmatrix}^T,\]
such that 
\[u^f_{{2}}=T^{{2}}u. \]

\section*{QIV}
The projection matrices for quadrant IV are the same except that the quadrant IV basis function are used and the change of variables are instead $\xi_x(\zeta_x)=\frac{\zeta_x+1}{2}$ and $\xi_y(\zeta_y)=\frac{\zeta_y-1}{2}$.

\section*{QIII}
The cases of quadrants I and III are similar in that a change of variables must be applied to split the integrals over the quadrant into a pair of {integrals} over diagonals of that quadrant. Beginning with the integral over the whole quadrant, we have
\begin{align*}
    \bar{u}_{kx,ky}^{2i-1,2j-1}&= \langle u^{\star}_h,\phi_{kx}^i\phi_{ky}^j \rangle_{QIII}\\
  &=\int_{-1}^1\int_{-1}^1u_h^{\star}|_{QIII}(\xi_x(\zeta_x),\xi_y(\zeta_y))\phi^i_{kx}(\zeta_x)\phi^j_{ky}(\zeta_y)\;d\zeta_y\;d\zeta_x,
\end{align*}
where $\xi_x=\frac{\zeta_x-1}{2}$ and $\xi_y=\frac{\zeta_y-1}{2}$. Breaking this up into {integrals} over $3T$ and $3B$ we have 
\begin{align*}
     \bar{u}_{kx,ky}^{2i-1,2j-1} &=\int_{-1}^1\int_{\zeta_x}^1u^{\star}_{3T}(\xi_x,\xi_y)\phi^i_{kx}(\zeta_x)\phi^j_{ky}(\zeta_y)\;d\zeta_y\;d\zeta_x\\
     &+\int_{-1}^1\int^{\zeta_x}_{-1}u^{\star}_{3B}(\xi_x,\xi_y)\phi^i_{kx}(\zeta_x)\phi^j_{ky}(\zeta_y)\;d\zeta_y\;d\zeta_x {.}
\end{align*}
Now to use Gauss-Legendre quadrature{,} we must scale the bounds of the inner integral to $[-1,1]$. To do so introduce the change of variables $\zeta_y=\frac{1-\zeta_x}{2}\alpha_T+\frac{1+\zeta_x}{2}$ for the first integral and  $\zeta_y=\frac{\zeta_x-1}{2}\alpha_B+\frac{1+\zeta_x}{2}$ for the second. This allows us to rewrite the expression as
\begin{align*}
       \bar{u}_{kx,ky}^{2i-1,2j-1}&= \int_{-1}^1\int_{-1}^1u^{\star}_{3T}(\xi_x,\xi_y(\zeta_y(\alpha_T))\phi^i_{kx}(\zeta_x)\phi^j_{ky}\Big(\frac{1}{2}\big[\alpha_T(1-\zeta_x)+1+\zeta_x\big]\Big)\frac{1-\zeta_x}{2}\;d\alpha_T\;d\zeta_x\\
     &+\int_{-1}^1\int^{1}_{-1}u^{\star}_{3B}(\xi_x,\xi_y(\zeta_y(\alpha_B))\phi^i_{kx}(\zeta_x)\phi^j_{ky}\Big(\frac{1}{2}\big[\alpha_B(\zeta_x-1)+1+\zeta_x\big]\Big)\frac{\zeta_x-1}{2}\;d\alpha_B\;d\zeta_x.\\
     &=\mathscr{I}^T+\mathscr{I}^B.
\end{align*}
We now need to split the $\mathscr{I}$ terms into a sum of integrals over the basis functions and then multiply by the appropriate modal coefficients.
Defining \[P^{3T}_{(\_)}(kx,ky,m,r)=\int_{-1}^1\int_{-1}^1\chi^{m,r}_{3T,(\_)}(\xi_x,\xi_y(\alpha_T))\phi_{kx}(\zeta_x)\phi_{ky}(\zeta_y(\alpha_T))\;d\zeta_x\;d\alpha_T,\]
where $(\_)\in\{L,M,R\}$, we construct the projection matrices $P^{3T}_L,P^{3T}_M,P^{3T}_R$ in the same manner as in quadrant II. Now for the second set of integrals define \[P^{3B}_{(\_)}(kx,ky,m,r)=\int_{-1}^1\int_{-1}^1\chi^{kx,ky}_{3B,(\_)}(\xi_x,\xi_y(\alpha_B))\phi_{kx}(\zeta_x)\phi_{ky}(\zeta_y(\alpha_B))\;d\zeta_x\;d\alpha_B,\]
and analogously construct the projection matrices $P^{3B}_L,P^{3B}_M,P^{3B}_R$. We then have
\begin{align*}
u^f_{2i-1,2j-1}&=\phantom{+}P_L^{3T}\vec{a}^{i,j}(-1,-1)+P_M^{3T}\vec{a}^{i,j}(-1,0)+P_R^{3T}\vec{a}^{i,j}(0,0)\\
&\phantom{=}+P_L^{3B}\vec{a}^{i,j}(-1,-1)+P_M^{3B}\vec{a}^{i,j}(0,-1)+P_R^{3B}\vec{a}^{i,j}(0,0).
\end{align*}
Now $\vec{a}^{i,j}(q_x,q_y)=K^{i,j}_{q_x,q_y}\vec{u},$ thus 
\[u^f_{2i-1,2j-1}=\Big[(P_L^{3T}+P_L^{3B})K_{-1,-1}^{i,j}+(P^{3T}_R+P^{3B}_R)K_{0,0}^{i,j}+P_M^{3T}K^{i,j}_{-1,0}+P_M^{3B}K^{i,j}_{0,-1}\Big]\vec{u}.\]
Denoting this composition by $T^{3}_{i,j}$ we can concatenate these matrices together to obtain the global transition matrix $T^{3}$ defined by
$$T^{3}=\begin{bmatrix}
T^{3}_{1,1}&\hdots& T^{3}_{N,1} &|\hdots|&T^{3}_{1,N}&\hdots T^{3}_{N,N} 
\end{bmatrix}^T,$$
such that 
$$u^f_{3}=T^{3}u. $$

\section*{QI}
The case of quadrant one is that same as quadrant three except now $\xi_x=\frac{\zeta_x+1}{2}$, $\xi_y=\frac{\zeta_y+1}{2}$, and the $\chi$ functions used are those defined on quadrant one.

\section*{Coding Comment:}
{The LSIAC-MRA procedure was implemented in MATLAB$^{TM}$ R2021a. The authors make no claim to the computational efficiency of the implementation but simply state that the memory requirements are on the order of the number of modal values of the finer mesh $4N^2(p+1)^2$. An algorithm, {whereby only the local transition matrices $T^Q_{i,j}$ are formed to perform the enhancement} is provided in the supplementary materials. }

%% file: Pictures/matrix_figure/total.tex
  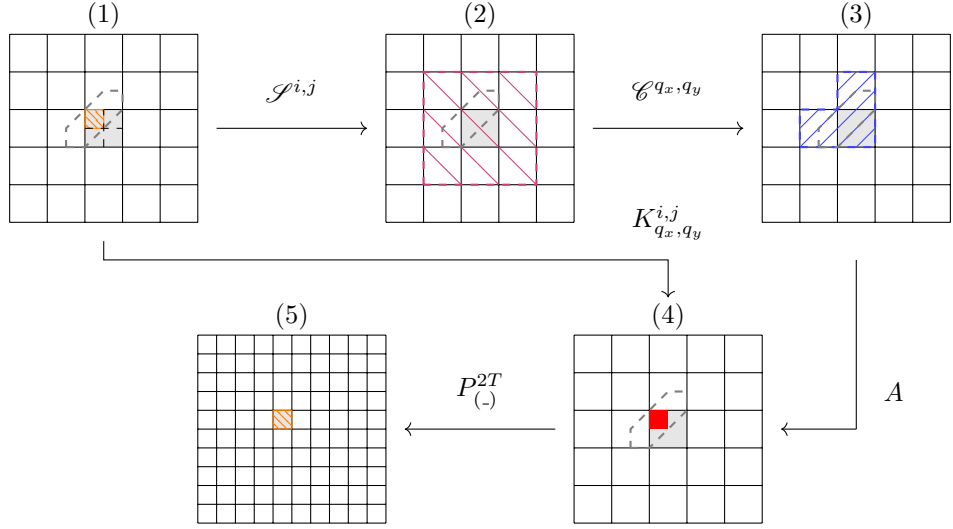
\begin{figure}[btp!]

    \begin{center}
\begin{tikzpicture}[every node/.style={minimum size=.5cm-\pgflinewidth, outer sep=0pt}]

    \begin{scope}
        \input{Pictures/matrix_figure/s1}
    \end{scope}

\draw [->] (2.75, 1.25) to (4.75, 1.25);
         \node(circle) at (3.75, 1.75)  (a) {$\mathscr{S}^{i,j}$};
\node(circle) at (1.25, 2.75)  (a) {$(1)$};
    \begin{scope}[xshift=5cm]
        \input{Pictures/matrix_figure/s2}
    \end{scope}
    
\draw [->] (7.75, 1.25) to (9.75, 1.25);
         \node(circle) at (8.75, 1.75)  (a) {$\mathscr{C}^{q_x,q_y}$};
\node(circle) at (6.25, 2.75)  (a) {$(2)$};

    \begin{scope}[xshift=10cm]
        \input{Pictures/matrix_figure/s3}
    \end{scope}
    

\draw [-] (11.25, -0.5) to (11.25, -2.75);
\draw [->] (11.25, -2.75) to (10.25, -2.75);
   \node(circle) at (11.75, -2.25)  (a) {$A$};
   \node(circle) at (11.25, 2.75)  (a) {$(3)$};

    \begin{scope}[xshift=7.5cm,yshift=-4cm]
        \input{Pictures/matrix_figure/s4}
    \end{scope}

\draw [->] (7.25, -2.75) to (5.25, -2.75);
   \node(circle) at (6.25, -2.25)  (a) {$P^{2T}_{(\_)}$};
   \node(circle) at (8.75, -1.25)  (a) {$(4)$};

    \begin{scope}[xshift=2.5cm,yshift=-4cm]
        \input{Pictures/matrix_figure/s5}
    \end{scope}
   \node(circle) at (3.75, -1.25)  (a) {$(5)$};

   \draw[->] (1.25,-0.25) to (1.25,-0.5) to (8.75,-0.5) to (8.75,-1.0); 
      \node(circle) at (8.75, 0.0)  (a) {$K^{i,j}_{q_x,q_y}$};

\end{tikzpicture}

    \end{center}   
      \caption{{Purpose of transition matrix components during refinement procedure: $(1)\rightarrow (2)$ Select modes from elements relevant to refining element $(i,j)$, $(2)\rightarrow (3)$ Select only the modes relevant to filtering a given region, in this case $2T$, $(3)\rightarrow (4)$ Use that modal information to obtain filtered approximation coefficients, $(4)\rightarrow(5)$ Obtain fine mesh modes by projecting the filtered approximation onto the refined mesh.  Notice that $K^{i,j}_{q_x,q_y}=A\mathscr{C}^{q_x,q_y}\mathscr{S}^{i,j}.$}}
         \label{fig:matrix_Description}
 \end{figure}

%% file: Pictures/matrix_figure/s1.tex

    \fill[color = gray!20] (1.005,1.005) rectangle (1.505,1.505);

    \draw[step=0.5cm,color=black] (0,0) grid (2.5,2.5);

    \draw[ color =orange ,pattern=north west lines, pattern color=orange] (1, 1.25) rectangle (1.25, 1.5);

    \draw[dashed, thin] (1.25,1)--(1.25,1.5);
\draw[dashed, thin] (1,1.25)--(1.5,1.25);

\draw[dashed,thick, color=gray] (0.75,1.0)--(1.0,1.0)--(1.5,1.5)--(1.5,1.75)--(1.25,1.75)--(0.75,1.25)--(0.75,1.0);

%% file: Pictures/matrix_figure/s2.tex
    \fill[color = gray!20] (1.005,1.005) rectangle (1.505,1.505);

    \draw[step=0.5cm,color=black] (0,0) grid (2.5,2.5);

    
        
     \draw[ color= purple!75] (0.5,1.0)--(1,0.5);
     \draw[ color= purple!75] (0.5,1.5)--(1.5,0.5);
     \draw[ color= purple!75] (0.5,2.0)--(2.0,0.5);
     \draw[ color= purple!75] (1.0,2.0)--(2.0,1.0);
     \draw[ color= purple!75] (1.5,2.0)--(2.0,1.5);
 \draw[thick, dashed, color= purple!75] (0.5, 0.5) rectangle (2,2);

\draw[dashed,thick, color=gray] (0.75,1.0)--(1.0,1.0)--(1.5,1.5)--(1.5,1.75)--(1.25,1.75)--(0.75,1.25)--(0.75,1.0);



%% file: Pictures/matrix_figure/s3.tex
    \fill[color = gray!20] (1.005,1.005) rectangle (1.505,1.505);

    \draw[step=0.5cm,color=black] (0,0) grid (2.5,2.5);

        

     \draw[ color= blue!75] (0.75,1.0)--(1.5,1.75);
     \draw[ color= blue!75] (1.25,1.0)--(1.5,1.25);
     \draw[ color= blue!75] (0.5,1.25)--(0.75,1.5);
     \draw[ color= blue!75] (1.0,1.75)--(1.25,2.0);
  \draw[ color= blue!75] (0.5,1.0)--(1.5,2.0);
     \draw[ color= blue!75] (1.0,1.0)--(1.5,1.5);
 \draw[thick, dashed, color= blue!75] (0.5, 1.0)--(1.5,1.0)--(1.5,2.0)--(1.0,2.0)--(1.0,1.5)--(0.5,1.5)--(0.5, 1.0) ;

\draw[dashed,thick, color=gray] (0.75,1.0)--(1.0,1.0)--(1.5,1.5)--(1.5,1.75)--(1.25,1.75)--(0.75,1.25)--(0.75,1.0);

  

%% file: Pictures/matrix_figure/s4.tex
    \fill[color = gray!20] (1.005,1.005) rectangle (1.505,1.505);

    \draw[step=0.5cm,color=black] (0,0) grid (2.5,2.5);
   
        


\draw[dashed,thick, color=gray] (0.75,1.0)--(1.0,1.0)--(1.5,1.5)--(1.5,1.75)--(1.25,1.75)--(0.75,1.25)--(0.75,1.0);

    \fill[ color =red ] (1, 1.25) rectangle (1.25, 1.5);

%% file: Pictures/matrix_figure/s5.tex
    \fill[color = gray!20] (1.005,1.245) rectangle (1.245,1.495);

    \draw[step=0.25cm,color=black] (0,0) grid (2.5,2.5);

    \draw[ color =orange ,pattern=north west lines, pattern color=orange] (1, 1.25) rectangle (1.25, 1.5);

%% file: Sections_LSIAC/Numerical_Results.tex
\section{Numerical Results} 
\label{sec:NumericalResults}

The goal of LSIAC-MRA is to produce a fine-mesh modal approximation from coarse-mesh data that has lower error than repeated $L^2$-projection alone. The measures used for comparison are the $L^2-$ and $L^{\infty}-$errors defined respectively by
\begin{equation}
       \norm{u_{exact}-u_{approx}}_0=\sqrt{\frac{1}{|\Omega|}\int_{\Omega}|u_{exact}-u_{approx}|^2\;d\Omega}\;,\qquad
       \norm{u_{exact}-u_{approx}}_{\infty}=\sup_{x\in\Omega}|u_{exact}-u_{approx}|.
   \end{equation}
For our simulations, the $L^2$ error is evaluated through Gauss-Legendre quadrature at $6^d$ nodes per element. Similarly, the $L^{\infty}$ is taken to be the maximum {absolute} error over these nodes. {To ensure a consistent standard for comparison, we compute the errors on the finest mesh occurring during the enhancement procedure. It is very important to use a standard mesh for computing the errors. Not doing so can lead to counter-intuitive behaviors in the error. An example of this is displayed in Figure \ref{fig:Error_Mesh} where repeated projections of the same polynomial approximation have varying errors, contradicting the reproduction property of the projection. The cause is the discrete approximation of the error-norms. By performing the error computation over a fixed sufficiently fine mesh, we no longer encounter these discretization artifacts.}
In all the simulations considered below, periodic boundary conditions are assumed to allow the application of the filter in regions where the kernel support would surpass the domain boundaries. {The first two columns of the tables detail the projections errors of the $L^2$ approximation subject to scale transition by projection without filtering. This represents the standard which our filtered approximations must outperform to justify the computational expense. The third and fourth columns detail the errors of an approximation when LSIAC-MRA is applied only once and then $L^2-$projection is used. The final two columns detail the errors of the LSIAC-MRA procedure applied at each refinement. As a result, the first rows are equal. Similarly, the Enhanced Once and Enhanced Each Refinement columns will agree on the second {mesh simply} because only one filtering procedure has been applied at that point. These equalities are emphasized by the shading of entries in the tables (gray and medium gray, respectively). We have performed numerous two- and three-dimensional experiments to gauge the ability of the LSIAC-MRA filtering-projection procedure, but below we will only focus on the most illuminating examples. Note that we only expect reduced errors, not improvement in order of convergence.}

\subsection{2D Test Problems}
\begin{figure}[bp!]

        \centering

         \scalebox{0.75}{
         \begin{tabular}{c c}
              
       \includegraphics[width=0.6\linewidth]{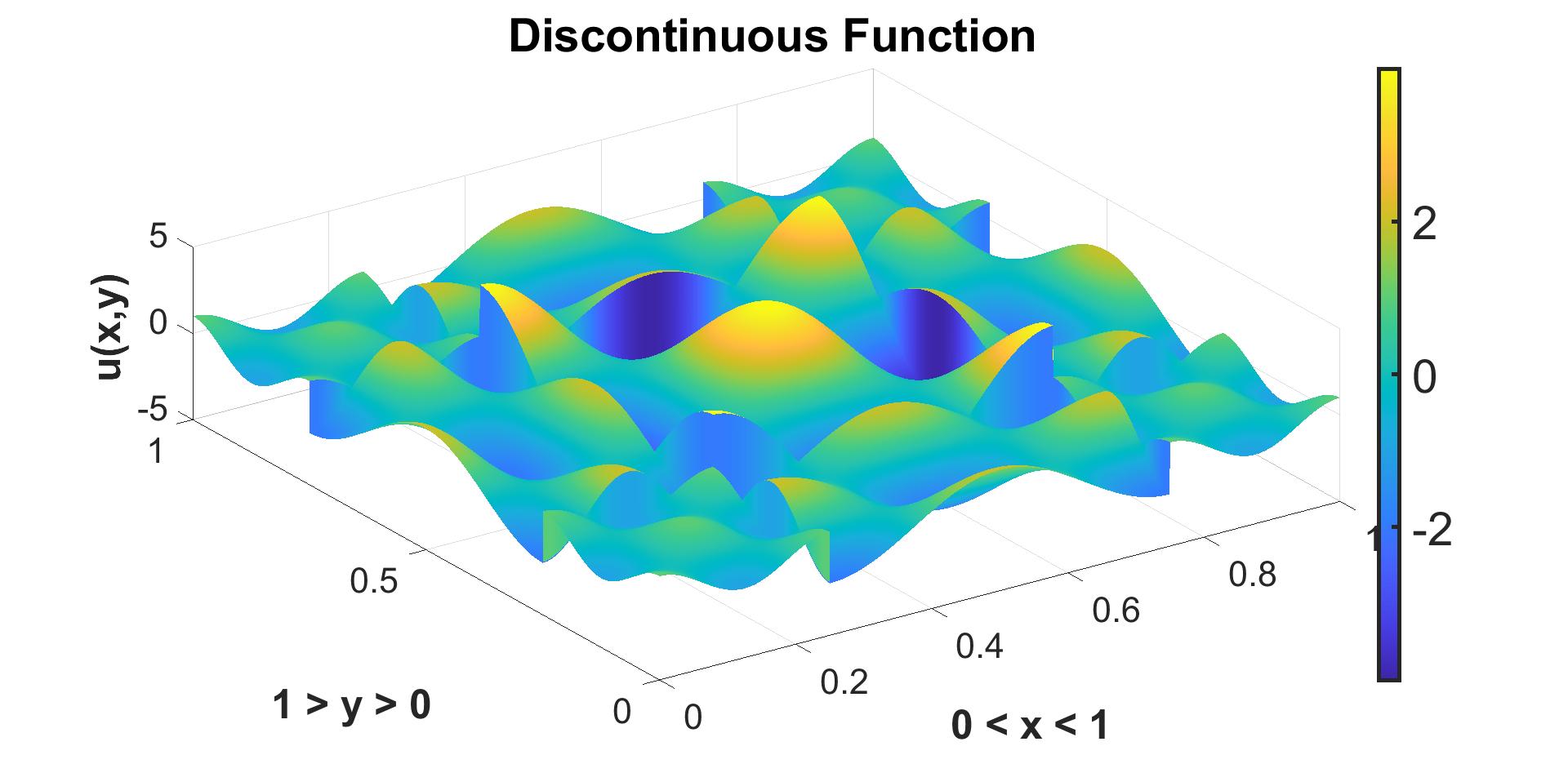}&
             \includegraphics[width=0.6\linewidth]{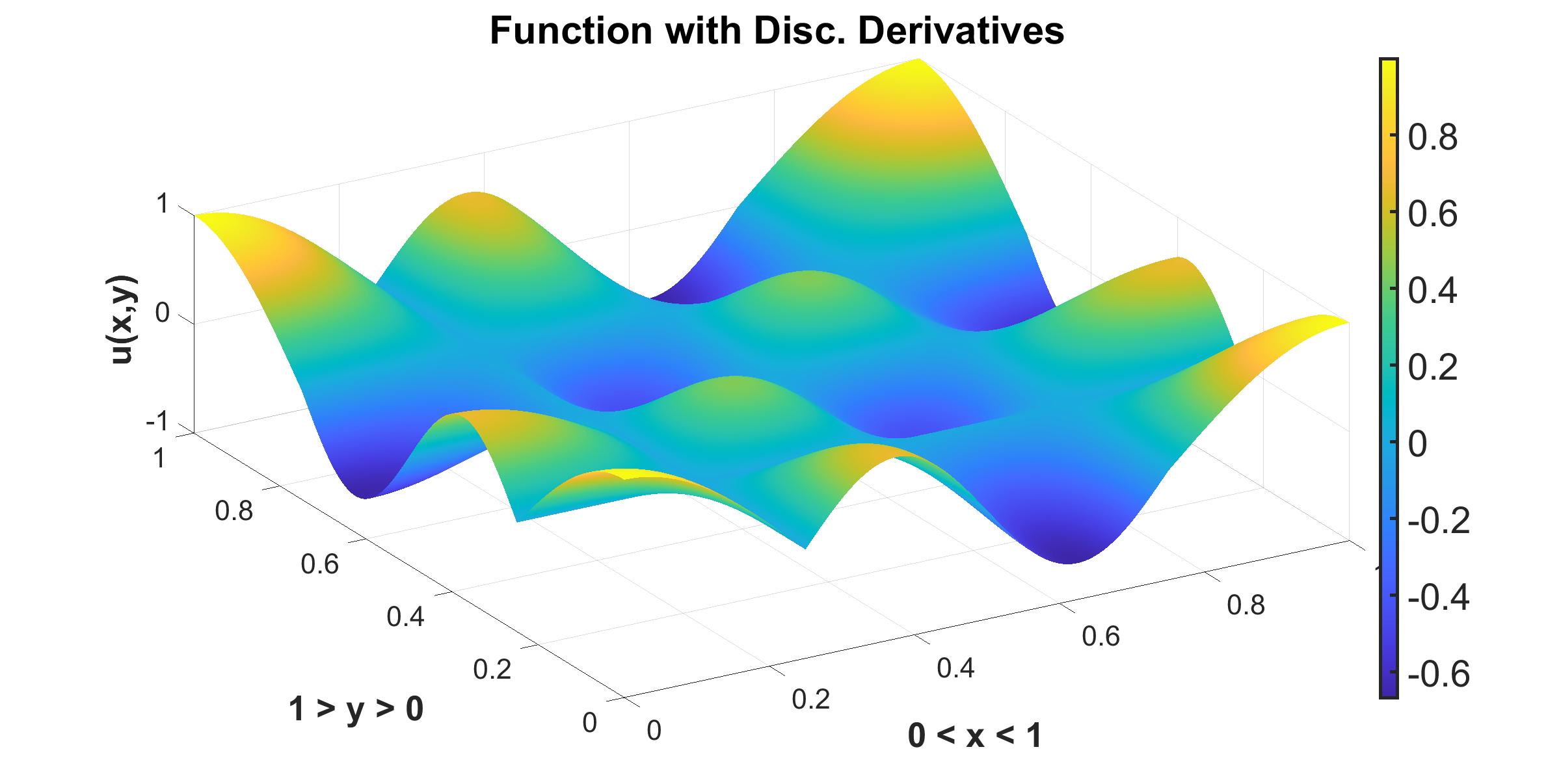}
                 \end{tabular}
        } \caption{{Surface plots of $u_0(x,y)=f_1(x)f_1(y)$ (Left) and $u_0(x,y)=f_2(x)f_2(y)$ (Right). The initial condition defined as a tensor product of a discontinuous function $f_1$ represents a ``strong discontinuity", while that using a tensor product of a continuous but not globally differentiable function $f_2$ represents a continuous function with discontinuous derivatives $\frac{\partial u_0}{\partial x}$ and $\frac{\partial u_0}{\partial y}$.}}
                    \label{fig:2D_surfplot}
    \end{figure}

Consider the following functions on the domain $\Omega=\{(x,y)\in [0,1]^2\}:$
\begin{enumerate}
    \item     $u_0(x,y)=\sin(2\pi (x+y))$
    \item   { $u_0(x,y)=\sin(10\pi x)\sin(10\pi y)$}
        \item   { $u_0(x,y)=f_1(x)f_1(y)$}
            \item   { $u_0(x,y)=f_2(x)f_2(y)$}
\end{enumerate}
where
\begin{align*}
   {  f_1(x)}&=\begin{cases}
    2\cos(2\pi(2x-1)),   &1/4\leq x \leq 3/4\\
    \cos(4\pi(2x-1)),   &Else
    \end{cases},\\
      { f_2(x)}&=\begin{cases}
          2/3\sin(2\pi(2x-1)),     &1/4\leq x \leq 3/4\\
    \cos(\pi(2x-1)),     &Else
    \end{cases}.
\end{align*}
{These tests aid us in evaluating the effectiveness of the LSIAC-MRA procedure for different function types. The first function can be viewed as a combination of a product and summation of the coordinate directions. The second represents a higher frequency function. The third is discontinuous, while the fourth has discontinuous derivatives for $x$ and $y$ at $1/4$ and $3/4$. Depictions of the third and fourth functions are provided in Figure \ref{fig:2D_surfplot}.}

\begin{figure}[bp!]

        \centering

         \scalebox{0.75}{
       \includegraphics[width=1.3\linewidth]{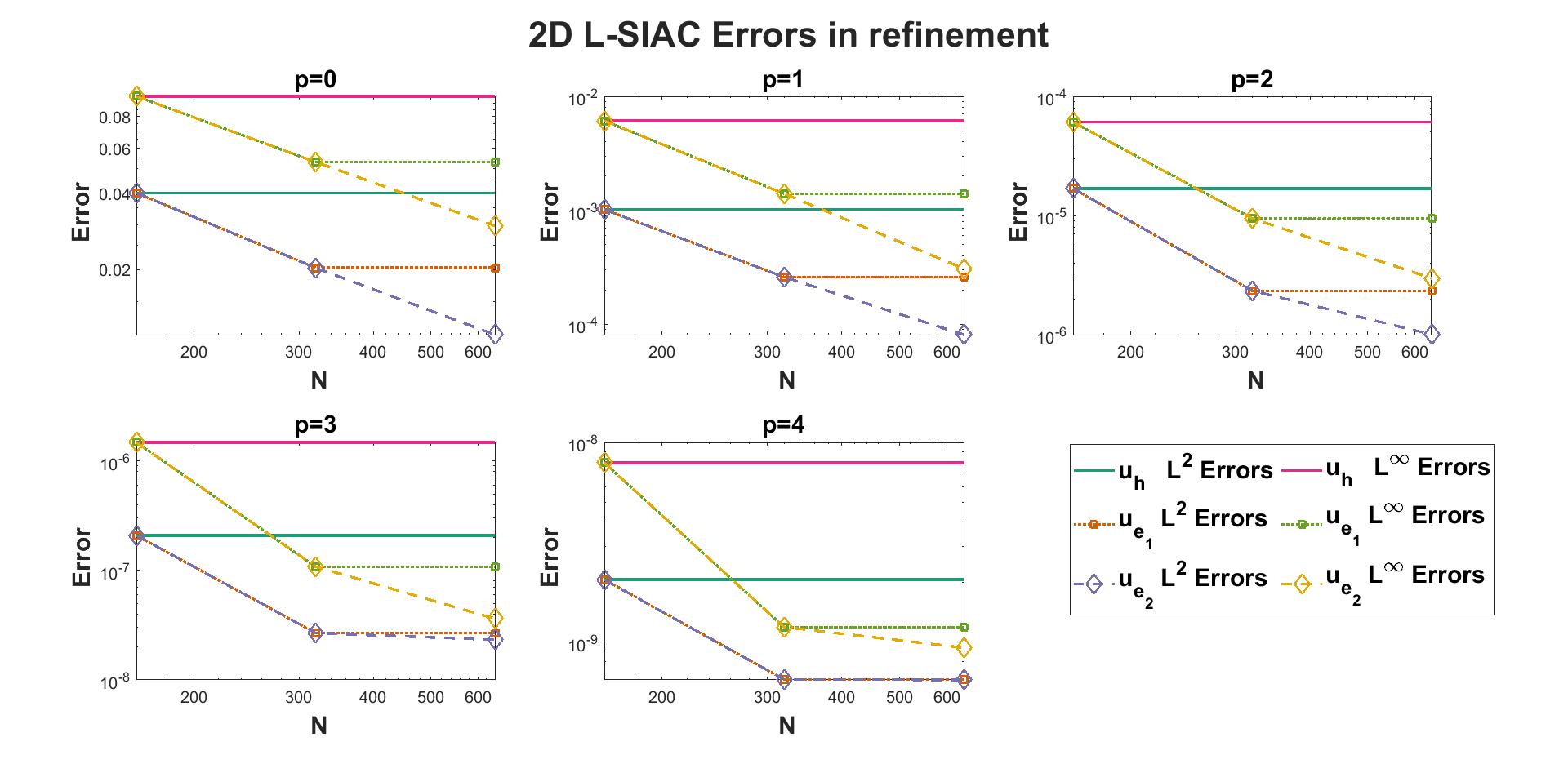}
      
        } \caption{{Log-log plots of $N$ vs. Error for a high frequency function $u_{{0}}(x,y)=\sin(10\pi x)\sin(10\pi y)$ with wavenumber $k=5$ on a series of three meshes using \textbf{2D LSIAC-MRA}. We see that the effectiveness somewhat degrades with increasing polynomial order.}}
                    \label{fig:2D_IC4}
    \end{figure}

{Performing the LSIAC-MRA on functions with different wave number $k$, we observe that different initial mesh resolutions are necessary for improvement. The first and second functions behave similarly, except that the first corresponds to a wavenumber $k=1$ and has repeated improvement from LSIAC-MRA on a relatively coarse mesh ($35\times 35,$ see Figure \ref{fig:2D_contour_IC3_p=4}). This contrasts with the higher frequency second function $k=5$, which requires a much finer starting resolution  of $160\times160$ for improvement. Even then, the resolution requirements increase with approximation order and so for coarser meshes the LSIAC-MRA procedure is more effective for lower order polynomial approximations, see Figure \ref{fig:2D_IC4} and Table \ref{tab:2DLSIACIC4}}.

\begin{figure}[bp!]

        \centering

         \scalebox{0.75}{
         \begin{tabular}{c}
              
       \includegraphics[width=1.29\linewidth]{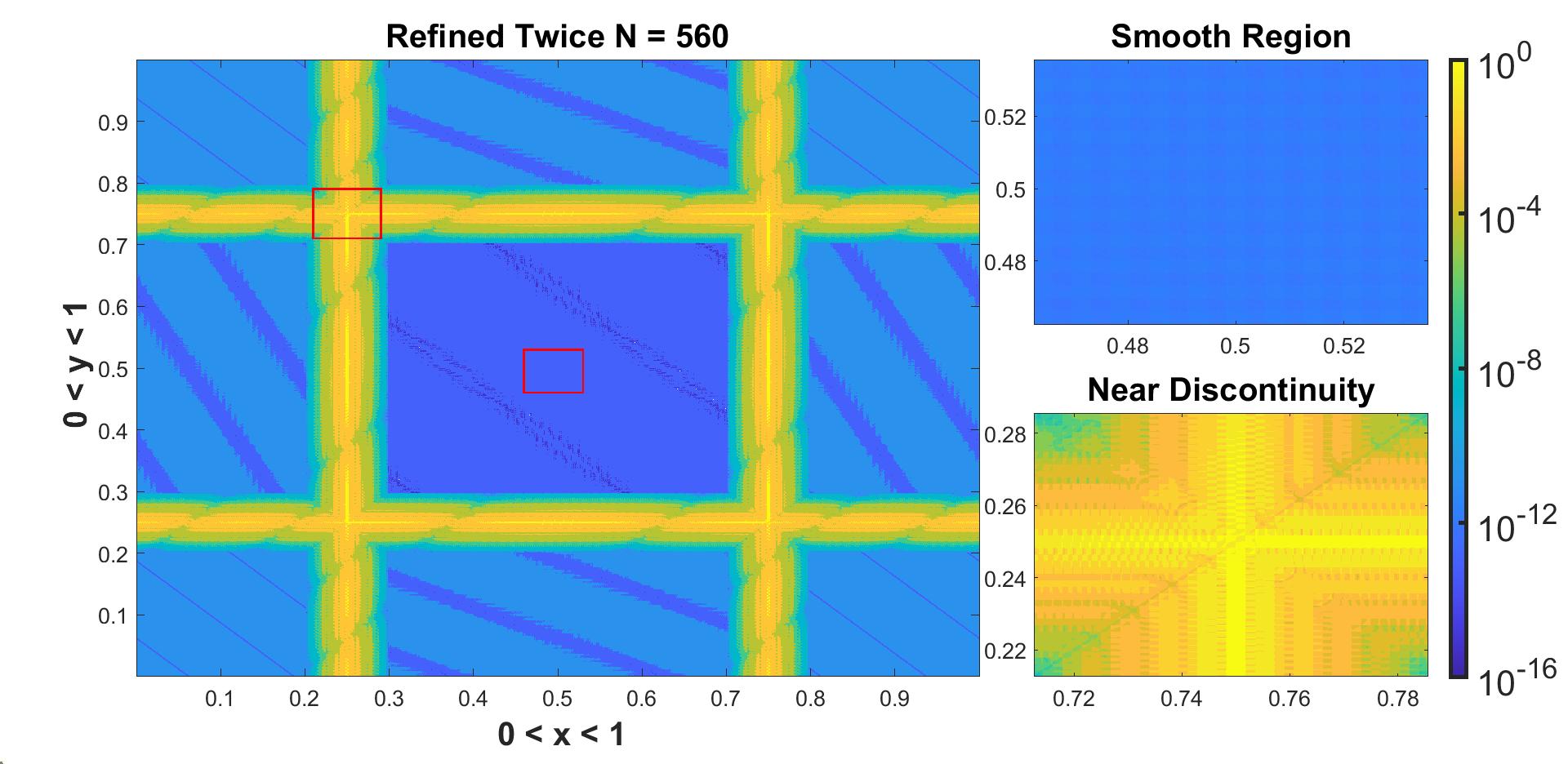}\\
             \includegraphics[width=1.29\linewidth]{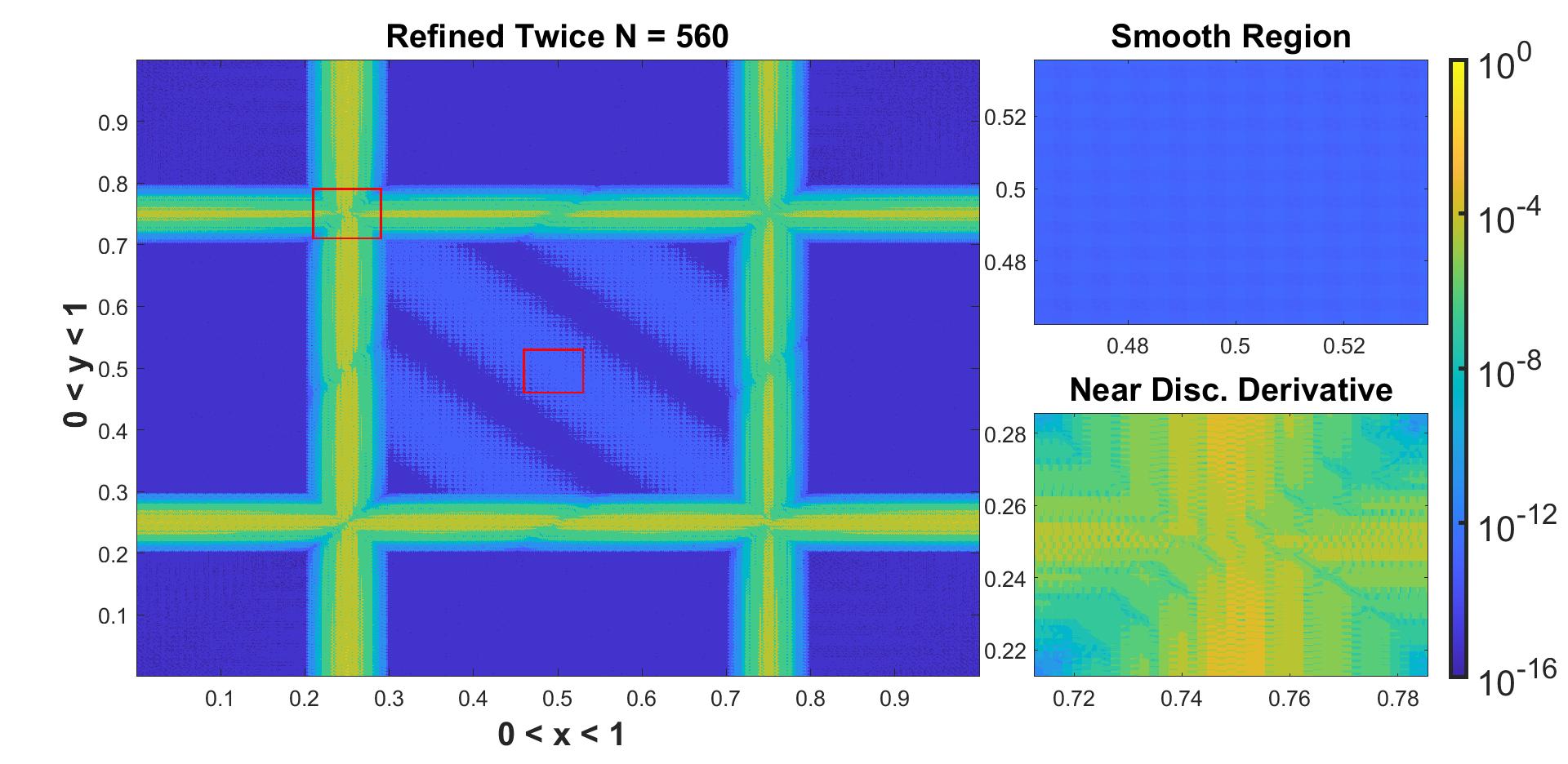}
                 \end{tabular}
        } \caption{{Zoomed in contour error plots for the twice enhanced $p=4$ approximation of a discontinuous function (Top) and a continuous function that is discontinuous in its $\frac{\partial}{\partial x}$ and $\frac{\partial}{\partial y}$ derivatives (Bottom). Both approximations were initialized on a $140\times 140$ mesh.}}
                    \label{fig:2D_zoom}
    \end{figure}

{For the third and fourth functions we begin with their projection onto a piecewise orthonormal Legendre basis defined over a $70\times 70$ uniform quadrilateral mesh. We observe that the refinement procedure provides error reduction for lower polynomial degrees, though performance degrades for increasing $p$. While LSIAC-MRA is effective in smooth regions, it is not as globally effective at error reduction as compared to the analytic examples above. This is due to the lack of smoothness of the functions under consideration. In cases 3 and 4, application of the SIAC filter near a discontinuity will smear the discontinuity by introducing artificial smoothness, and the pollution region will increase with each application of the filter. However, for these functions with non-smooth components occurring along lines, the pollution region will reach a finite extent of $h(2p+3/2)$ in any Cartesian direction from the point of discontinuity, where $h$ is the uniform coarse mesh size. If we look at the contour plot given in Figure \ref{fig:2D_contour_IC6_p=4} or the zoomed in plots of Figure \ref{fig:2D_zoom}, we can in fact observe that near the discontinuous derivative of the function at $x,y=1/4,3/4$, increased errors occur with refinement. This error growth is lesser than compared to the discontinuous test case 3. We do not include refinement contour plots for test case 3 as the error distribution is similar to test case 4. Away from the polluted regions, significant error reduction occurs with refinement. The $L^2-$ and $L^{\infty}-$errors of the enhancement procedure \textit{excluding the polluted regions} for these less smooth functions are given in Tables \ref{tab:2DLSIACIC5} and  \ref{tab:2DLSIACIC6}. We do not include the polluted errors because, as expected, the $L^{\infty}-$errors stay relatively constant while $L^2-$errors grow as the polluted region covers a larger portion of the domain. The pollution regions for tests 3 and 4 are the same, and so we use an identical procedure for determining the pollution regions and computing the errors. The procedure is relatively simple, if an element on the coarse mesh contains a non-smooth feature caused by lack of smoothness in the initial condition, we exclude that element. When we enhance, we wish to exclude any element where the filtered approximation on that element drew information from polluted regions. As the non-smooth components of these functions occur along the lines $x,y=1/4,3/4$ we can exclude any element within $h(p+1/2)$ in any Cartesian direction of the non-smooth coarse elements on the refined mesh after application of LSIAC-MRA. Finally, following the second and final application of LSIAC-MRA, we exclude any element on the finest mesh within $\frac{3}{2}h(p+1/2)$ in any Cartesian direction of the non-smooth coarse elements.} {We observe that application of the filtering procedure at each stage continues to provide error reductions with each refinement and outperforms only applying the filtering procedure once, except in the $p=4$ case for the third function. This is caused by insufficient resolution in the initial mesh, and increasing the initial resolution enables continued error reduction.}

\input{Log-Log_Plots/comp_new_old}

\begin{figure}[bp!]

        \centering

         \scalebox{0.75}{
 
             \includegraphics[width=1.29\linewidth]{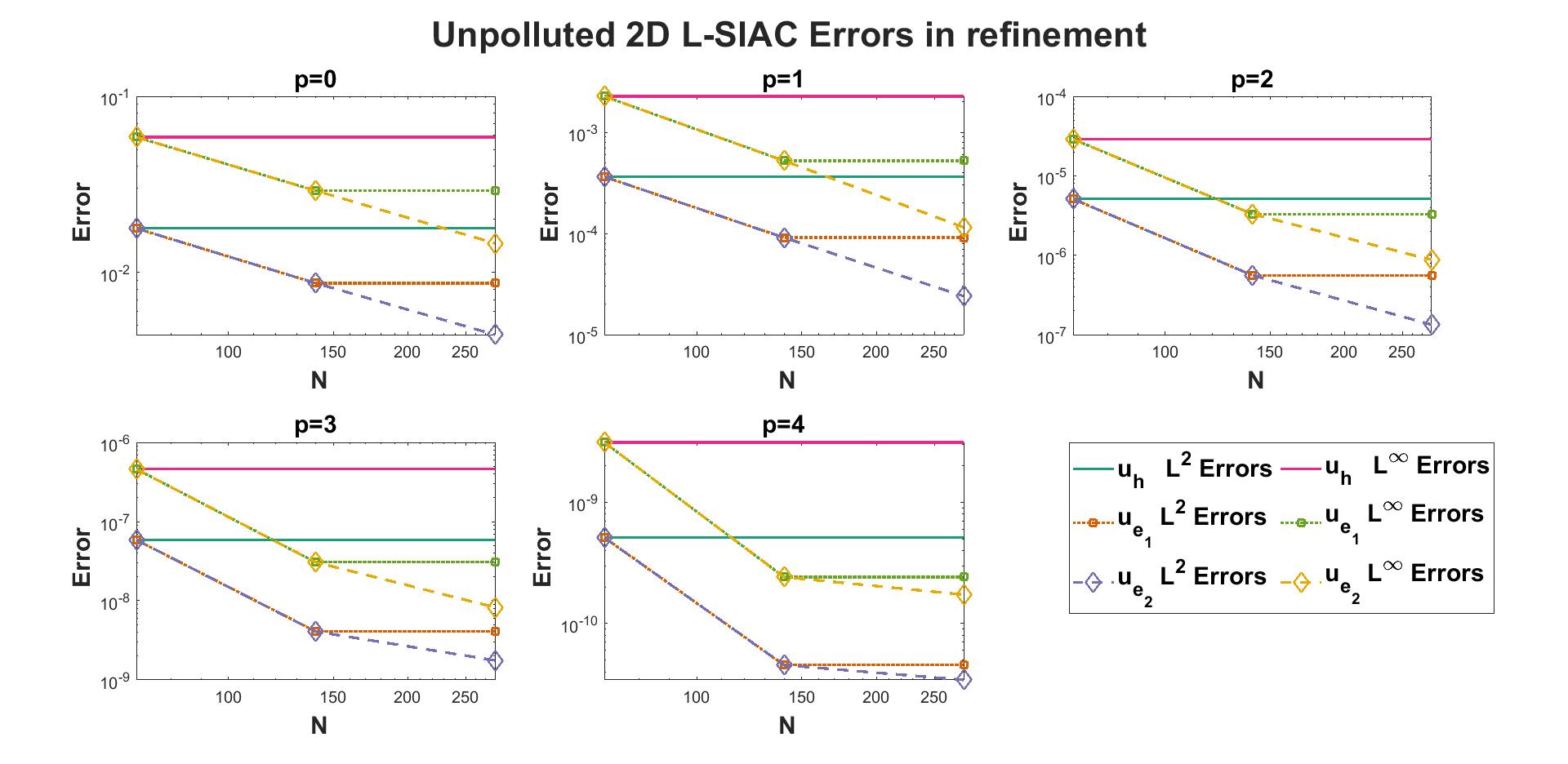}
        } \caption{{Log-log plots of $N$ vs. Error for a continuous function that has discontinuous $\frac{\partial}{\partial x}$ and $\frac{\partial}{\partial y}$ derivatives on a series of three meshes using \textbf{2D LSIAC-MRA}. Errors were calculated by excluding the polluted regions around the discontinuous derivatives. Notice that performance is better than for the discontinuous function.}}
                    \label{fig:2D_IC6}
    \end{figure}

{As depicted in the log-log plots in Figures \ref{fig:2D_IC4} and \ref{fig:2D_IC6}, LSIAC-MRA can indeed provide error reduction for sufficiently smooth functions; however, this technique does not guaranteed an improvement in the order of the method. The contour plots in Figures \ref{fig:2D_contour_IC3_p=4} and \ref{fig:2D_contour_IC6_p=4} display that repeated application of LSIAC-MRA will alter not only the magnitude of errors, but also their distribution throughout the domain.}

\input{2D_Contour_Plots/Group_Contour_Plots_p=4_IC=3}
\input{2D_Contour_Plots/Group_Contour_Plots_p=4_IC=6}

\input{Refinement_Error_Tables/2D_L-SIAC_Error_Tables/2D_L-SIAC_IC=4}
\input{Refinement_Error_Tables/2D_L-SIAC_Error_Tables/2D_L-SIAC_IC=5}
\input{Refinement_Error_Tables/2D_L-SIAC_Error_Tables/2D_L-SIAC_IC=6}





   \clearpage

  \subsection{3D Test Problems}
We now investigate the application of 3D LSIAC-MRA for three-dimensional functions. We have tested the procedure on other functions as well, but present results only for the function $u_0(x,y,z)=\sin(2\pi x)+\sin(2\pi y)+\sin(2\pi z)$ on the domain $\Omega=\{(x,y)\in[0,1]^3\}.$

\begin{figure}[hbt!]
        \centering
   
         \scalebox{0.75}{
            \includegraphics[width=1.0\linewidth]{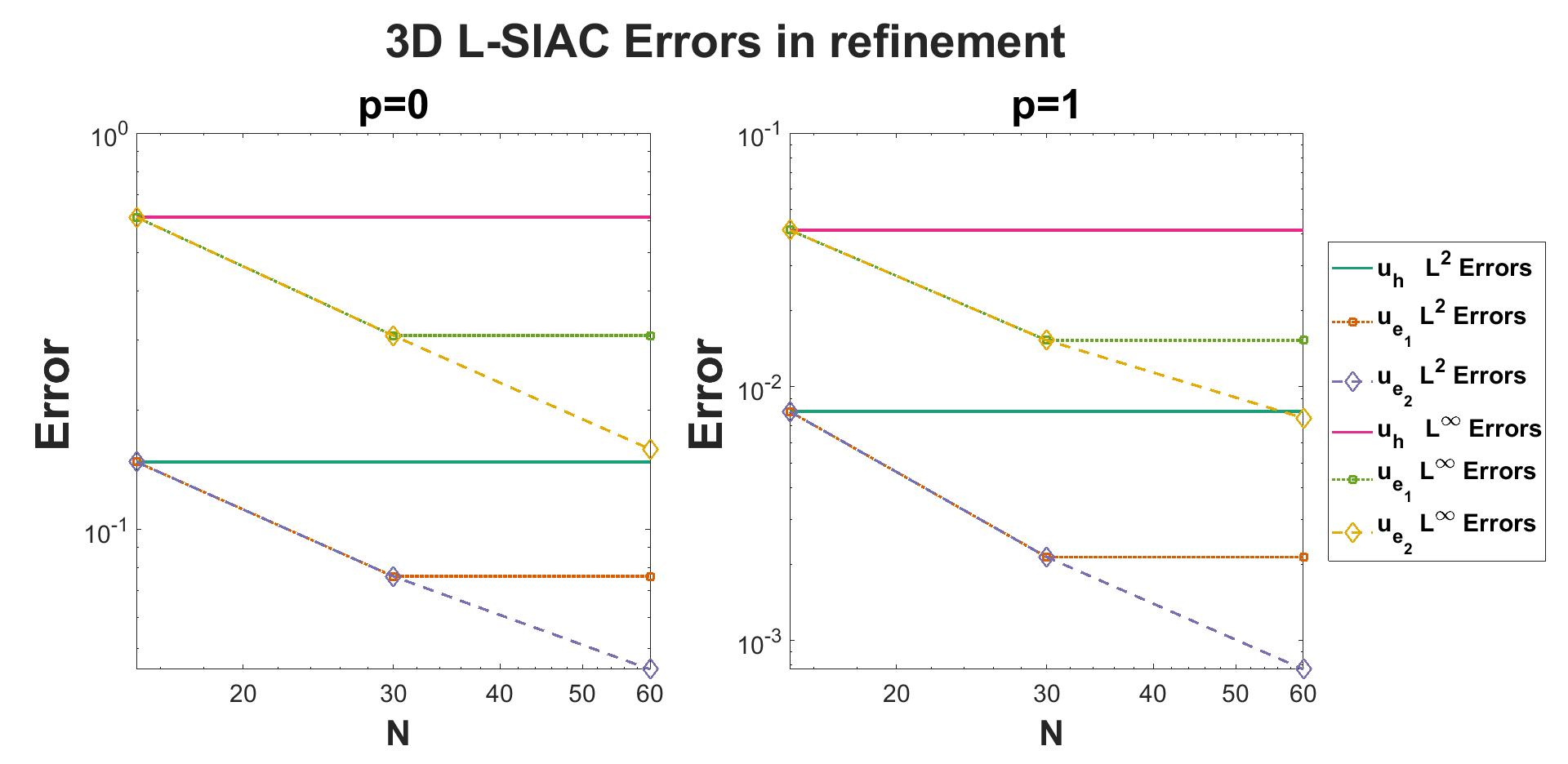}
            }
                 \caption{Log-log plots of $N$ vs. {Error for the function $u_0(x,y,z)=\sin(2\pi x)+\sin(2\pi y)+\sin(2\pi z) $ on a series of three meshes using \textbf{3D LSIAC-MRA}. Notice that procedure improves the approximation and provides error reduction.}  }
        \label{fig:3D_IC2}
    \end{figure}

Beginning with a projection of this function onto a piecewise orthonormal Legendre basis on a $15\times 15\times 15$ uniform hexahedral mesh, the approximation is then enhanced using the LSIAC-MRA procedure. This procedure is applied two times, resulting in a final approximation on a $60\times 60\times 60$ mesh. The $L^2-$ and $L^{\infty}-$errors of those approximations with and without the filtering enhancement are given in {Table} \ref{tab:3DLSIACIC2}.  Notice that we have included the results for a piecewise constant approximation, $p=0$.  The theory of (L)SIAC does not generally extend to piecewise constants. However, it is evident from the results in these tables that the LSIAC-MRA procedure is effective for piecewise constant approximations. For piecewise linear approximations, $p=1$, we can see that it is less effective, but still provides an error reduction. {Similar to the two-dimnesional case, we} speculate that for higher degree polynomial approximations it is necessary to start with a higher resolution in order for the LSIAC-MRA procedure to be more effective. {The} log-log plots for the results of the 3D LSIAC-MRA procedure are given in {Figure} \ref{fig:3D_IC2}.

\input{Refinement_Error_Tables/3D_L-SIAC_Error_Tables/3D_L-SIAC_IC=2}

 \clearpage

%% file: Log-Log_Plots/comp_new_old.tex
    \begin{figure}[bp!]
        \centering
                 \scalebox{0.75}{ \includegraphics[width=1.0\linewidth]{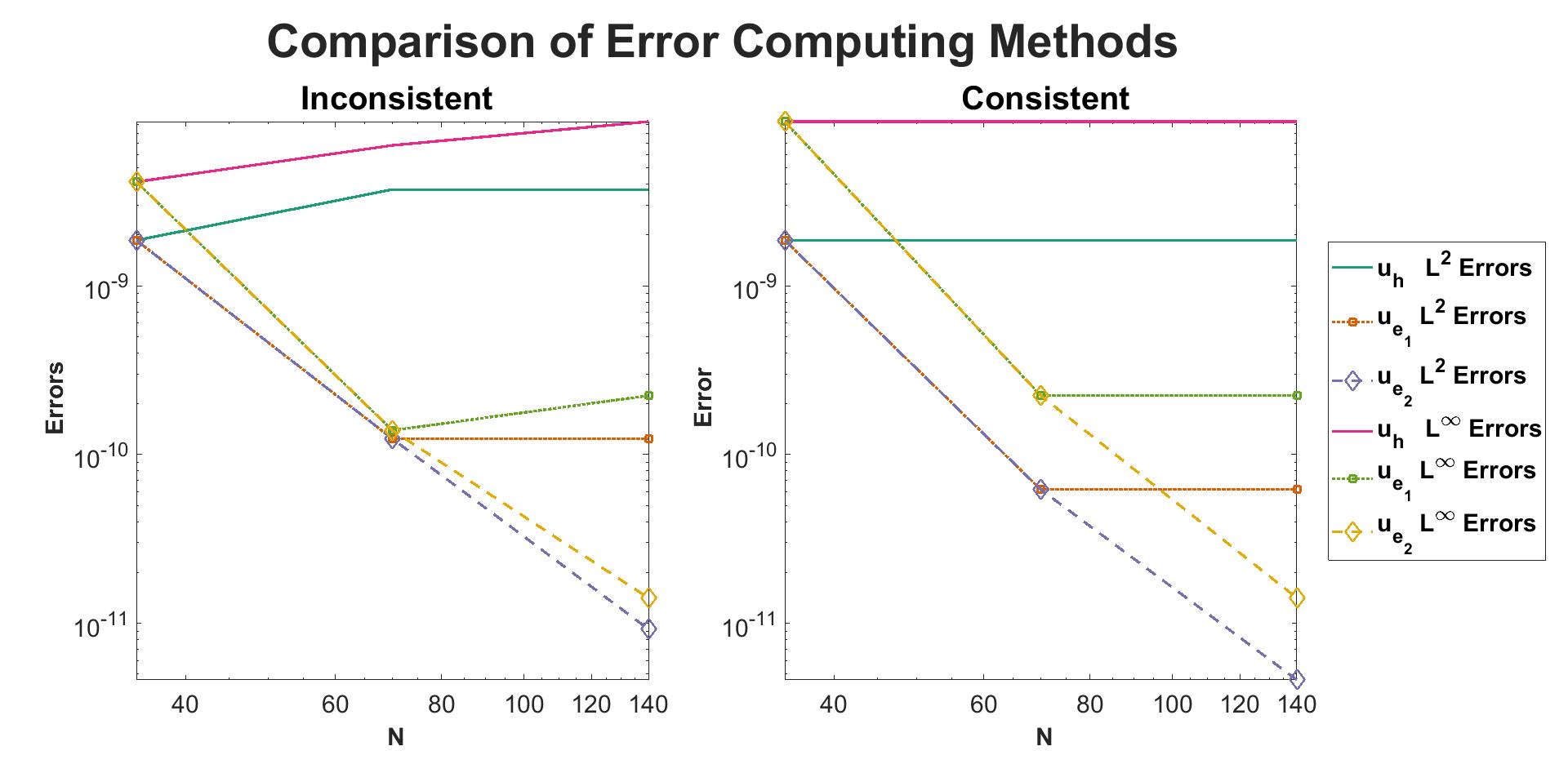}
                 }
        \caption{{Comparison of the log-log plots of the $L^2$ and $L^{\infty}$ errors of the function $u_0(x,y)=\sin(2\pi x)+\sin(2\pi y)$ with $p=4$ under refinement. In the left plot, the errors are computed on the mesh level of the iteration, i.e on the mesh size indicated by the data points abscissa. In the right plot, the errors are computed on the $140\times 140$ mesh which is the finest mesh encountered during the procedure. We observe that this latter approach allows for the errors to remain constant under $L^2$-projection alone, which agrees with the intuition of the approximation being reproduced.}}
                \label{fig:Error_Mesh}

    \end{figure}

%% file: 2D_Contour_Plots/Group_Contour_Plots_p=4_IC=3.tex
    \begin{figure}[bp!]

        \centering
      
       \includegraphics[width=0.97\linewidth]{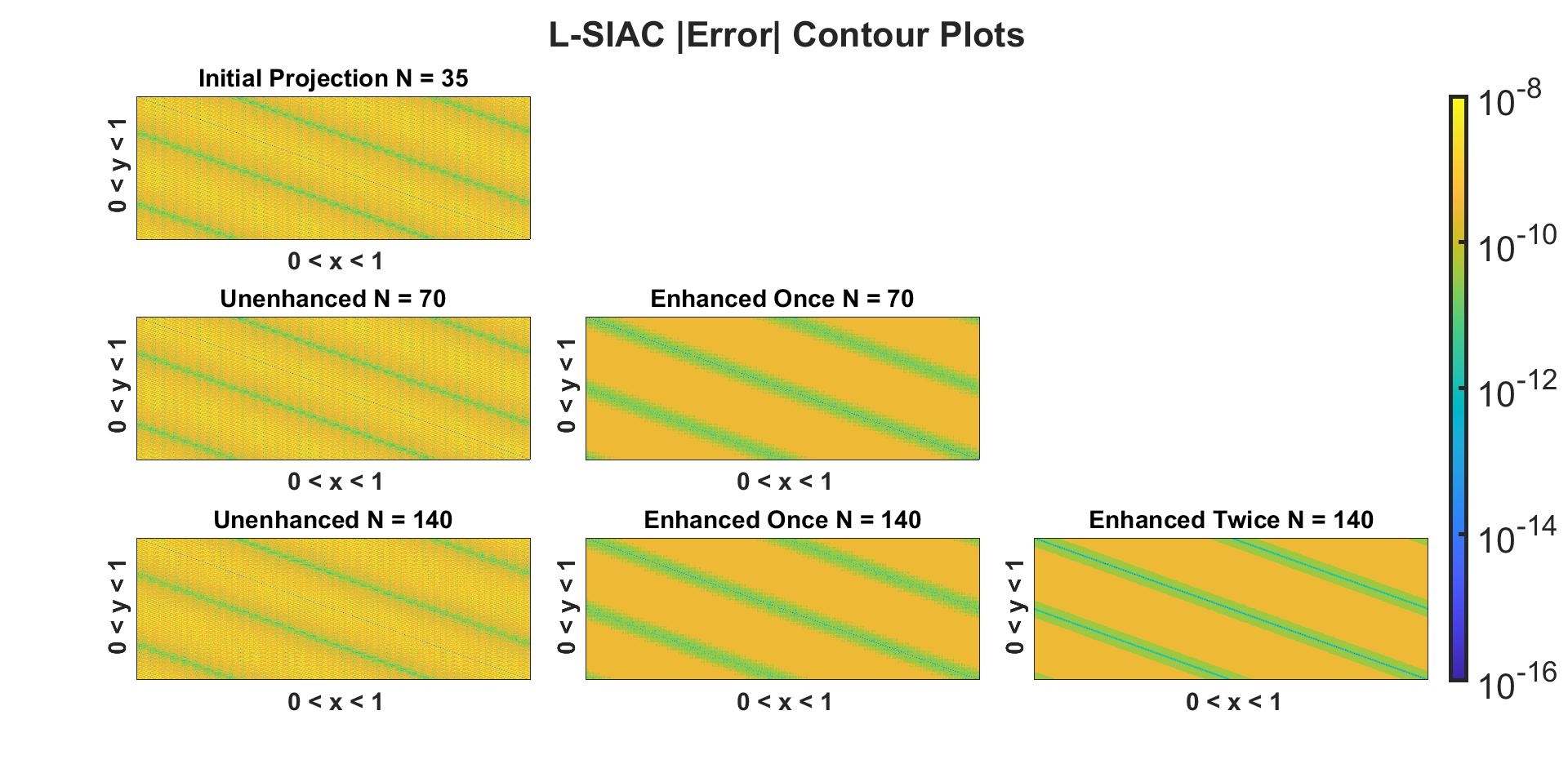}
          \caption{Contour plots depicting pointwise $|$Error$|$ for initial condition $u_{{0}}(x,y)=\sin(2\pi (x+y))$ on a series of three meshes using {\textbf{2D LSIAC-MRA} for $p=4$. Notice that the filtering procedure regularizes and reduces the errors when it is applied at each refinement.}}
        \label{fig:2D_contour_IC3_p=4}
    \end{figure}

%% file: 2D_Contour_Plots/Group_Contour_Plots_p=4_IC=6.tex
    \begin{figure}[bp!]

        \centering
      
       \includegraphics[width=0.97\linewidth]{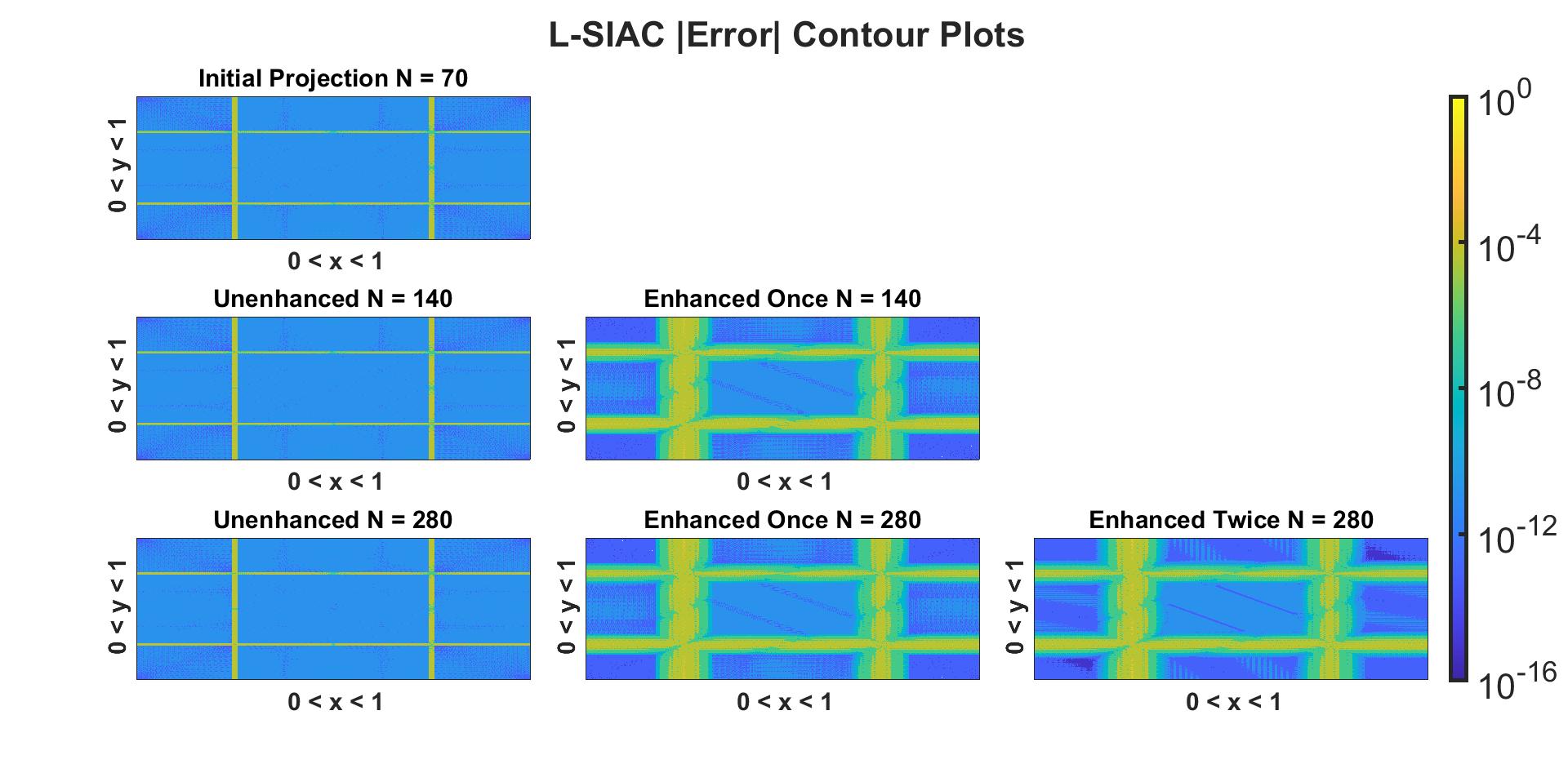}
          \caption{{Contour plots depicting pointwise $|$Error$|$ for the initial condition $u_0(x,y)=f_2(x)f_2(y)$ which contains a discontinuity in its first derivatives derivatives $\frac{\partial u_0}{\partial x}$ and $\frac{\partial u_0}{\partial y}$. 2D LSIAC-MRA was performed over a series of three meshes for $p=4$. Notice that owing to the jump discontinuity in the derivatives at $x,y=1/4,3/4$, the polluted region where filtering adds artificial smoothness grows with each enhancement. Away from these locations, we still have error improvement. A zoomed in view is given in Figure \ref{fig:2D_zoom}.}}
        \label{fig:2D_contour_IC6_p=4}
    \end{figure}

%% file: Refinement_Error_Tables/2D_L-SIAC_Error_Tables/2D_L-SIAC_IC=4.tex
\begin{table}[hbt!]\centering
  \small
    \caption{\label{tab:2DLSIACIC4} {Table of $L^2$ and $L^{\infty}$ errors for the high frequency initial condition $u(x,y)=\sin(10\pi x)\sin(10\pi y)$ on a series of three meshes using \textbf{2D LSIAC-MRA}. Because it is high frequency, we have to initialize our approximation on a more refined grid. The shaded regions of the same shade represent the same approximation. For higher wave numbers the enhancement procedure requires greater initial mesh resolution for error reduction.}}
\begin{tabular}{|c||cc|cc|cc|}  
\hline &\multicolumn{2}{c|}{\textbf{Projection Error}}&\multicolumn{2}{c|}{\textbf{Enhanced Once}} &\multicolumn{2}{c|}{\textbf{Enhanced Each Refinement}}\\ \hline \textbf{N}& \multicolumn{1}{c }{\textbf{$L^2$-Error}}&\multicolumn{1}{c| }{\textbf{$L^{\infty}$-Error}}&\multicolumn{1}{c }{\textbf{$L^2$-Error}}&\multicolumn{1}{c }{\textbf{$L^{\infty}$-Error}}&\multicolumn{1}{|c }{\textbf{$L^2$-Error}}&\multicolumn{1}{c| }{\textbf{$L^{\infty}$-Error}}\\ \hline  &\multicolumn{6}{c|}{\textbf{$\mathbb{P}^{\mathbf{0}}$}}\\ \hline$\mathbf{160}$	&\cellcolor[gray]{0.9}	$4.01e-02$	&\cellcolor[gray]{0.9}	$9.59e-02$	&\cellcolor[gray]{0.9}	$4.01e-02$	&\cellcolor[gray]{0.9}	$9.59e-02$	&\cellcolor[gray]{0.9}	$4.01e-02$	&\cellcolor[gray]{0.9}	$9.59e-02$\\
	\hline
$\mathbf{320}$	&	$4.01e-02$	&	$9.59e-02$	&\cellcolor[gray]{0.85}	$2.04e-02$	&\cellcolor[gray]{0.85}	$5.30e-02$	&\cellcolor[gray]{0.85}	$2.04e-02$	&\cellcolor[gray]{0.85}	$5.30e-02$\\
	\hline
$\mathbf{640}$	&	$4.01e-02$	&	$9.59e-02$	&	$2.04e-02$	&	$5.30e-02$	&\cellcolor[gray]{0.8}	$1.11e-02$	&\cellcolor[gray]{0.8}	$2.97e-02$\\
	\hline
  &\multicolumn{6}{c|}{\textbf{$\mathbb{P}^{\mathbf{1}}$}}\\ \hline$\mathbf{160}$	&\cellcolor[gray]{0.9}	$1.02e-03$	&\cellcolor[gray]{0.9}	$6.10e-03$	&\cellcolor[gray]{0.9}	$1.02e-03$	&\cellcolor[gray]{0.9}	$6.10e-03$	&\cellcolor[gray]{0.9}	$1.02e-03$	&\cellcolor[gray]{0.9}	$6.10e-03$\\
	\hline
$\mathbf{320}$	&	$1.02e-03$	&	$6.10e-03$	&\cellcolor[gray]{0.85}	$2.60e-04$	&\cellcolor[gray]{0.85}	$1.40e-03$	&\cellcolor[gray]{0.85}	$2.60e-04$	&\cellcolor[gray]{0.85}	$1.40e-03$\\
	\hline
$\mathbf{640}$	&	$1.02e-03$	&	$6.10e-03$	&	$2.60e-04$	&	$1.40e-03$	&\cellcolor[gray]{0.8}	$8.11e-05$	&\cellcolor[gray]{0.8}	$3.08e-04$\\
	\hline
  &\multicolumn{6}{c|}{\textbf{$\mathbb{P}^{\mathbf{2}}$}}\\ \hline$\mathbf{160}$	&\cellcolor[gray]{0.9}	$1.69e-05$	&\cellcolor[gray]{0.9}	$6.07e-05$	&\cellcolor[gray]{0.9}	$1.69e-05$	&\cellcolor[gray]{0.9}	$6.07e-05$	&\cellcolor[gray]{0.9}	$1.69e-05$	&\cellcolor[gray]{0.9}	$6.07e-05$\\
	\hline
$\mathbf{320}$	&	$1.69e-05$	&	$6.07e-05$	&\cellcolor[gray]{0.85}	$2.33e-06$	&\cellcolor[gray]{0.85}	$9.44e-06$	&\cellcolor[gray]{0.85}	$2.33e-06$	&\cellcolor[gray]{0.85}	$9.44e-06$\\
	\hline
$\mathbf{640}$	&	$1.69e-05$	&	$6.07e-05$	&	$2.33e-06$	&	$9.44e-06$	&\cellcolor[gray]{0.8}	$1.00e-06$	&\cellcolor[gray]{0.8}	$2.97e-06$\\
	\hline
  &\multicolumn{6}{c|}{\textbf{$\mathbb{P}^{\mathbf{3}}$}}\\ \hline$\mathbf{160}$	&\cellcolor[gray]{0.9}	$2.08e-07$	&\cellcolor[gray]{0.9}	$1.48e-06$	&\cellcolor[gray]{0.9}	$2.08e-07$	&\cellcolor[gray]{0.9}	$1.48e-06$	&\cellcolor[gray]{0.9}	$2.08e-07$	&\cellcolor[gray]{0.9}	$1.48e-06$\\
	\hline
$\mathbf{320}$	&	$2.08e-07$	&	$1.48e-06$	&\cellcolor[gray]{0.85}	$2.69e-08$	&\cellcolor[gray]{0.85}	$1.08e-07$	&\cellcolor[gray]{0.85}	$2.69e-08$	&\cellcolor[gray]{0.85}	$1.08e-07$\\
	\hline
$\mathbf{640}$	&	$2.08e-07$	&	$1.48e-06$	&	$2.69e-08$	&	$1.08e-07$	&\cellcolor[gray]{0.8}	$2.32e-08$	&\cellcolor[gray]{0.8}	$3.63e-08$\\
	\hline
  &\multicolumn{6}{c|}{\textbf{$\mathbb{P}^{\mathbf{4}}$}}\\ \hline$\mathbf{160}$	&\cellcolor[gray]{0.9}	$2.06e-09$	&\cellcolor[gray]{0.9}	$7.91e-09$	&\cellcolor[gray]{0.9}	$2.06e-09$	&\cellcolor[gray]{0.9}	$7.91e-09$	&\cellcolor[gray]{0.9}	$2.06e-09$	&\cellcolor[gray]{0.9}	$7.91e-09$\\
	\hline
$\mathbf{320}$	&	$2.06e-09$	&	$7.91e-09$	&\cellcolor[gray]{0.85}	$6.51e-10$	&\cellcolor[gray]{0.85}	$1.19e-09$	&\cellcolor[gray]{0.85}	$6.51e-10$	&\cellcolor[gray]{0.85}	$1.19e-09$\\
	\hline
$\mathbf{640}$	&	$2.06e-09$	&	$7.91e-09$	&	$6.51e-10$	&	$1.19e-09$	&\cellcolor[gray]{0.8}	$6.47e-10$	&\cellcolor[gray]{0.8}	$9.39e-10$\\
	\hline
\end{tabular} 
\end{table} 

%% file: Refinement_Error_Tables/2D_L-SIAC_Error_Tables/2D_L-SIAC_IC=5.tex
\begin{table}[hbt!]\centering
  \small
  \caption{\label{tab:2DLSIACIC5} {Table of $L^2$ and $L^{\infty}$ errors for an initial condition that has discontinuities along $x,y=1/4,3/4$, $u(x,y)=f_1(x)f_1(y)$. The \textbf{2D LSIAC-MRA} procedure was performed on a series of three meshes. Errors were calculated by excluding the polluted regions around the discontinuities. The LSIAC-MRA procedure reduces errors in these unpolluted regions.}}
\begin{tabular}{|c||cc|cc|cc|}  
\hline &\multicolumn{2}{c|}{\textbf{Projection Error}}&\multicolumn{2}{c|}{\textbf{Enhanced Once}} &\multicolumn{2}{c|}{\textbf{Enhanced Each Refinement}}\\ \hline \textbf{N}& \multicolumn{1}{c }{\textbf{$L^2$-Error}}&\multicolumn{1}{c| }{\textbf{$L^{\infty}$-Error}}&\multicolumn{1}{c }{\textbf{$L^2$-Error}}&\multicolumn{1}{c }{\textbf{$L^{\infty}$-Error}}&\multicolumn{1}{|c }{\textbf{$L^2$-Error}}&\multicolumn{1}{c| }{\textbf{$L^{\infty}$-Error}}\\ \hline  &\multicolumn{6}{c|}{\textbf{$\mathbb{P}^{\mathbf{0}}$}}\\ \hline$\mathbf{70}$	&\cellcolor[gray]{0.9}	$1.12e-01$	&\cellcolor[gray]{0.9}	$3.52e-01$	&\cellcolor[gray]{0.9}	$1.12e-01$	&\cellcolor[gray]{0.9}	$3.52e-01$	&\cellcolor[gray]{0.9}	$1.12e-01$	&\cellcolor[gray]{0.9}	$3.52e-01$\\
	\hline
$\mathbf{140}$	&	$1.12e-01$	&	$3.52e-01$	&\cellcolor[gray]{0.85}	$5.58e-02$	&\cellcolor[gray]{0.85}	$1.92e-01$	&\cellcolor[gray]{0.85}	$5.58e-02$	&\cellcolor[gray]{0.85}	$1.92e-01$\\
	\hline
$\mathbf{280}$	&	$1.12e-01$	&	$3.52e-01$	&	$5.58e-02$	&	$1.92e-01$	&\cellcolor[gray]{0.8}	$3.05e-02$	&\cellcolor[gray]{0.8}	$1.07e-01$\\
	\hline
  &\multicolumn{6}{c|}{\textbf{$\mathbb{P}^{\mathbf{1}}$}}\\ \hline$\mathbf{70}$	&\cellcolor[gray]{0.9}	$4.00e-03$	&\cellcolor[gray]{0.9}	$2.55e-02$	&\cellcolor[gray]{0.9}	$4.00e-03$	&\cellcolor[gray]{0.9}	$2.55e-02$	&\cellcolor[gray]{0.9}	$4.00e-03$	&\cellcolor[gray]{0.9}	$2.55e-02$\\
	\hline
$\mathbf{140}$	&	$4.00e-03$	&	$2.55e-02$	&\cellcolor[gray]{0.85}	$9.08e-04$	&\cellcolor[gray]{0.85}	$5.71e-03$	&\cellcolor[gray]{0.85}	$9.08e-04$	&\cellcolor[gray]{0.85}	$5.71e-03$\\
	\hline
$\mathbf{280}$	&	$4.00e-03$	&	$2.55e-02$	&	$9.08e-04$	&	$5.71e-03$	&\cellcolor[gray]{0.8}	$3.55e-04$	&\cellcolor[gray]{0.8}	$1.35e-03$\\
	\hline
  &\multicolumn{6}{c|}{\textbf{$\mathbb{P}^{\mathbf{2}}$}}\\ \hline$\mathbf{70}$	&\cellcolor[gray]{0.9}	$1.15e-04$	&\cellcolor[gray]{0.9}	$6.95e-04$	&\cellcolor[gray]{0.9}	$1.15e-04$	&\cellcolor[gray]{0.9}	$6.95e-04$	&\cellcolor[gray]{0.9}	$1.15e-04$	&\cellcolor[gray]{0.9}	$6.95e-04$\\
	\hline
$\mathbf{140}$	&	$1.15e-04$	&	$6.95e-04$	&\cellcolor[gray]{0.85}	$1.94e-05$	&\cellcolor[gray]{0.85}	$1.04e-04$	&\cellcolor[gray]{0.85}	$1.94e-05$	&\cellcolor[gray]{0.85}	$1.04e-04$\\
	\hline
$\mathbf{280}$	&	$1.15e-04$	&	$6.95e-04$	&	$1.94e-05$	&	$1.04e-04$	&\cellcolor[gray]{0.8}	$1.38e-05$	&\cellcolor[gray]{0.8}	$6.15e-05$\\
	\hline
  &\multicolumn{6}{c|}{\textbf{$\mathbb{P}^{\mathbf{3}}$}}\\ \hline$\mathbf{70}$	&\cellcolor[gray]{0.9}	$2.47e-06$	&\cellcolor[gray]{0.9}	$1.76e-05$	&\cellcolor[gray]{0.9}	$2.47e-06$	&\cellcolor[gray]{0.9}	$1.76e-05$	&\cellcolor[gray]{0.9}	$2.47e-06$	&\cellcolor[gray]{0.9}	$1.76e-05$\\
	\hline
$\mathbf{140}$	&	$2.47e-06$	&	$1.76e-05$	&\cellcolor[gray]{0.85}	$1.09e-06$	&\cellcolor[gray]{0.85}	$4.63e-06$	&\cellcolor[gray]{0.85}	$1.09e-06$	&\cellcolor[gray]{0.85}	$4.63e-06$\\
	\hline
$\mathbf{280}$	&	$2.47e-06$	&	$1.76e-05$	&	$1.09e-06$	&	$4.63e-06$	&\cellcolor[gray]{0.8}	$9.30e-07$	&\cellcolor[gray]{0.8}	$3.86e-06$\\
	\hline
  &\multicolumn{6}{c|}{\textbf{$\mathbb{P}^{\mathbf{4}}$}}\\ \hline$\mathbf{70}$	&\cellcolor[gray]{0.9}	$4.57e-08$	&\cellcolor[gray]{0.9}	$3.00e-07$	&\cellcolor[gray]{0.9}	$4.57e-08$	&\cellcolor[gray]{0.9}	$3.00e-07$	&\cellcolor[gray]{0.9}	$4.57e-08$	&\cellcolor[gray]{0.9}	$3.00e-07$\\
	\hline
$\mathbf{140}$	&	$4.57e-08$	&	$3.00e-07$	&\cellcolor[gray]{0.85}	$8.91e-08$	&\cellcolor[gray]{0.85}	$3.55e-07$	&\cellcolor[gray]{0.85}	$8.91e-08$	&\cellcolor[gray]{0.85}	$3.55e-07$\\
	\hline
$\mathbf{280}$	&	$4.57e-08$	&	$3.00e-07$	&	$8.91e-08$	&	$3.55e-07$	&\cellcolor[gray]{0.8}	$7.17e-08$	&\cellcolor[gray]{0.8}	$3.49e-07$\\
	\hline
\end{tabular} 
\end{table} 

%% file: Refinement_Error_Tables/2D_L-SIAC_Error_Tables/2D_L-SIAC_IC=6.tex
\begin{table}[hbt!]\centering
  \small
  \caption{\label{tab:2DLSIACIC6} {Shown are the $L^2$ and $L^{\infty}$ errors for a function $u_0(x,y)=f_2(x)f_2(y)$ that has discontinuous derivatives $\frac{\partial u_0}{\partial x}$ and $\frac{\partial u_0}{\partial y}$. The \textbf{2D LSIAC-MRA} procedure was performed on a series of three meshes. Errors were calculated by excluding the polluted regions around the discontinuous derivatives. We observe that the LSIAC-MRA is effective for error reduction in unpolluted regions.}}
\begin{tabular}{|c||cc|cc|cc|}  
\hline &\multicolumn{2}{c|}{\textbf{Projection Error}}&\multicolumn{2}{c|}{\textbf{Enhanced Once}} &\multicolumn{2}{c|}{\textbf{Enhanced Each Refinement}}\\ \hline \textbf{N}& \multicolumn{1}{c }{\textbf{$L^2$-Error}}&\multicolumn{1}{c| }{\textbf{$L^{\infty}$-Error}}&\multicolumn{1}{c }{\textbf{$L^2$-Error}}&\multicolumn{1}{c }{\textbf{$L^{\infty}$-Error}}&\multicolumn{1}{|c }{\textbf{$L^2$-Error}}&\multicolumn{1}{c| }{\textbf{$L^{\infty}$-Error}}\\ \hline  &\multicolumn{6}{c|}{\textbf{$\mathbb{P}^{\mathbf{0}}$}}\\ \hline$\mathbf{70}$	&\cellcolor[gray]{0.9}	$1.78e-02$	&\cellcolor[gray]{0.9}	$5.87e-02$	&\cellcolor[gray]{0.9}	$1.78e-02$	&\cellcolor[gray]{0.9}	$5.87e-02$	&\cellcolor[gray]{0.9}	$1.78e-02$	&\cellcolor[gray]{0.9}	$5.87e-02$\\
	\hline
$\mathbf{140}$	&	$1.78e-02$	&	$5.87e-02$	&\cellcolor[gray]{0.85}	$8.71e-03$	&\cellcolor[gray]{0.85}	$2.92e-02$	&\cellcolor[gray]{0.85}	$8.71e-03$	&\cellcolor[gray]{0.85}	$2.92e-02$\\
	\hline
$\mathbf{280}$	&	$1.78e-02$	&	$5.87e-02$	&	$8.71e-03$	&	$2.92e-02$	&\cellcolor[gray]{0.8}	$4.44e-03$	&\cellcolor[gray]{0.8}	$1.47e-02$\\
	\hline
  &\multicolumn{6}{c|}{\textbf{$\mathbb{P}^{\mathbf{1}}$}}\\ \hline$\mathbf{70}$	&\cellcolor[gray]{0.9}	$3.63e-04$	&\cellcolor[gray]{0.9}	$2.27e-03$	&\cellcolor[gray]{0.9}	$3.63e-04$	&\cellcolor[gray]{0.9}	$2.27e-03$	&\cellcolor[gray]{0.9}	$3.63e-04$	&\cellcolor[gray]{0.9}	$2.27e-03$\\
	\hline
$\mathbf{140}$	&	$3.63e-04$	&	$2.27e-03$	&\cellcolor[gray]{0.85}	$9.10e-05$	&\cellcolor[gray]{0.85}	$5.23e-04$	&\cellcolor[gray]{0.85}	$9.10e-05$	&\cellcolor[gray]{0.85}	$5.23e-04$\\
	\hline
$\mathbf{280}$	&	$3.63e-04$	&	$2.27e-03$	&	$9.10e-05$	&	$5.23e-04$	&\cellcolor[gray]{0.8}	$2.40e-05$	&\cellcolor[gray]{0.8}	$1.15e-04$\\
	\hline
  &\multicolumn{6}{c|}{\textbf{$\mathbb{P}^{\mathbf{2}}$}}\\ \hline$\mathbf{70}$	&\cellcolor[gray]{0.9}	$5.10e-06$	&\cellcolor[gray]{0.9}	$2.89e-05$	&\cellcolor[gray]{0.9}	$5.10e-06$	&\cellcolor[gray]{0.9}	$2.89e-05$	&\cellcolor[gray]{0.9}	$5.10e-06$	&\cellcolor[gray]{0.9}	$2.89e-05$\\
	\hline
$\mathbf{140}$	&	$5.10e-06$	&	$2.89e-05$	&\cellcolor[gray]{0.85}	$5.59e-07$	&\cellcolor[gray]{0.85}	$3.30e-06$	&\cellcolor[gray]{0.85}	$5.59e-07$	&\cellcolor[gray]{0.85}	$3.30e-06$\\
	\hline
$\mathbf{280}$	&	$5.10e-06$	&	$2.89e-05$	&	$5.59e-07$	&	$3.30e-06$	&\cellcolor[gray]{0.8}	$1.34e-07$	&\cellcolor[gray]{0.8}	$8.75e-07$\\
	\hline
  &\multicolumn{6}{c|}{\textbf{$\mathbb{P}^{\mathbf{3}}$}}\\ \hline$\mathbf{70}$	&\cellcolor[gray]{0.9}	$5.86e-08$	&\cellcolor[gray]{0.9}	$4.59e-07$	&\cellcolor[gray]{0.9}	$5.86e-08$	&\cellcolor[gray]{0.9}	$4.59e-07$	&\cellcolor[gray]{0.9}	$5.86e-08$	&\cellcolor[gray]{0.9}	$4.59e-07$\\
	\hline
$\mathbf{140}$	&	$5.86e-08$	&	$4.59e-07$	&\cellcolor[gray]{0.85}	$4.09e-09$	&\cellcolor[gray]{0.85}	$3.07e-08$	&\cellcolor[gray]{0.85}	$4.09e-09$	&\cellcolor[gray]{0.85}	$3.07e-08$\\
	\hline
$\mathbf{280}$	&	$5.86e-08$	&	$4.59e-07$	&	$4.09e-09$	&	$3.07e-08$	&\cellcolor[gray]{0.8}	$1.77e-09$	&\cellcolor[gray]{0.8}	$8.27e-09$\\
	\hline
  &\multicolumn{6}{c|}{\textbf{$\mathbb{P}^{\mathbf{4}}$}}\\ \hline$\mathbf{70}$	&\cellcolor[gray]{0.9}	$5.12e-10$	&\cellcolor[gray]{0.9}	$3.12e-09$	&\cellcolor[gray]{0.9}	$5.12e-10$	&\cellcolor[gray]{0.9}	$3.12e-09$	&\cellcolor[gray]{0.9}	$5.12e-10$	&\cellcolor[gray]{0.9}	$3.12e-09$\\
	\hline
$\mathbf{140}$	&	$5.12e-10$	&	$3.12e-09$	&\cellcolor[gray]{0.85}	$4.55e-11$	&\cellcolor[gray]{0.85}	$2.42e-10$	&\cellcolor[gray]{0.85}	$4.55e-11$	&\cellcolor[gray]{0.85}	$2.42e-10$\\
	\hline
$\mathbf{280}$	&	$5.12e-10$	&	$3.12e-09$	&	$4.55e-11$	&	$2.42e-10$	&\cellcolor[gray]{0.8}	$3.43e-11$	&\cellcolor[gray]{0.8}	$1.73e-10$\\
	\hline
\end{tabular} 
\end{table} 

%% file: Refinement_Error_Tables/3D_L-SIAC_Error_Tables/3D_L-SIAC_IC=2.tex
\begin{table}[hbt!]\centering
\small
 \caption{\label{tab:3DLSIACIC2}Table of $L^2$ and $L^{\infty}$ errors for initial condition $u(x,y,z)=\sin(2\pi x)+\sin(2\pi y)+\sin(2\pi z) $ on a series of three meshes using \textbf{3D LSIAC-MRA}. We observe error reduction with each enhancement.}
\begin{tabular}{|c||cc|cc|cc|} 
\hline &\multicolumn{2}{c|}{\textbf{Projection Error}}&\multicolumn{2}{c|}{\textbf{Enhanced Once}} &\multicolumn{2}{c|}{\textbf{Enhanced Each Refinement}}\\ \hline $\mathbf{N}$& \multicolumn{1}{c }{$\mathbf{L^2}$\textbf{-Error}}&\multicolumn{1}{c| }{$\mathbf{L^{\infty}}$\textbf{-Error}}&\multicolumn{1}{c }{$\mathbf{L^2}$\textbf{-Error}}&\multicolumn{1}{c }{$\mathbf{L^{\infty}}$\textbf{-Error}}&\multicolumn{1}{|c }{$\mathbf{L^2}$\textbf{-Error}}&\multicolumn{1}{c| }{$\mathbf{L^{\infty}}$\textbf{-Error}}\\ \hline  &\multicolumn{6}{c|}{$\mathbf{\mathbb{P}^0}$}\\ \hline$\mathbf{15}$	&\cellcolor[gray]{0.9}	$1.48e-01$	&\cellcolor[gray]{0.9}	$6.13e-01$	&\cellcolor[gray]{0.9}	$1.48e-01$	&\cellcolor[gray]{0.9}	$6.13e-01$	&\cellcolor[gray]{0.9}	$1.48e-01$	&\cellcolor[gray]{0.9}	$6.13e-01$\\
	\hline
$\mathbf{30}$	&	$1.48e-01$	&	$6.13e-01$	&\cellcolor[gray]{0.85}	$7.61e-02$	&\cellcolor[gray]{0.85}	$3.08e-01$	&\cellcolor[gray]{0.85}	$7.61e-02$	&\cellcolor[gray]{0.85}	$3.08e-01$\\
	\hline
$\mathbf{60}$	&	$1.48e-01$	&	$6.13e-01$	&	$7.61e-02$	&	$3.08e-01$	&\cellcolor[gray]{0.8}	$4.43e-02$	&\cellcolor[gray]{0.8}	$1.59e-01$\\
	\hline
  &\multicolumn{6}{c|}{$\mathbf{\mathbb{P}^1}$}\\ \hline$\mathbf{15}$	&\cellcolor[gray]{0.9}	$7.99e-03$	&\cellcolor[gray]{0.9}	$4.14e-02$	&\cellcolor[gray]{0.9}	$7.99e-03$	&\cellcolor[gray]{0.9}	$4.14e-02$	&\cellcolor[gray]{0.9}	$7.99e-03$	&\cellcolor[gray]{0.9}	$4.14e-02$\\
	\hline
$\mathbf{30}$	&	$7.99e-03$	&	$4.14e-02$	&\cellcolor[gray]{0.85}	$2.13e-03$	&\cellcolor[gray]{0.85}	$1.53e-02$	&\cellcolor[gray]{0.85}	$2.13e-03$	&	\cellcolor[gray]{0.85}$1.53e-02$\\
	\hline
$\mathbf{60}$	&	$7.99e-03$	&	$4.14e-02$	&	$2.13e-03$	&	$1.53e-02$	&\cellcolor[gray]{0.8}	$7.70e-04$	&	\cellcolor[gray]{0.8}$7.53e-03$\\
	\hline
\end{tabular} 
\end{table}

%% file: Sections_LSIAC/Conclusions.tex
\section{Conclusions and Future Work}
\label{sec:conclusions}

In this article, we have introduced an improved multi-resolution analysis scheme for multi-dimensional applications, LSIAC-MRA. This scheme utilizes the Line Smoothness-Increasing Accuracy-Conserving filter which {post}-processes multi-dimensional data using a one-dimensional support{.} {This allows for} approximating the difference coefficients in the multi-wavelet representation {and} allows for {error reduction} with mesh refinement. {We have provided the underlying operational framework and demonstrated that LSIAC-MRA is effective for two- and three-dimensional applications}.  Furthermore, LSIAC-MRA is shown to be effective for multi-dimensional piecewise constant data as well. {Though the method does not have translational invariance with respect to arbitrary translations owing to the discrete projection, it is invariant with respect to translations by the uniform mesh scaling, $h$.} In the future, we plan on investigating the requirements for the points per wavelength of the initial data as well as allowed non-uniformities in the data in order to leverage LSIAC-MRA for turbulence modeling data or improved identification of undersampled signals.  

%% file: Sections_LSIAC/Acknowledgements.tex
    \section{Acknowledgements}
 This work is partially supported by AFOSR under grant number FA9550-20-1-0166. {We} would like to thank Dr. Ayaboe Edoh and Dr. Julia Docampo-S\'{a}nchez for their valuable comments in the presentation of this paper. {Finally, we would like to thank the reviewers for their useful comments that improved the clarity of the presentation of the material contained in this paper.}

%% file: Supplementary_Material/1D_Basis_Table.tex
\begin{table}[h!]
 \caption{Table of basis functions and associated coefficients resulting from the analytical evaluation of the convolution in $u_h^{\star}(x)$.}
 \label{tab:1Dchi}
    \centering
    \begin{tabular}{|c||c||c|} \hline
        ${\bf k}$    &  ${\bf d^j_k}$ &  ${\bf \tilde{d}^j_k}$  \\ \hline\hline
        $ {\bf 0}$   & $\frac{\sqrt{2}}{2}\big\{\sum_{\gamma=-p}^p c_{\gamma}u^{j+\gamma}_0\big\} $ &  $\frac{\sqrt{2}}{2}\big\{\sum_{\gamma=-p}^p c_{\gamma}u^{j+\gamma}_0\big\}$ \\ \hline
        $ {\bf 1}$   & $\frac{\sqrt{2}}{4}\big\{\sum_{\gamma=-p}^p c_{\gamma}(u^{j+\gamma}_0-u^{j+\gamma-1}_0)\big\}$ & $\frac{\sqrt{2}}{4}\big\{\sum_{\gamma=-p}^p c_{\gamma}(u^{j+\gamma+1}_0-u^{j+\gamma}_0)\big\}$ \\ \hline
      $ {\bf k>1}$ & $\frac{1}{2}\big\{\sum_{\gamma=-p}^p c_{\gamma}(u^{j+\gamma}_k-u^{j+\gamma-1}_k)\big\}$ &  $\frac{1}{2}\big\{\sum_{\gamma=-p}^p c_{\gamma}(u^{j+\gamma+1}_k-u^{j+\gamma}_k)\big\}$ \\ \hline\hline
                          &   ${\bf \chi^k_L(\zeta)}$ &  ${\bf \chi^k_R(\zeta)}$  \\ \hline
     $ {\bf 0}$.     &  $1$&$1$\\ \hline
    $ {\bf 1}$.      & $ \zeta$ &  $\zeta$\\ \hline
$  {\bf k>1}$.     & $\frac{\sqrt{k-1/2}}{2k-1}\Big[P_k(1+\zeta)-P_{k-2}(1+\zeta)\Big]$& $\frac{\sqrt{k-1/2}}{2k-1}\Big[P_k(\zeta-1)-P_{k-2}(\zeta-1)\Big]$\\ \hline
    \end{tabular} 
   
\end{table}

%% file: Supplementary_Material/2D_L-SIAC_Basis_Table.tex
\begin{landscape}
\begin{table}[h!]
    \centering
    \small
   
    \caption{2D L-SIAC basis functions. For notational convenience suppose $\mathbf{x}$ belongs to element $(i,j)$. }
     \label{tab:2D_LSIAC_basis}
\begin{tabular}{|M{1.5cm}||M{6.5cm}|M{10.5cm}|}\hline
{\bf Region} $S$  & {\bf Basis Function} & {\bf Definition}\\ \hline\hline
  {\bf 1T} &$\chi^{m,r}_{1T,L}(\zeta_x,\zeta_y)$ & $\frac{2-\zeta_y}{2}\int_{-1}^1\phi_{kx}\Big(\zeta_x-\frac{\zeta_y}{2}+\big(1-\frac{\zeta_y}{2}\big)\eta\Big)\phi_{ky}\Big(\big(1-\frac{\zeta_y}{2}\big)\eta+\frac{\zeta_y}{2}\Big)\;d\eta$  \\
   
   &$\chi^{m,r}_{1T,M}(\zeta_x,\zeta_y)$ & $\frac{\zeta_y-\zeta_x}{2}\int_{-1}^1\phi_{kx}\Big(\frac{\zeta_y-\zeta_x}{2}\eta+\frac{\zeta_x-\zeta_y}{2}+1\Big)\phi_{ky}\Big(\frac{\zeta_y-\zeta_x}{2}\eta+\frac{\zeta_y-\zeta_x}{2}-1\Big)\;d\eta$ \\

     &$\chi^{m,r}_{1T,R}(\zeta_x,\zeta_y)$ & $\frac{\zeta_x}{2}\int_{-1}^1\phi_{kx}\Big(\frac{\zeta_x}{2}\eta+\frac{\zeta_x}{2}-1\Big)\phi_{ky}\Big(\frac{\zeta_x}{2}\eta+\zeta_y-1-\frac{\zeta_x}{2}\Big)\;d\eta$ \\ \hline\hline
      
     {\bf 1B} &$\chi^{m,r}_{1B,L}(\zeta_x,\zeta_y)$ & $\frac{2-\zeta_x}{2}\int_{-1}^1\phi_{kx}\Big(\big(1-\frac{\zeta_x}{2}\big)\eta+\frac{\zeta_x}{2}\Big)\phi_{ky}\Big(\zeta_y-\frac{\zeta_x}{2}+\big(1-\frac{\zeta_x}{2}\big)\eta\Big)\;d\eta$  \\
   
     &$\chi^{m,r}_{1B,M}(\zeta_x,\zeta_y)$ & $\frac{\zeta_x-\zeta_y}{2}\int_{-1}^1\phi_{kx}\Big(\frac{\zeta_x-\zeta_y}{2}\eta+\frac{\zeta_x-\zeta_y}{2}-1\Big)\phi_{ky}\Big(\frac{\zeta_x-\zeta_y}{2}\eta+\frac{\zeta_y-\zeta_x}{2}+1\Big)\;d\eta$ \\

       &$\chi^{m,r}_{1B,R}(\zeta_x,\zeta_y)$ & $\frac{\zeta_y}{2}\int_{-1}^1\phi_{kx}\Big(\frac{\zeta_y}{2}\eta+\zeta_x-\frac{\zeta_y}{2}-1\Big)\phi_{ky}\Big(\frac{\zeta_y}{2}\eta-1+\frac{\zeta_y}{2}\Big)\;d\eta$ \\ \hline\hline
      
       {\bf 2T} &$\chi^{m,r}_{2T,L}(\zeta_x,\zeta_y)$ & $-\frac{\zeta_x}{2}\int_{-1}^1\phi_{kx}\Big(-\frac{\zeta_x}{2}\eta+\frac{\zeta_x}{2}+1\Big)\phi_{ky}\Big(-\frac{\zeta_x}{2}\eta+\zeta_y-1-\frac{\zeta_x}{2}\Big)\;d\eta$  \\
   
     &$\chi^{m,r}_{2T,M}(\zeta_x,\zeta_y)$ & $\frac{2+\zeta_x-\zeta_y}{2}\int_{-1}^1\phi_{kx}\Big(\frac{2+\zeta_x-\zeta_y}{2}\eta+\frac{\zeta_x-\zeta_y}{2}\Big)\phi_{ky}\Big(\frac{2+\zeta_x-\zeta_y}{2}\eta+\frac{\zeta_y-\zeta_x}{2}\Big)\;d\eta$ \\

       &$\chi^{m,r}_{2T,R}(\zeta_x,\zeta_y)$ &  $\frac{\zeta_y}{2}\int_{-1}^1\phi_{kx}\Big(\frac{\zeta_y}{2}\eta+\zeta_x-\frac{\zeta_y}{2}+1\Big)\phi_{ky}\Big(\frac{\zeta_y}{2}\eta+\frac{\zeta_y}{2}-1\Big)\;d\eta$ \\ \hline\hline

      {\bf 3T} &$\chi^{m,r}_{3T,L}(\zeta_x,\zeta_y)$ & $-\frac{\zeta_y}{2}\int_{-1}^1\phi_{kx}\Big(-\frac{\zeta_y}{2}\eta-\frac{\zeta_y}{2}+\zeta_x+1\Big)\phi_{ky}\Big(-\frac{\zeta_y}{2}\eta+1+\frac{\zeta_y}{2}\Big)\;d\eta$  \\
   
     &$\chi^{m,r}_{3T,M}(\zeta_x,\zeta_y)$ & $\frac{\zeta_y-\zeta_x}{2}\int_{-1}^1\phi_{kx}\Big(\frac{\zeta_y-\zeta_x}{2}\eta+\frac{\zeta_x-\zeta_y}{2}+1\Big)\phi_{ky}\Big(\frac{\zeta_y-\zeta_x}{2}\eta+\frac{\zeta_y-\zeta_x}{2}-1\Big)\;d\eta$ \\

       &$\chi^{m,r}_{3T,R}(\zeta_x,\zeta_y)$ &  $\frac{2+\zeta_x}{2}\int_{-1}^1\phi_{kx}\Big(\big(\frac{2+\zeta_x}{2}\big)\eta+\frac{\zeta_x}{2}\Big)\phi_{ky}\Big(\big(\frac{2+\zeta_x}{2}\big)\eta-\frac{\zeta_x}{2}+\zeta_y\Big)\;d\eta$ \\ \hline\hline

       {\bf 3B} &$\chi^{m,r}_{3B,L}(\zeta_x,\zeta_y)$ & $-\frac{\zeta_x}{2}\int_{-1}^1\phi_{kx}\Big(-\frac{\zeta_x}{2}\eta+\frac{\zeta_x}{2}+1\Big)\phi_{ky}\Big(-\frac{\zeta_x}{2}\eta+1+\zeta_y-\frac{\zeta_x}{2}\Big)\;d\eta$  \\
   
        &$\chi^{m,r}_{3B,M}(\zeta_x,\zeta_y)$ &  $\frac{\zeta_x-\zeta_y}{2}\int_{-1}^1\phi_{kx}\Big(\frac{\zeta_x-\zeta_y}{2}\eta+\frac{\zeta_x-\zeta_y}{2}-1\Big)\phi_{ky}\Big(\frac{\zeta_x-\zeta_y}{2}\eta+\frac{\zeta_y-\zeta_x}{2}+1\Big)\;d\eta$ \\

        &$\chi^{m,r}_{3B,R}(\zeta_x,\zeta_y)$ &  $\frac{2+\zeta_y}{2}\int_{-1}^1\phi_{kx}\Big(\big(\frac{2+\zeta_y}{2}\big)\eta+\zeta_x-\frac{\zeta_y}{2}\Big)\phi_{ky}\Big(\big(\frac{2+\zeta_y}{2}\big)\eta+\frac{\zeta_y}{2}\Big)\;d\eta$ \\ \hline

     {\bf  4B} &$\chi^{m,r}_{4B,L}(\zeta_x,\zeta_y)$ & $-\frac{\zeta_y}{2}\int_{-1}^1\phi_{kx}\Big(-\frac{\zeta_y}{2}\eta-\frac{\zeta_y}{2}+\zeta_x-1\Big)\phi_{ky}\Big(-\frac{\zeta_y}{2}\eta+1+\frac{\zeta_y}{2}\Big)\;d\eta$  \\
   
       &$\chi^{m,r}_{4B,M}(\zeta_x,\zeta_y)$ & $\frac{\zeta_y-\zeta_x+2}{2}\int_{-1}^1\phi_{kx}\Big(\big(1+\frac{\zeta_y-\zeta_x}{2}\big)\eta+\frac{\zeta_x-\zeta_y}{2}\Big)\phi_{ky}\Big(\big(1+\frac{\zeta_y-\zeta_x}{2}\big)\eta+\frac{\zeta_y-\zeta_x}{2}\Big)\;d\eta$ \\

       &$\chi^{m,r}_{4B,R}(\zeta_x,\zeta_y)$ &  $\frac{\zeta_x}{2}\int_{-1}^1\phi_{kx}\Big(\frac{\zeta_x}{2}\eta+\frac{\zeta_x}{2}-1\Big)\phi_{ky}\Big(\frac{\zeta_x}{2}\eta+1+\zeta_y-\frac{\zeta_x}{2}\Big)\;d\eta$ \\ \hline

\end{tabular}

\end{table}
\end{landscape}

%% file: Supplementary_Material/2D_L-SIAC_coefficients.tex
  \newpage
  
  \begin{landscape}
  \begin{table}[h!]
      \centering
       
      \caption{Table describing how the shifting of the arguments of the $a$ coefficients vary by region.}
      \label{tab:2D_LSIAC_coefficient}
   \begin{tabular}{|M{7.5cm}||M{3cm}|M{3cm}|M{3cm}|M{3cm}||}\hline
{\bf Region $S$ $\Big{\backslash}$ Arguments of $a_{\alpha}(q_x,q_y)$} & ${\bf  P=L}$&${\bf P=M}$&${\bf P=R}$\\ \hline\hline
{\bf 1T}  &$(0,0)$   &  $(0,1)$&  $(1,1)$   \\ \hline
{\bf 1B}  &$(0,0)$   &  $(1,0)$&  $(1,1)$   \\ \hline
{\bf 2T}  &$(-1,0)$   &  $(0,0)$&  $(0,1)$   \\ \hline
{\bf 3T}  &$(-1,-1)$   &  $(-1,0)$&  $(0,0)$   \\ \hline
{\bf 3B}  &$(-1,-1)$   &  $(0,-1)$&  $(0,0)$   \\ \hline
{\bf 4B}  &$(0,-1)$   &  $(0,0)$&  $(1,0)$   \\ \hline
\end{tabular}

  \end{table}

\input{Supplementary_Material/basis_function_domains_LSIAC}

  \end{landscape}
  

%% file: Supplementary_Material/basis_function_domains_LSIAC.tex
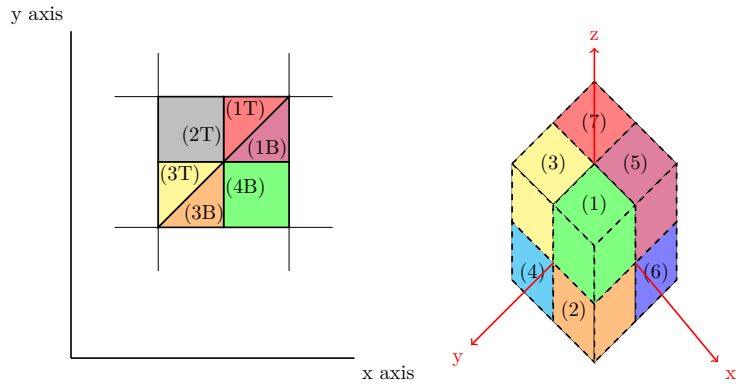
\begin{figure}[bp!]
\centering
    \begin{tabular}{c r}

  \begin{tikzpicture}[scale =0.75]
 \usetikzlibrary{calc}

 \coordinate (Origin) at (0,0);

\draw (0,1) -- (5,1);
\draw (0,4) -- (5,4);
\draw (1,0) -- (1,5);
\draw (4,0) -- (4,5);
\draw (1,1) -- (4,4);

\draw (1,2.5) -- (4,2.5);
\draw (2.5,1) -- (2.5,4);

\draw[thick] (-1,-2) -- (5.5,-2) node[anchor=north west] {x axis };
\draw[thick] (-1,-2) -- (-1,5.5) node[anchor=south east] {y axis};

\filldraw [thick, fill=red!50]  (2.5,2.5) -- (2.5,4) -- (4,4) -- cycle;
   \node at (3.0,3.7){(1T)};
\filldraw [thick, fill=gray!50]  (1,2.5) -- (2.5,2.5) -- (2.5,4) -- (1,4) -- cycle;
   \node at (2.0,3.1){(2T)};
\filldraw [thick, fill=yellow!50]  (1,2.5) -- (2.5,2.5) -- (1,1) -- cycle;
   \node at (1.5,2.2){(3T)};
   
   \filldraw [thick, fill=purple!50]  (2.5,2.5) -- (4,2.5) -- (4,4) -- cycle;
   \node at (3.5,2.8){(1B)};
\filldraw [thick, fill=green!50]  (2.5,1) -- (4,1) -- (4,2.5) -- (2.5,2.5) -- cycle;
   \node at (3.0,1.9){(4B)};
\filldraw [thick, fill=orange!50]  (1,1) -- (2.5,2.5) -- (2.5,1) -- cycle;
   \node at (2.05,1.3){(3B)};
\end{tikzpicture}
      
    
    &
       \begin{tikzpicture}
 \usetikzlibrary{calc}

 \coordinate (Origin) at (0,0);


\draw (0,0) -- ({2*cos(45)},{2*sin(45)}) -- ({4*cos(45)},0) -- ({4*cos(45)},-2) -- ({2*cos(45)},{-2-2*sin(45)}) -- (0,{-2}) -- cycle; 
\draw (0,0) -- ({2*cos(45)},{-2*sin(45)}) -- ({4*cos(45)},0);

\draw  ({2*cos(45)},{-2*sin(45)}) --  ({2*cos(45)},{-2-2*sin(45)});
\draw  ({cos(45)},{-sin(45)}) --  ({cos(45)},{-2-1*sin(45)});
\draw  ({3*cos(45)},{-sin(45)}) --  ({3*cos(45)},{-2-1*sin(45)});

\draw  (0,{-1}) --  ({2*cos(45)},{-1-2*sin(45)});
\draw  ({2*cos(45)},{-1-2*sin(45)}) --  ({4*cos(45)},{-1});

\draw  ({cos(45)},{-sin(45)}) --  ({3*cos(45)},{sin(45)});
\draw  ({cos(45)},{sin(45)}) --  ({3*cos(45)},{-sin(45)});

\filldraw [thick,dashed, fill=yellow!50]  (0,0) -- ({cos(45)},{sin(45)}) -- ({2*cos(45)},{0}) -- ({cos(45)},{-sin(45)}) -- cycle;
\filldraw [thick,dashed, fill=yellow!50]  (0,0) -- ({cos(45)},{-sin(45)}) -- ({cos(45)},{-1-sin(45)}) -- ({0},{-1}) -- cycle;
   \node at ({cos(45)},0){(3)};

\filldraw [thick,dashed, fill=green!50]  ({2*cos(45)},{0}) -- ({cos(45)},{-sin(45)}) -- ({2*cos(45)},{-2*sin(45)}) -- ({3*cos(45)},{-sin(45)}) -- cycle;
\filldraw [thick,dashed, fill=green!50]  ({cos(45)},{-sin(45)}) -- ({cos(45)},{-sin(45)-1}) -- ({2*cos(45)},{-1-2*sin(45)}) -- ({2*cos(45)},{-2*sin(45)}) -- cycle;
\filldraw [thick,dashed, fill=green!50]  ({2*cos(45)},{-2*sin(45)}) -- ({2*cos(45)},{-2*sin(45)-1}) -- ({3*cos(45)},{-1-sin(45)}) -- ({3*cos(45)},{-sin(45)}) -- cycle;
   \node at ({2*cos(45)},{-sin(45)}){(1)};

   \filldraw [thick,dashed, fill=purple!50]  ({2*cos(45)},{0}) -- ({3*cos(45)},{-sin(45)}) -- ({4*cos(45)},{0}) -- ({3*cos(45)},{sin(45)}) -- cycle;
\filldraw [thick,dashed, fill=purple!50]  ({3*cos(45)},{-sin(45)}) -- ({3*cos(45)},{-sin(45)-1}) -- ({4*cos(45)},{-1}) -- ({4*cos(45)},0) -- cycle;
   \node at ({3*cos(45)},{0}){(5)};

\filldraw [thick,dashed, fill=red!50]  ({cos(45)},{sin(45)}) -- ({2*cos(45)},{0}) -- ({3*cos(45)},{sin(45)}) -- ({2*cos(45)},{2*sin(45)}) -- cycle;
   \node at ({2*cos(45)},{sin(45)}){(7)};

\filldraw [thick,dashed, fill=cyan!50]  ({0},{-2}) -- ({cos(45)},{-2-sin(45)}) -- ({cos(45)},{-1-sin(45)}) -- ({0},{-1}) -- cycle;
   \node at ({0.5*cos(45)},{-1-1.25*sin(45)}){ (4)};

\filldraw [thick,dashed, fill=orange!50]  ({cos(45)},{-2-sin(45)}) -- ({cos(45)},{-1-sin(45)}) -- ({2*cos(45)},{-1-2*sin(45)}) -- ({2*cos(45)},{-2-2*sin(45)}) -- cycle;
\filldraw [thick,dashed, fill=orange!50]  ({2*cos(45)},{-1-2*sin(45)}) -- ({2*cos(45)},{-2-2*sin(45)}) -- ({3*cos(45)},{-2-sin(45)}) -- ({3*cos(45)},{-1-sin(45)}) -- cycle; 
   \node at ({1.5*cos(45)},{-2-0.75*sin(45)}){ (2)};

\filldraw [thick,dashed, fill=blue!50]  ({3*cos(45)},{-2-sin(45)}) -- ({3*cos(45)},{-1-sin(45)}) -- ({4*cos(45)},{-1}) -- ({4*cos(45)},{-2}) -- cycle; 
   \node at ({3.5*cos(45)},{-1-1.25*sin(45)}){ (6)};

\draw[thick,red,->] ({2*cos(45)},0) -- ({2*cos(45)},2) node[anchor=south] {z};
\draw[thick,red,->] ({3*cos(45)},{-1-sin(45)}) -- ({5*cos(45)},{-2-2*sin(45)}) node[anchor= north west] {x};
\draw[thick,red,->] ({cos(45)},{-1-sin(45)}) -- ({-cos(45)},{-1-3*sin(45)}) node[anchor= north east] {y};

\end{tikzpicture}
      \label{fig:Octant_Ordering}

    \end{tabular}
      \caption{Partition of domain for both  2D L-SIAC (left) and the octant numbering scheme in 3D (right).}
         \label{fig:2D_basis_function_domain}
 \end{figure}
 

%% file: Supplementary_Material/2D_LSIAC_algorithm.tex
\begin{algorithm}
\caption{Calculate $u^f$ from $u$ using L-SIAC kernel}
\label{alg:LSIAC}

\begin{algorithmic}[1]

 \State \textbf{Input:}\text{ Number of elements in one direction $N$, and polynomial degree $p$.}\;
    
  \State \textbf{Output:}\text{ Fine mesh modal coefficients $u^f$}\;\\
  
  \State \%  \text{Compute modal projection on coarse mesh}
  \State  $u=\text{Initialize}(N,p,\Omega,u_0)$\;     \\

\State   \% \text{Compute projection matrices by region}\;
  \State  $\hphantom{[P1BL,P1BM,P1BR]}\mathllap{[P1TL,P1TM,P1TR]}\phantom{]}=\text{ProjMatrixL}(N,p,Q1T)$\;
  \State  $\hphantom{[[P1BL,P1BM,P1BR]}\mathllap{[P1BL,P1BM,P1BR]}=\text{ProjMatrixL}(N,p,Q1B)$\;
  \State  $\hphantom{[P1BL,P1BM,P1BR]}\mathllap{[P2TL,P2TM,P2TR]}\phantom{]}=\text{ProjMatrixL}(N,p,Q2T)$\;
 \State   $\hphantom{[P1BL,P1BM,P1BR]}\mathllap{[P3TL,P3TM,P3TR]}\phantom{]}=\text{ProjMatrixL}(N,p,Q3T)$\;
  \State  $\hphantom{[[P1BL,P1BM,P1BR]}\mathllap{[P3BL,P3BM,P3BR]}=\text{ProjMatrixL}(N,p,Q3B)$\;
  \State  $\hphantom{[[P1BL,P1BM,P1BR]}\mathllap{[P4BL,P4BM,P4BR]}=\text{ProjMatrixL}(N,p,Q4B)$\;\\
    
    \State \% \text{Loop through mesh elements}\;
\For{$i\gets 1, N$}
    \For{$j\gets 1, N$}   \\

   \State \% \text{Construct filtering/sifting matrices}\;
    
 \State   $\hphantom{Kn1n1}\mathllap{K11}\phantom{nn}=\text{FiltMatrixL}(p,1,1,i,j)$\;
\State    $\hphantom{Kn1n1}\mathllap{K01}\phantom{nn}=\text{FiltMatrixL}(p,0,1,i,j)$\;
 \State   $\hphantom{Kn1n1}\mathllap{K10}\phantom{nn}=\text{FiltMatrixL}(p,1,0,i,j)$\; 
 \State   $\hphantom{Kn1n1}\mathllap{K00}\phantom{nn}=\text{FiltMatrixL}(p,0,0,i,j)$\; 
 \State   $\hphantom{nKn1n1}\mathllap{Kn10}\phantom{n}=\text{FiltMatrixL}(p,-1,0,i,j)$\; 
\State    $\hphantom{nKn1n1}\mathllap{K0n1}\phantom{n}=\text{FiltMatrixL}(p,0,-1,i,j)$\; 
 \State   $\hphantom{nnKn1n1}\mathllap{Kn1n1}=\text{FiltMatrixL}(p,-1,-1,i,j)$\; \\
    
     \State \% \text{Compute local fine mesh coefficients}\;
   \State $u^f_{tensor}(:,2i,2j)=$
   \newline \hspace*{7em}$((P1TL+P1BL)*K00+$
   \newline \hspace*{7em} $(P1TR+P1BR)*K11+$
   \newline \hspace*{7em}$(P1TM*K01+P1BM*K10)*u$\; 

  \State  $u^f_{tensor}(:,2i-1,2j)=$
  \newline \hspace*{7em}$(P2L*Kn10+P2M*K00+P2R*K01)*u$\;
    
 \State   $u^f_{tensor}(:,2i-1,2j-1)=$
 \newline \hspace*{7em}$((P3TL+P3BL)*Kn1n1+$
  \newline \hspace*{7em}$(P3TR+P3BR)*K00+$
 \newline \hspace*{7em} $P3TM*Kn10+P3BM*K0n1)*u$\; 

  \State  $u^f_{tensor}(:,2i,2j-1)=$
  \newline \hspace*{7em}$(P4L*K0n1+P4M*K00+P4R*K10)*u$\;
        
    \EndFor
\EndFor
\\
      \State \% \text{Reshape fine modal array into original array form}\;
   \State  $u^f=\text{ReshapeArray}(p,N,u^f_{tensor})$\;
          \end{algorithmic}

\end{algorithm}

%% file: Supplementary_Material/3D_L-SIAC_Classification_Scheme.tex
\begin{table}[hp!]
\begin{minipage}[b]{.45\textwidth}
  \centering
          \caption{ Regions present by octant}
          \label{tab:3D_Regions_Present} 
      \begin{tabular}{|c|c|} \hline
           \textbf{Octant}&\textbf{Regions Present} \\ \hline
          \textbf{1} & 1A, 1B, 1C, 1D, 1E, 1F\\ \hline
          \textbf{2} & 2E, 2F\\ \hline
          \textbf{3} & 3A, 3B\\ \hline
          \textbf{4} & 4B, 4F\\ \hline
          \textbf{5}& 5C, 5D\\ \hline
          \textbf{6} & 6D, 6E\\ \hline
          \textbf{7} & 7A, 7C\\ \hline
          \textbf{8} & 8A, 8B, 8C, 8D, 8E, 8F\\ \hline
      \end{tabular}

\end{minipage}\qquad
\begin{minipage}[b]{0.45\textwidth}
  \centering
  \caption{Classification of regions.}
  \label{tab:3D_Classification}
      \begin{tabular}{|c|c|} \hline
           \textbf{Class}&\textbf{Filtering Coordinate Inequality} \\ \hline
          \textbf{(A)} & $\zeta_x\leq\zeta_y\leq\zeta_z$\\ \hline
          \textbf{(B)} & $\zeta_x\leq\zeta_z<\zeta_y$\\ \hline
          \textbf{(C)} & $\zeta_y<\zeta_x\leq\zeta_z$\\ \hline
          \textbf{(D)} & $\zeta_y\leq\zeta_z<\zeta_x$\\ \hline
          \textbf{(E)} & $\zeta_z<\zeta_y\leq\zeta_x$\\ \hline
          \textbf{(F)} & $\zeta_z<\zeta_x<\zeta_y$\\ \hline
      \end{tabular}
\end{minipage}
\end{table}

%% file: Supplementary_Material/3D_L-SIAC_Basis_Table.tex
\begin{landscape}
\begin{table}[hp!]

\caption{Three-dimensional L-SIAC basis functions. For notational convenience suppose that $\mathbf{x}$ belong to element $(i,j,k)$ and that the region $S$ and position $P$ of the function are as indicated in the table rows.  }
\label{tab:3D_LSIAC_Basis}
\begin{adjustbox}{width=1.0\linewidth}\centering
\begin{tabular}{||c|c|c|c|c|c||} \hline 
\textbf{Region $S$} & \multicolumn{5}{c||}{\textbf{\textbf{Basis Function} Parameters}} \\ \hline 
	  & \textbf{Basis Function} & \textbf{Coefficient} & $\textbf{\text{arg}}_x$ & $\textbf{\text{arg}}_y$&$\textbf{\text{arg}}_z$ \\ 
	\hline
\textbf{1A}&$\chi^{S,L}$&$(2 - \zeta_z)/4$	&	$w + \zeta_x - \zeta_z/2 - (w*\zeta_z)/2$	&	$w + \zeta_y - \zeta_z/2 - (w*\zeta_z)/2$	&	$w - ((-1 + w)*\zeta_z)/2$\\
	\hline
\textbf{1A}&$\chi^{S,LM}$&$(-\zeta_y + \zeta_z)/4$	&	$(2 + 2*\zeta_x - (1 + w)*\zeta_y - \zeta_z + w*\zeta_z)/2$	&	$(2 + \zeta_y - w*\zeta_y + (-1 + w)*\zeta_z)/2$	&	$(-2 - (1 + w)*\zeta_y + \zeta_z + w*\zeta_z)/2$\\
	\hline
\textbf{1A}&$\chi^{S,RM}$&$(-\zeta_x + \zeta_y)/4$	&	$(2 + \zeta_x - w*\zeta_x + (-1 + w)*\zeta_y)/2$	&	$(-2 - (1 + w)*\zeta_x + \zeta_y + w*\zeta_y)/2$	&	$-1 - ((1 + w)*\zeta_x)/2 + ((-1 + w)*\zeta_y)/2 + \zeta_z$\\
	\hline
\textbf{1A}&$\chi^{S,R}$&$\zeta_x/4$	&	$(-2 + \zeta_x + w*\zeta_x)/2$	&	$-1 + ((-1 + w)*\zeta_x)/2 + \zeta_y$	&	$-1 + ((-1 + w)*\zeta_x)/2 + \zeta_z$\\
	\hline \hline
\textbf{Region} & \multicolumn{5}{c||}{\textbf{\textbf{Basis Function} Parameters}} \\ \hline 
	
	  & \textbf{Basis Function} & \textbf{Coefficient} & $\textbf{\text{arg}}_x$ & $\textbf{\text{arg}}_y$&$\textbf{\text{arg}}_z$ \\ 
	\hline
\textbf{1B}&$\chi^{S,L}$&$(2 - \zeta_y)/4$	&	$w + \zeta_x - \zeta_y/2 - (w*\zeta_y)/2$	&	$(-(w*(-2 + \zeta_y)) + \zeta_y)/2$	&	$w - \zeta_y/2 - (w*\zeta_y)/2 + \zeta_z$\\
	\hline
\textbf{1B}&$\chi^{S,LM}$&$(\zeta_y - \zeta_z)/4$	&	$(2 + 2*\zeta_x + (-1 + w)*\zeta_y - \zeta_z - w*\zeta_z)/2$	&	$(-2 + (1 + w)*\zeta_y - (1 + w)*\zeta_z)/2$	&	$(2 + (-1 + w)*\zeta_y + \zeta_z - w*\zeta_z)/2$\\
	\hline
\textbf{1B}&$\chi^{S,RM}$&$(-\zeta_x + \zeta_z)/4$	&	$(2 + \zeta_x - w*\zeta_x + (-1 + w)*\zeta_z)/2$	&	$(-2 - (1 + w)*\zeta_x + 2*\zeta_y - \zeta_z + w*\zeta_z)/2$	&	$(-2 - (1 + w)*\zeta_x + \zeta_z + w*\zeta_z)/2$\\
	\hline
\textbf{1B}&$\chi^{S,R}$&$\zeta_x/4$	&	$(-2 + \zeta_x + w*\zeta_x)/2$	&	$-1 + ((-1 + w)*\zeta_x)/2 + \zeta_y$	&	$-1 + ((-1 + w)*\zeta_x)/2 + \zeta_z$\\
	\hline 	\hline
\textbf{Region} & \multicolumn{5}{c||}{\textbf{\textbf{Basis Function} Parameters}} \\ \hline 

	  & \textbf{Basis Function} & \textbf{Coefficient} & $\textbf{\text{arg}}_x$ & $\textbf{\text{arg}}_y$&$\textbf{\text{arg}}_z$ \\ 
	\hline
\textbf{1C}&$\chi^{S,L}$&$(2 - \zeta_z)/4$	&	$w + \zeta_x - \zeta_z/2 - (w*\zeta_z)/2$	&	$w + \zeta_y - \zeta_z/2 - (w*\zeta_z)/2$	&	$w - ((-1 + w)*\zeta_z)/2$\\
	\hline
\textbf{1C}&$\chi^{S,LM}$&$(-\zeta_x + \zeta_z)/4$	&	$(2 + \zeta_x - w*\zeta_x + (-1 + w)*\zeta_z)/2$	&	$(2 - (1 + w)*\zeta_x + 2*\zeta_y - \zeta_z + w*\zeta_z)/2$	&	$(-2 - (1 + w)*\zeta_x + \zeta_z + w*\zeta_z)/2$\\
	\hline
\textbf{1C}&$\chi^{S,RM}$&$(\zeta_x - \zeta_y)/4$	&	$(-2 + (1 + w)*\zeta_x - (1 + w)*\zeta_y)/2$	&	$(2 + (-1 + w)*\zeta_x + \zeta_y - w*\zeta_y)/2$	&	$-1 + ((-1 + w)*\zeta_x)/2 - ((1 + w)*\zeta_y)/2 + \zeta_z$\\
	\hline
\textbf{1C}&$\chi^{S,R}$&$\zeta_y/4$	&	$-1 + \zeta_x + ((-1 + w)*\zeta_y)/2$	&	$(-2 + \zeta_y + w*\zeta_y)/2$	&	$-1 + ((-1 + w)*\zeta_y)/2 + \zeta_z$\\
	\hline 	\hline
\textbf{Region} & \multicolumn{5}{c||}{\textbf{\textbf{Basis Function} Parameters}} \\ \hline 

	  & \textbf{Basis Function} & \textbf{Coefficient} & $\textbf{\text{arg}}_x$ & $\textbf{\text{arg}}_y$&$\textbf{\text{arg}}_z$ \\ 
	\hline
\textbf{1D}&$\chi^{S,L}$&$(2 - \zeta_x)/4$	&	$(-(w*(-2 + \zeta_x)) + \zeta_x)/2$	&	$w - \zeta_x/2 - (w*\zeta_x)/2 + \zeta_y$	&	$w - \zeta_x/2 - (w*\zeta_x)/2 + \zeta_z$\\
	\hline
\textbf{1D}&$\chi^{S,LM}$&$(\zeta_x - \zeta_z)/4$	&	$(-2 + (1 + w)*\zeta_x - (1 + w)*\zeta_z)/2$	&	$(2 + (-1 + w)*\zeta_x + 2*\zeta_y - \zeta_z - w*\zeta_z)/2$	&	$(2 + (-1 + w)*\zeta_x + \zeta_z - w*\zeta_z)/2$\\
	\hline
\textbf{1D}&$\chi^{S,RM}$&$(-\zeta_y + \zeta_z)/4$	&	$(-2 + 2*\zeta_x - (1 + w)*\zeta_y - \zeta_z + w*\zeta_z)/2$	&	$(2 + \zeta_y - w*\zeta_y + (-1 + w)*\zeta_z)/2$	&	$(-2 - (1 + w)*\zeta_y + \zeta_z + w*\zeta_z)/2$\\
	\hline
\textbf{1D}&$\chi^{S,R}$&$\zeta_y/4$	&	$-1 + \zeta_x + ((-1 + w)*\zeta_y)/2$	&	$(-2 + \zeta_y + w*\zeta_y)/2$	&	$-1 + ((-1 + w)*\zeta_y)/2 + \zeta_z$\\
	\hline 	\hline
\textbf{Region} & \multicolumn{5}{c||}{\textbf{\textbf{Basis Function} Parameters}} \\ \hline 

	  & \textbf{Basis Function} & \textbf{Coefficient} & $\textbf{\text{arg}}_x$ & $\textbf{\text{arg}}_y$&$\textbf{\text{arg}}_z$ \\ 
	\hline
\textbf{1E}&$\chi^{S,L}$&$(2 - \zeta_x)/4$	&	$(-(w*(-2 + \zeta_x)) + \zeta_x)/2$	&	$w - \zeta_x/2 - (w*\zeta_x)/2 + \zeta_y$	&	$w - \zeta_x/2 - (w*\zeta_x)/2 + \zeta_z$\\
	\hline
\textbf{1E}&$\chi^{S,LM}$&$(\zeta_x - \zeta_y)/4$	&	$(-2 + (1 + w)*\zeta_x - (1 + w)*\zeta_y)/2$	&	$(2 + (-1 + w)*\zeta_x + \zeta_y - w*\zeta_y)/2$	&	$1 + ((-1 + w)*\zeta_x)/2 - ((1 + w)*\zeta_y)/2 + \zeta_z$\\
	\hline
\textbf{1E}&$\chi^{S,RM}$&$(\zeta_y - \zeta_z)/4$	&	$(-2 + 2*\zeta_x + (-1 + w)*\zeta_y - \zeta_z - w*\zeta_z)/2$	&	$(-2 + (1 + w)*\zeta_y - (1 + w)*\zeta_z)/2$	&	$(2 + (-1 + w)*\zeta_y + \zeta_z - w*\zeta_z)/2$\\
	\hline
\textbf{1E}&$\chi^{S,R}$&$\zeta_z/4$	&	$-1 + \zeta_x + ((-1 + w)*\zeta_z)/2$	&	$-1 + \zeta_y + ((-1 + w)*\zeta_z)/2$	&	$(-2 + \zeta_z + w*\zeta_z)/2$\\
	\hline 	\hline
\textbf{Region} & \multicolumn{5}{c||}{\textbf{\textbf{Basis Function} Parameters}} \\ \hline 

	  & \textbf{Basis Function} & \textbf{Coefficient} & $\textbf{\text{arg}}_x$ & $\textbf{\text{arg}}_y$&$\textbf{\text{arg}}_z$ \\ 
	\hline
\textbf{1F}&$\chi^{S,L}$&$(2 - \zeta_y)/4$	&	$w + \zeta_x - \zeta_y/2 - (w*\zeta_y)/2$	&	$(-(w*(-2 + \zeta_y)) + \zeta_y)/2$	&	$w - \zeta_y/2 - (w*\zeta_y)/2 + \zeta_z$\\
	\hline
\textbf{1F}&$\chi^{S,LM}$&$(-\zeta_x + \zeta_y)/4$	&	$(2 + \zeta_x - w*\zeta_x + (-1 + w)*\zeta_y)/2$	&	$(-2 - (1 + w)*\zeta_x + \zeta_y + w*\zeta_y)/2$	&	$1 - ((1 + w)*\zeta_x)/2 + ((-1 + w)*\zeta_y)/2 + \zeta_z$\\
	\hline
\textbf{1F}&$\chi^{S,RM}$&$(\zeta_x - \zeta_z)/4$	&	$(-2 + (1 + w)*\zeta_x - (1 + w)*\zeta_z)/2$	&	$(-2 + (-1 + w)*\zeta_x + 2*\zeta_y - \zeta_z - w*\zeta_z)/2$	&	$(2 + (-1 + w)*\zeta_x + \zeta_z - w*\zeta_z)/2$\\
	\hline
\textbf{1F}&$\chi^{S,R}$&$\zeta_z/4$	&	$-1 + \zeta_x + ((-1 + w)*\zeta_z)/2$	&	$-1 + \zeta_y + ((-1 + w)*\zeta_z)/2$	&	$(-2 + \zeta_z + w*\zeta_z)/2$\\
	\hline
\end{tabular} 
\end{adjustbox}
\end{table}
\end{landscape}
\newpage
\begin{landscape}
    \centering
\begin{adjustbox}{width=1.0\linewidth}\centering
\begin{tabular}{||c|c|c|c|c|c||} \hline 
\textbf{Region} & \multicolumn{5}{c||}{\textbf{\textbf{Basis Function} Parameters}} \\ \hline 
	
	  & \textbf{Basis Function} & \textbf{Coefficient} & $\textbf{\text{arg}}_x$ & $\textbf{\text{arg}}_y$&$\textbf{\text{arg}}_z$ \\ 
	\hline
\textbf{2E}&$\chi^{S,L}$&$-\zeta_z/4$	&	$-1 + \zeta_x - ((1 + w)*\zeta_z)/2$	&	$-1 + \zeta_y - ((1 + w)*\zeta_z)/2$	&	$(2 + \zeta_z - w*\zeta_z)/2$\\
	\hline
\textbf{2E}&$\chi^{S,LM}$&$(2 - \zeta_x + \zeta_z)/4$	&	$(\zeta_x - \zeta_z + w*(2 - \zeta_x + \zeta_z))/2$	&	$(-\zeta_x + 2*\zeta_y - \zeta_z + w*(2 - \zeta_x + \zeta_z))/2$	&	$(-\zeta_x + \zeta_z + w*(2 - \zeta_x + \zeta_z))/2$\\
	\hline
\textbf{2E}&$\chi^{S,RM}$&$(\zeta_x - \zeta_y)/4$	&	$(-2 + (1 + w)*\zeta_x - (1 + w)*\zeta_y)/2$	&	$(2 + (-1 + w)*\zeta_x + \zeta_y - w*\zeta_y)/2$	&	$1 + ((-1 + w)*\zeta_x)/2 - ((1 + w)*\zeta_y)/2 + \zeta_z$\\
	\hline
\textbf{2E}&$\chi^{S,R}$&$\zeta_y/4$	&	$-1 + \zeta_x + ((-1 + w)*\zeta_y)/2$	&	$(-2 + \zeta_y + w*\zeta_y)/2$	&	$1 + ((-1 + w)*\zeta_y)/2 + \zeta_z$\\
	\hline 	\hline

\textbf{Region} & \multicolumn{5}{c||}{\textbf{\textbf{Basis Function} Parameters}} \\ \hline 
	  & \textbf{Basis Function} & \textbf{Coefficient} & $\textbf{\text{arg}}_x$ & $\textbf{\text{arg}}_y$&$\textbf{\text{arg}}_z$ \\ 
	\hline
\textbf{2F}&$\chi^{S,L}$&$-\zeta_z/4$	&	$-1 + \zeta_x - ((1 + w)*\zeta_z)/2$	&	$-1 + \zeta_y - ((1 + w)*\zeta_z)/2$	&	$(2 + \zeta_z - w*\zeta_z)/2$\\
	\hline
\textbf{2F}&$\chi^{S,LM}$&$(2 - \zeta_y + \zeta_z)/4$	&	$(2*\zeta_x - \zeta_y - \zeta_z + w*(2 - \zeta_y + \zeta_z))/2$	&	$(\zeta_y - \zeta_z + w*(2 - \zeta_y + \zeta_z))/2$	&	$(-\zeta_y + \zeta_z + w*(2 - \zeta_y + \zeta_z))/2$\\
	\hline
\textbf{2F}&$\chi^{S,RM}$&$(-\zeta_x + \zeta_y)/4$	&	$(2 + \zeta_x - w*\zeta_x + (-1 + w)*\zeta_y)/2$	&	$(-2 - (1 + w)*\zeta_x + \zeta_y + w*\zeta_y)/2$	&	$1 - ((1 + w)*\zeta_x)/2 + ((-1 + w)*\zeta_y)/2 + \zeta_z$\\
	\hline
\textbf{2F}&$\chi^{S,R}$&$\zeta_x/4$	&	$(-2 + \zeta_x + w*\zeta_x)/2$	&	$-1 + ((-1 + w)*\zeta_x)/2 + \zeta_y$	&	$1 + ((-1 + w)*\zeta_x)/2 + \zeta_z$\\
	\hline \hline
\textbf{Region} & \multicolumn{5}{c||}{\textbf{\textbf{Basis Function} Parameters}} \\ \hline 
	
	  & \textbf{Basis Function} & \textbf{Coefficient} & $\textbf{\text{arg}}_x$ & $\textbf{\text{arg}}_y$&$\textbf{\text{arg}}_z$ \\ 
	\hline
\textbf{3A}&$\chi^{S,L}$&$-\zeta_x/4$	&	$(2 + \zeta_x - w*\zeta_x)/2$	&	$-1 - ((1 + w)*\zeta_x)/2 + \zeta_y$	&	$-1 - ((1 + w)*\zeta_x)/2 + \zeta_z$\\
	\hline
\textbf{3A}&$\chi^{S,LM}$&$(2 + \zeta_x - \zeta_z)/4$	&	$(\zeta_x + w*(2 + \zeta_x - \zeta_z) - \zeta_z)/2$	&	$(-\zeta_x + 2*\zeta_y + w*(2 + \zeta_x - \zeta_z) - \zeta_z)/2$	&	$(-\zeta_x + w*(2 + \zeta_x - \zeta_z) + \zeta_z)/2$\\
	\hline
\textbf{3A}&$\chi^{S,RM}$&$(-\zeta_y + \zeta_z)/4$	&	$(2 + 2*\zeta_x - (1 + w)*\zeta_y - \zeta_z + w*\zeta_z)/2$	&	$(2 + \zeta_y - w*\zeta_y + (-1 + w)*\zeta_z)/2$	&	$(-2 - (1 + w)*\zeta_y + \zeta_z + w*\zeta_z)/2$\\
	\hline
\textbf{3A}&$\chi^{S,R}$&$\zeta_y/4$	&	$1 + \zeta_x + ((-1 + w)*\zeta_y)/2$	&	$(-2 + \zeta_y + w*\zeta_y)/2$	&	$-1 + ((-1 + w)*\zeta_y)/2 + \zeta_z$\\
	\hline \hline
\textbf{Region} & \multicolumn{5}{c||}{\textbf{\textbf{Basis Function} Parameters}} \\ \hline 
	
	  & \textbf{Basis Function} & \textbf{Coefficient} & $\textbf{\text{arg}}_x$ & $\textbf{\text{arg}}_y$&$\textbf{\text{arg}}_z$ \\ 
	\hline
\textbf{3B}&$\chi^{S,L}$&$-\zeta_x/4$	&	$(2 + \zeta_x - w*\zeta_x)/2$	&	$-1 - ((1 + w)*\zeta_x)/2 + \zeta_y$	&	$-1 - ((1 + w)*\zeta_x)/2 + \zeta_z$\\
	\hline
\textbf{3B}&$\chi^{S,LM}$&$(2 + \zeta_x - \zeta_y)/4$	&	$(\zeta_x + w*(2 + \zeta_x - \zeta_y) - \zeta_y)/2$	&	$(-\zeta_x + w*(2 + \zeta_x - \zeta_y) + \zeta_y)/2$	&	$(-\zeta_x + w*(2 + \zeta_x - \zeta_y) - \zeta_y + 2*\zeta_z)/2$\\
	\hline
\textbf{3B}&$\chi^{S,RM}$&$(\zeta_y - \zeta_z)/4$	&	$(2 + 2*\zeta_x + (-1 + w)*\zeta_y - \zeta_z - w*\zeta_z)/2$	&	$(-2 + (1 + w)*\zeta_y - (1 + w)*\zeta_z)/2$	&	$(2 + (-1 + w)*\zeta_y + \zeta_z - w*\zeta_z)/2$\\
	\hline
\textbf{3B}&$\chi^{S,R}$&$\zeta_z/4$	&	$1 + \zeta_x + ((-1 + w)*\zeta_z)/2$	&	$-1 + \zeta_y + ((-1 + w)*\zeta_z)/2$	&	$(-2 + \zeta_z + w*\zeta_z)/2$\\
	\hline \hline
\textbf{Region} & \multicolumn{5}{c||}{\textbf{\textbf{Basis Function} Parameters}} \\ \hline 
	
	  & \textbf{Basis Function} & \textbf{Coefficient} & $\textbf{\text{arg}}_x$ & $\textbf{\text{arg}}_y$&$\textbf{\text{arg}}_z$ \\ 
	\hline
\textbf{4B}&$\chi^{S,L}$&$-\zeta_z/4$	&	$1 + \zeta_x - ((1 + w)*\zeta_z)/2$	&	$-1 + \zeta_y - ((1 + w)*\zeta_z)/2$	&	$(2 + \zeta_z - w*\zeta_z)/2$\\
	\hline
\textbf{4B}&$\chi^{S,LM}$&$(-\zeta_x + \zeta_z)/4$	&	$(2 + \zeta_x - w*\zeta_x + (-1 + w)*\zeta_z)/2$	&	$(-2 - (1 + w)*\zeta_x + 2*\zeta_y - \zeta_z + w*\zeta_z)/2$	&	$(-2 - (1 + w)*\zeta_x + \zeta_z + w*\zeta_z)/2$\\
	\hline
\textbf{4B}&$\chi^{S,RM}$&$(2 + \zeta_x - \zeta_y)/4$	&	$(\zeta_x + w*(2 + \zeta_x - \zeta_y) - \zeta_y)/2$	&	$(-\zeta_x + w*(2 + \zeta_x - \zeta_y) + \zeta_y)/2$	&	$(-\zeta_x + w*(2 + \zeta_x - \zeta_y) - \zeta_y + 2*\zeta_z)/2$\\
	\hline
\textbf{4B}&$\chi^{S,R}$&$\zeta_y/4$	&	$1 + \zeta_x + ((-1 + w)*\zeta_y)/2$	&	$(-2 + \zeta_y + w*\zeta_y)/2$	&	$1 + ((-1 + w)*\zeta_y)/2 + \zeta_z$\\
	\hline
\textbf{Region} & \multicolumn{5}{c||}{\textbf{\textbf{Basis Function} Parameters}} \\ \hline 
	\hline
	  & \textbf{Basis Function} & \textbf{Coefficient} & $\textbf{\text{arg}}_x$ & $\textbf{\text{arg}}_y$&$\textbf{\text{arg}}_z$ \\ 
	\hline
\textbf{4F}&$\chi^{S,L}$&$-\zeta_x/4$	&	$(2 + \zeta_x - w*\zeta_x)/2$	&	$-1 - ((1 + w)*\zeta_x)/2 + \zeta_y$	&	$1 - ((1 + w)*\zeta_x)/2 + \zeta_z$\\
	\hline
\textbf{4F}&$\chi^{S,LM}$&$(\zeta_x - \zeta_z)/4$	&	$(-2 + (1 + w)*\zeta_x - (1 + w)*\zeta_z)/2$	&	$(-2 + (-1 + w)*\zeta_x + 2*\zeta_y - \zeta_z - w*\zeta_z)/2$	&	$(2 + (-1 + w)*\zeta_x + \zeta_z - w*\zeta_z)/2$\\
	\hline
\textbf{4F}&$\chi^{S,RM}$&$(2 - \zeta_y + \zeta_z)/4$	&	$(2*\zeta_x - \zeta_y - \zeta_z + w*(2 - \zeta_y + \zeta_z))/2$	&	$(\zeta_y - \zeta_z + w*(2 - \zeta_y + \zeta_z))/2$	&	$(-\zeta_y + \zeta_z + w*(2 - \zeta_y + \zeta_z))/2$\\
	\hline
\textbf{4F}&$\chi^{S,R}$&$\zeta_y/4$	&	$1 + \zeta_x + ((-1 + w)*\zeta_y)/2$	&	$(-2 + \zeta_y + w*\zeta_y)/2$	&	$1 + ((-1 + w)*\zeta_y)/2 + \zeta_z$\\
	\hline

\end{tabular} 
\end{adjustbox}

\end{landscape}

\newpage
\begin{landscape}

\centering
 \begin{adjustbox}{width=1.0\linewidth}\centering
\begin{tabular}{||c|c|c|c|c|c||} \hline 
\textbf{Region} & \multicolumn{5}{c||}{\textbf{\textbf{Basis Function} Parameters}} \\ \hline 
	  & \textbf{Basis Function} & \textbf{Coefficient} & $\textbf{\text{arg}}_x$ & $\textbf{\text{arg}}_y$&$\textbf{\text{arg}}_z$ \\ 
	\hline
\textbf{5C}&$\chi^{S,L}$&$-\zeta_y/4$	&	$-1 + \zeta_x - ((1 + w)*\zeta_y)/2$	&	$(2 + \zeta_y - w*\zeta_y)/2$	&	$-1 - ((1 + w)*\zeta_y)/2 + \zeta_z$\\
	\hline
\textbf{5C}&$\chi^{S,LM}$&$(2 + \zeta_y - \zeta_z)/4$	&	$(2*\zeta_x - \zeta_y + w*(2 + \zeta_y - \zeta_z) - \zeta_z)/2$	&	$(\zeta_y + w*(2 + \zeta_y - \zeta_z) - \zeta_z)/2$	&	$(-\zeta_y + w*(2 + \zeta_y - \zeta_z) + \zeta_z)/2$\\
	\hline
\textbf{5C}&$\chi^{S,RM}$&$(-\zeta_x + \zeta_z)/4$	&	$(2 + \zeta_x - w*\zeta_x + (-1 + w)*\zeta_z)/2$	&	$(2 - (1 + w)*\zeta_x + 2*\zeta_y - \zeta_z + w*\zeta_z)/2$	&	$(-2 - (1 + w)*\zeta_x + \zeta_z + w*\zeta_z)/2$\\
	\hline
\textbf{5C}&$\chi^{S,R}$&$\zeta_x/4$	&	$(-2 + \zeta_x + w*\zeta_x)/2$	&	$1 + ((-1 + w)*\zeta_x)/2 + \zeta_y$	&	$-1 + ((-1 + w)*\zeta_x)/2 + \zeta_z$\\
	\hline 	\hline
\textbf{Region} & \multicolumn{5}{c||}{\textbf{\textbf{Basis Function} Parameters}} \\ \hline 

	  & \textbf{Basis Function} & \textbf{Coefficient} & $\textbf{\text{arg}}_x$ & $\textbf{\text{arg}}_y$&$\textbf{\text{arg}}_z$ \\ 
	\hline
\textbf{5D}&$\chi^{S,L}$&$-\zeta_y/4$	&	$-1 + \zeta_x - ((1 + w)*\zeta_y)/2$	&	$(2 + \zeta_y - w*\zeta_y)/2$	&	$-1 - ((1 + w)*\zeta_y)/2 + \zeta_z$\\
	\hline
\textbf{5D}&$\chi^{S,LM}$&$(2 - \zeta_x + \zeta_y)/4$	&	$(\zeta_x - \zeta_y + w*(2 - \zeta_x + \zeta_y))/2$	&	$(-\zeta_x + \zeta_y + w*(2 - \zeta_x + \zeta_y))/2$	&	$(-\zeta_x - \zeta_y + w*(2 - \zeta_x + \zeta_y) + 2*\zeta_z)/2$\\
	\hline
\textbf{5D}&$\chi^{S,RM}$&$(\zeta_x - \zeta_z)/4$	&	$(-2 + (1 + w)*\zeta_x - (1 + w)*\zeta_z)/2$	&	$(2 + (-1 + w)*\zeta_x + 2*\zeta_y - \zeta_z - w*\zeta_z)/2$	&	$(2 + (-1 + w)*\zeta_x + \zeta_z - w*\zeta_z)/2$\\
	\hline
\textbf{5D}&$\chi^{S,R}$&$\zeta_z/4$	&	$-1 + \zeta_x + ((-1 + w)*\zeta_z)/2$	&	$1 + \zeta_y + ((-1 + w)*\zeta_z)/2$	&	$(-2 + \zeta_z + w*\zeta_z)/2$\\
	\hline 	\hline
\textbf{Region} & \multicolumn{5}{c||}{\textbf{\textbf{Basis Function} Parameters}} \\ \hline 

	  & \textbf{Basis Function} & \textbf{Coefficient} & $\textbf{\text{arg}}_x$ & $\textbf{\text{arg}}_y$&$\textbf{\text{arg}}_z$ \\ 
	\hline
\textbf{6D}&$\chi^{S,L}$&$-\zeta_z/4$	&	$-1 + \zeta_x - ((1 + w)*\zeta_z)/2$	&	$1 + \zeta_y - ((1 + w)*\zeta_z)/2$	&	$(2 + \zeta_z - w*\zeta_z)/2$\\
	\hline
\textbf{6D}&$\chi^{S,LM}$&$(-\zeta_y + \zeta_z)/4$	&	$(-2 + 2*\zeta_x - (1 + w)*\zeta_y - \zeta_z + w*\zeta_z)/2$	&	$(2 + \zeta_y - w*\zeta_y + (-1 + w)*\zeta_z)/2$	&	$(-2 - (1 + w)*\zeta_y + \zeta_z + w*\zeta_z)/2$\\
	\hline
\textbf{6D}&$\chi^{S,RM}$&$(2 - \zeta_x + \zeta_y)/4$	&	$(\zeta_x - \zeta_y + w*(2 - \zeta_x + \zeta_y))/2$	&	$(-\zeta_x + \zeta_y + w*(2 - \zeta_x + \zeta_y))/2$	&	$(-\zeta_x - \zeta_y + w*(2 - \zeta_x + \zeta_y) + 2*\zeta_z)/2$\\
	\hline
\textbf{6D}&$\chi^{S,R}$&$\zeta_x/4$	&	$(-2 + \zeta_x + w*\zeta_x)/2$	&	$1 + ((-1 + w)*\zeta_x)/2 + \zeta_y$	&	$1 + ((-1 + w)*\zeta_x)/2 + \zeta_z$\\
	\hline \hline
\textbf{Region} & \multicolumn{5}{c||}{\textbf{\textbf{Basis Function} Parameters}} \\ \hline 
	
	  & \textbf{Basis Function} & \textbf{Coefficient} & $\textbf{\text{arg}}_x$ & $\textbf{\text{arg}}_y$&$\textbf{\text{arg}}_z$ \\ 
	\hline
\textbf{6E}&$\chi^{S,L}$&$-\zeta_y/4$	&	$-1 + \zeta_x - ((1 + w)*\zeta_y)/2$	&	$(2 + \zeta_y - w*\zeta_y)/2$	&	$1 - ((1 + w)*\zeta_y)/2 + \zeta_z$\\
	\hline
\textbf{6E}&$\chi^{S,LM}$&$(\zeta_y - \zeta_z)/4$	&	$(-2 + 2*\zeta_x + (-1 + w)*\zeta_y - \zeta_z - w*\zeta_z)/2$	&	$(-2 + (1 + w)*\zeta_y - (1 + w)*\zeta_z)/2$	&	$(2 + (-1 + w)*\zeta_y + \zeta_z - w*\zeta_z)/2$\\
	\hline
\textbf{6E}&$\chi^{S,RM}$&$(2 - \zeta_x + \zeta_z)/4$	&	$(\zeta_x - \zeta_z + w*(2 - \zeta_x + \zeta_z))/2$	&	$(-\zeta_x + 2*\zeta_y - \zeta_z + w*(2 - \zeta_x + \zeta_z))/2$	&	$(-\zeta_x + \zeta_z + w*(2 - \zeta_x + \zeta_z))/2$\\
	\hline
\textbf{6E}&$\chi^{S,R}$&$\zeta_x/4$	&	$(-2 + \zeta_x + w*\zeta_x)/2$	&	$1 + ((-1 + w)*\zeta_x)/2 + \zeta_y$	&	$1 + ((-1 + w)*\zeta_x)/2 + \zeta_z$\\
	\hline \hline
\textbf{Region} & \multicolumn{5}{c||}{\textbf{\textbf{Basis Function} Parameters}} \\ \hline 
	
	  & \textbf{Basis Function} & \textbf{Coefficient} & $\textbf{\text{arg}}_x$ & $\textbf{\text{arg}}_y$&$\textbf{\text{arg}}_z$ \\ 
	\hline
\textbf{7A}&$\chi^{S,L}$&$-\zeta_y/4$	&	$1 + \zeta_x - ((1 + w)*\zeta_y)/2$	&	$(2 + \zeta_y - w*\zeta_y)/2$	&	$-1 - ((1 + w)*\zeta_y)/2 + \zeta_z$\\
	\hline
\textbf{7A}&$\chi^{S,LM}$&$(-\zeta_x + \zeta_y)/4$	&	$(2 + \zeta_x - w*\zeta_x + (-1 + w)*\zeta_y)/2$	&	$(-2 - (1 + w)*\zeta_x + \zeta_y + w*\zeta_y)/2$	&	$-1 - ((1 + w)*\zeta_x)/2 + ((-1 + w)*\zeta_y)/2 + \zeta_z$\\
	\hline
\textbf{7A}&$\chi^{S,RM}$&$(2 + \zeta_x - \zeta_z)/4$	&	$(\zeta_x + w*(2 + \zeta_x - \zeta_z) - \zeta_z)/2$	&	$(-\zeta_x + 2*\zeta_y + w*(2 + \zeta_x - \zeta_z) - \zeta_z)/2$	&	$(-\zeta_x + w*(2 + \zeta_x - \zeta_z) + \zeta_z)/2$\\
	\hline
\textbf{7A}&$\chi^{S,R}$&$\zeta_z/4$	&	$1 + \zeta_x + ((-1 + w)*\zeta_z)/2$	&	$1 + \zeta_y + ((-1 + w)*\zeta_z)/2$	&	$(-2 + \zeta_z + w*\zeta_z)/2$\\
	\hline \hline
\textbf{Region} & \multicolumn{5}{c||}{\textbf{\textbf{Basis Function} Parameters}} \\ \hline 
	
	  & \textbf{Basis Function} & \textbf{Coefficient} & $\textbf{\text{arg}}_x$ & $\textbf{\text{arg}}_y$&$\textbf{\text{arg}}_z$ \\ 
	\hline
\textbf{7C}&$\chi^{S,L}$&$-\zeta_x/4$	&	$(2 + \zeta_x - w*\zeta_x)/2$	&	$1 - ((1 + w)*\zeta_x)/2 + \zeta_y$	&	$-1 - ((1 + w)*\zeta_x)/2 + \zeta_z$\\
	\hline
\textbf{7C}&$\chi^{S,LM}$&$(\zeta_x - \zeta_y)/4$	&	$(-2 + (1 + w)*\zeta_x - (1 + w)*\zeta_y)/2$	&	$(2 + (-1 + w)*\zeta_x + \zeta_y - w*\zeta_y)/2$	&	$-1 + ((-1 + w)*\zeta_x)/2 - ((1 + w)*\zeta_y)/2 + \zeta_z$\\
	\hline
\textbf{7C}&$\chi^{S,RM}$&$(2 + \zeta_y - \zeta_z)/4$	&	$(2*\zeta_x - \zeta_y + w*(2 + \zeta_y - \zeta_z) - \zeta_z)/2$	&	$(\zeta_y + w*(2 + \zeta_y - \zeta_z) - \zeta_z)/2$	&	$(-\zeta_y + w*(2 + \zeta_y - \zeta_z) + \zeta_z)/2$\\
	\hline
\textbf{7C}&$\chi^{S,R}$&$\zeta_z/4$	&	$1 + \zeta_x + ((-1 + w)*\zeta_z)/2$	&	$1 + \zeta_y + ((-1 + w)*\zeta_z)/2$	&	$(-2 + \zeta_z + w*\zeta_z)/2$\\
	\hline
\end{tabular} 
\end{adjustbox}

\end{landscape}

\newpage 
\begin{landscape}

\centering \begin{adjustbox}{width=1.0\linewidth}\centering
\begin{tabular}{||c|c|c|c|c|c||} \hline 
\textbf{Region} & \multicolumn{5}{c||}{\textbf{\textbf{Basis Function} Parameters}} \\ \hline 
	
	  & \textbf{Basis Function} & \textbf{Coefficient} & $\textbf{\text{arg}}_x$ & $\textbf{\text{arg}}_y$&$\textbf{\text{arg}}_z$ \\ 
	\hline
\textbf{8A}&$\chi^{S,L}$&$-\zeta_z/4$	&	$1 + \zeta_x - ((1 + w)*\zeta_z)/2$	&	$1 + \zeta_y - ((1 + w)*\zeta_z)/2$	&	$(2 + \zeta_z - w*\zeta_z)/2$\\
	\hline
\textbf{8A}&$\chi^{S,LM}$&$(-\zeta_y + \zeta_z)/4$	&	$(2 + 2*\zeta_x - (1 + w)*\zeta_y - \zeta_z + w*\zeta_z)/2$	&	$(2 + \zeta_y - w*\zeta_y + (-1 + w)*\zeta_z)/2$	&	$(-2 - (1 + w)*\zeta_y + \zeta_z + w*\zeta_z)/2$\\
	\hline
\textbf{8A}&$\chi^{S,RM}$&$(-\zeta_x + \zeta_y)/4$	&	$(2 + \zeta_x - w*\zeta_x + (-1 + w)*\zeta_y)/2$	&	$(-2 - (1 + w)*\zeta_x + \zeta_y + w*\zeta_y)/2$	&	$-1 - ((1 + w)*\zeta_x)/2 + ((-1 + w)*\zeta_y)/2 + \zeta_z$\\
	\hline
\textbf{8A}&$\chi^{S,R}$&$(2 + \zeta_x)/4$	&	$(\zeta_x + w*(2 + \zeta_x))/2$	&	$w - \zeta_x/2 + (w*\zeta_x)/2 + \zeta_y$	&	$w - \zeta_x/2 + (w*\zeta_x)/2 + \zeta_z$\\
	\hline \hline
\textbf{Region} & \multicolumn{5}{c||}{\textbf{\textbf{Basis Function} Parameters}} \\ \hline 
	
	  & \textbf{Basis Function} & \textbf{Coefficient} & $\textbf{\text{arg}}_x$ & $\textbf{\text{arg}}_y$&$\textbf{\text{arg}}_z$ \\ 
	\hline
\textbf{8B}&$\chi^{S,L}$&$-\zeta_y/4$	&	$1 + \zeta_x - ((1 + w)*\zeta_y)/2$	&	$(2 + \zeta_y - w*\zeta_y)/2$	&	$1 - ((1 + w)*\zeta_y)/2 + \zeta_z$\\
	\hline
\textbf{8B}&$\chi^{S,LM}$&$(\zeta_y - \zeta_z)/4$	&	$(2 + 2*\zeta_x + (-1 + w)*\zeta_y - \zeta_z - w*\zeta_z)/2$	&	$(-2 + (1 + w)*\zeta_y - (1 + w)*\zeta_z)/2$	&	$(2 + (-1 + w)*\zeta_y + \zeta_z - w*\zeta_z)/2$\\
	\hline
\textbf{8B}&$\chi^{S,RM}$&$(-\zeta_x + \zeta_z)/4$	&	$(2 + \zeta_x - w*\zeta_x + (-1 + w)*\zeta_z)/2$	&	$(-2 - (1 + w)*\zeta_x + 2*\zeta_y - \zeta_z + w*\zeta_z)/2$	&	$(-2 - (1 + w)*\zeta_x + \zeta_z + w*\zeta_z)/2$\\
	\hline
\textbf{8B}&$\chi^{S,R}$&$(2 + \zeta_x)/4$	&	$(\zeta_x + w*(2 + \zeta_x))/2$	&	$w - \zeta_x/2 + (w*\zeta_x)/2 + \zeta_y$	&	$w - \zeta_x/2 + (w*\zeta_x)/2 + \zeta_z$\\
	\hline \hline
\textbf{Region} & \multicolumn{5}{c||}{\textbf{\textbf{Basis Function} Parameters}} \\ \hline 

	  & \textbf{Basis Function} & \textbf{Coefficient} & $\textbf{\text{arg}}_x$ & $\textbf{\text{arg}}_y$&$\textbf{\text{arg}}_z$ \\ 
	\hline
\textbf{8C}&$\chi^{S,L}$&$-\zeta_z/4$	&	$1 + \zeta_x - ((1 + w)*\zeta_z)/2$	&	$1 + \zeta_y - ((1 + w)*\zeta_z)/2$	&	$(2 + \zeta_z - w*\zeta_z)/2$\\
	\hline
\textbf{8C}&$\chi^{S,LM}$&$(-\zeta_x + \zeta_z)/4$	&	$(2 + \zeta_x - w*\zeta_x + (-1 + w)*\zeta_z)/2$	&	$(2 - (1 + w)*\zeta_x + 2*\zeta_y - \zeta_z + w*\zeta_z)/2$	&	$(-2 - (1 + w)*\zeta_x + \zeta_z + w*\zeta_z)/2$\\
	\hline
\textbf{8C}&$\chi^{S,RM}$&$(\zeta_x - \zeta_y)/4$	&	$(-2 + (1 + w)*\zeta_x - (1 + w)*\zeta_y)/2$	&	$(2 + (-1 + w)*\zeta_x + \zeta_y - w*\zeta_y)/2$	&	$-1 + ((-1 + w)*\zeta_x)/2 - ((1 + w)*\zeta_y)/2 + \zeta_z$\\
	\hline
\textbf{8C}&$\chi^{S,R}$&$(2 + \zeta_y)/4$	&	$w + \zeta_x - \zeta_y/2 + (w*\zeta_y)/2$	&	$(\zeta_y + w*(2 + \zeta_y))/2$	&	$w - \zeta_y/2 + (w*\zeta_y)/2 + \zeta_z$\\
	\hline \hline
\textbf{Region} & \multicolumn{5}{c||}{\textbf{\textbf{Basis Function} Parameters}} \\ \hline 

	  & \textbf{Basis Function} & \textbf{Coefficient} & $\textbf{\text{arg}}_x$ & $\textbf{\text{arg}}_y$&$\textbf{\text{arg}}_z$ \\ 
	\hline
\textbf{8D}&$\chi^{S,L}$&$-\zeta_x/4$	&	$(2 + \zeta_x - w*\zeta_x)/2$	&	$1 - ((1 + w)*\zeta_x)/2 + \zeta_y$	&	$1 - ((1 + w)*\zeta_x)/2 + \zeta_z$\\
	\hline
\textbf{8D}&$\chi^{S,LM}$&$(\zeta_x - \zeta_z)/4$	&	$(-2 + (1 + w)*\zeta_x - (1 + w)*\zeta_z)/2$	&	$(2 + (-1 + w)*\zeta_x + 2*\zeta_y - \zeta_z - w*\zeta_z)/2$	&	$(2 + (-1 + w)*\zeta_x + \zeta_z - w*\zeta_z)/2$\\
	\hline
\textbf{8D}&$\chi^{S,RM}$&$(-\zeta_y + \zeta_z)/4$	&	$(-2 + 2*\zeta_x - (1 + w)*\zeta_y - \zeta_z + w*\zeta_z)/2$	&	$(2 + \zeta_y - w*\zeta_y + (-1 + w)*\zeta_z)/2$	&	$(-2 - (1 + w)*\zeta_y + \zeta_z + w*\zeta_z)/2$\\
	\hline
\textbf{8D}&$\chi^{S,R}$&$(2 + \zeta_y)/4$	&	$w + \zeta_x - \zeta_y/2 + (w*\zeta_y)/2$	&	$(\zeta_y + w*(2 + \zeta_y))/2$	&	$w - \zeta_y/2 + (w*\zeta_y)/2 + \zeta_z$\\
	\hline \hline
\textbf{Region} & \multicolumn{5}{c||}{\textbf{\textbf{Basis Function} Parameters}} \\ \hline 
	
	  & \textbf{Basis Function} & \textbf{Coefficient} & $\textbf{\text{arg}}_x$ & $\textbf{\text{arg}}_y$&$\textbf{\text{arg}}_z$ \\ 
	\hline
\textbf{8E}&$\chi^{S,L}$&$-\zeta_x/4$	&	$(2 + \zeta_x - w*\zeta_x)/2$	&	$1 - ((1 + w)*\zeta_x)/2 + \zeta_y$	&	$1 - ((1 + w)*\zeta_x)/2 + \zeta_z$\\
	\hline
\textbf{8E}&$\chi^{S,LM}$&$(\zeta_x - \zeta_y)/4$	&	$(-2 + (1 + w)*\zeta_x - (1 + w)*\zeta_y)/2$	&	$(2 + (-1 + w)*\zeta_x + \zeta_y - w*\zeta_y)/2$	&	$1 + ((-1 + w)*\zeta_x)/2 - ((1 + w)*\zeta_y)/2 + \zeta_z$\\
	\hline
\textbf{8E}&$\chi^{S,RM}$&$(\zeta_y - \zeta_z)/4$	&	$(-2 + 2*\zeta_x + (-1 + w)*\zeta_y - \zeta_z - w*\zeta_z)/2$	&	$(-2 + (1 + w)*\zeta_y - (1 + w)*\zeta_z)/2$	&	$(2 + (-1 + w)*\zeta_y + \zeta_z - w*\zeta_z)/2$\\
	\hline
\textbf{8E}&$\chi^{S,R}$&$(2 + \zeta_z)/4$	&	$w + \zeta_x - \zeta_z/2 + (w*\zeta_z)/2$	&	$w + \zeta_y - \zeta_z/2 + (w*\zeta_z)/2$	&	$w + ((1 + w)*\zeta_z)/2$\\
	\hline \hline
\textbf{Region} & \multicolumn{5}{c||}{\textbf{\textbf{Basis Function} Parameters}} \\ \hline 

	  & \textbf{Basis Function} & \textbf{Coefficient} & $\textbf{\text{arg}}_x$ & $\textbf{\text{arg}}_y$&$\textbf{\text{arg}}_z$ \\ 
	\hline
\textbf{8F}&$\chi^{S,L}$&$-\zeta_y/4$	&	$1 + \zeta_x - ((1 + w)*\zeta_y)/2$	&	$(2 + \zeta_y - w*\zeta_y)/2$	&	$1 - ((1 + w)*\zeta_y)/2 + \zeta_z$\\
	\hline
\textbf{8F}&$\chi^{S,LM}$&$(-\zeta_x + \zeta_y)/4$	&	$(2 + \zeta_x - w*\zeta_x + (-1 + w)*\zeta_y)/2$	&	$(-2 - (1 + w)*\zeta_x + \zeta_y + w*\zeta_y)/2$	&	$1 - ((1 + w)*\zeta_x)/2 + ((-1 + w)*\zeta_y)/2 + \zeta_z$\\
	\hline
\textbf{8F}&$\chi^{S,RM}$&$(\zeta_x - \zeta_z)/4$	&	$(-2 + (1 + w)*\zeta_x - (1 + w)*\zeta_z)/2$	&	$(-2 + (-1 + w)*\zeta_x + 2*\zeta_y - \zeta_z - w*\zeta_z)/2$	&	$(2 + (-1 + w)*\zeta_x + \zeta_z - w*\zeta_z)/2$\\
	\hline
\textbf{8F}&$\chi^{S,R}$&$(2 + \zeta_z)/4$	&	$w + \zeta_x - \zeta_z/2 + (w*\zeta_z)/2$	&	$w + \zeta_y - \zeta_z/2 + (w*\zeta_z)/2$	&	$w + ((1 + w)*\zeta_z)/2$\\
	\hline
\end{tabular} 
\end{adjustbox}

\end{landscape}

%% file: Supplementary_Material/3D_L-SIAC_coefficients.tex
  \newpage
  
  \begin{landscape}
  \begin{table}[h!]
      \centering

      \caption{Table describing how the shifting arguments ($q_x,q_y,q_z$) of the $a$ coefficients vary by region.}
      \label{tab:3D_LSIAC_coefficient}
   \begin{tabular}{||M{7.5cm}|M{3cm}|M{3cm}|M{3cm}|M{3cm}||}\hline
\textbf{Region} $\mathbf{S}$ $\mathbf{\Big{\backslash}}$ \textbf{Arguments of} $\mathbf{a_{\alpha}(q_x,q_y,q_z)}$ & $ \mathbf{P=L}$&$ \mathbf{P=LM}$&$ \mathbf{P=RM}$& $\mathbf{P=R}$\\ \hline
\textbf{1A}  &$(0,0,0)$   &  $(0,0,1)$&  $(0,1,1)$ &  $(1,1,1)$  \\ \hline
   
\textbf{1B}  &  $(0,0,0)$   &  $(0,1,0)$&  $(0,1,1)$ &  $(1,1,1)$  \\ \hline
    
\textbf{1C}  &  $(0,0,0)$   &  $(0,0,1)$&  $(1,0,1)$ &  $(1,1,1)$  \\ \hline
\textbf{1D}  &  $(0,0,0)$   &  $(1,0,0)$&  $(1,0,1)$ &  $(1,1,1)$  \\ \hline
\textbf{1E}  &  $(0,0,0)$   &  $(1,0,0)$&  $(1,1,0)$ &  $(1,1,1)$  \\ \hline
   
\textbf{1F}  &  $(0,0,0)$   &  $(0,1,0)$&  $(1,1,0)$ &  $(1,1,1)$  \\ \hline

\textbf{2E}  &  $(0,0,-1)$   &  $(0,0,0)$&  $(1,0,0)$ &  $(1,1,0)$  \\ \hline
   
\textbf{2F}  &  $(0,0,-1)$   &  $(0,0,0)$&  $(0,1,0)$ &  $(1,1,0)$  \\ \hline 

\textbf{3A}  &  $(-1,0,0)$   &  $(0,0,0)$&  $(0,0,1)$ &  $(0,1,1)$  \\ \hline
   
\textbf{3B}  &  $(-1,0,0)$   &  $(0,0,0)$&  $(0,1,0)$ &  $(0,1,1)$  \\ \hline 

\textbf{4B}  &  $(-1,0,-1)$   &  $(-1,0,0)$&  $(0,0,0)$ &  $(0,1,0)$  \\ \hline
   
\textbf{4F}  &  $(-1,0,-1)$   &  $(0,0,-1)$&  $(0,0,0)$ &  $(0,1,0)$  \\ \hline 

\textbf{5C}  &  $(0,-1,0)$   &  $(0,0,0)$&  $(0,0,1)$ &  $(1,0,1)$  \\ \hline
   
\textbf{5D}  &  $(0,-1,0)$   &  $(0,0,0)$&  $(1,0,0)$ &  $(1,0,1)$  \\ \hline 

\textbf{6D}  &  $(0,-1,-1)$   &  $(0,-1,0)$&  $(0,0,0)$ &  $(1,0,0)$  \\ \hline
   
\textbf{6E}  &  $(0,-1,-1)$   &  $(0,0,-1)$&  $(0,0,0)$ &  $(1,0,0)$  \\ \hline

\textbf{7A}  &  $(-1,-1,0)$   &  $(-1,0,0)$&  $(0,0,0)$ &  $(0,0,1)$  \\ \hline
   
\textbf{7C}  &  $(-1,-1,0)$   &  $(0,-1,0)$&  $(0,0,0)$ &  $(0,0,1)$  \\ \hline

\textbf{8A}  &  $(-1,-1,-1)$   &  $(-1,-1,0)$&  $(-1,0,0)$ &  $(0,0,0)$  \\ \hline
   
\textbf{8B}  &  $(-1,-1,-1)$   &  $(-1,0,-1)$&  $(-1,0,0)$ &  $(0,0,0)$  \\ \hline
    
\textbf{8C}  &  $(-1,-1,-1)$   &  $(-1,-1,0)$&  $(0,-1,0)$ &  $(0,0,0)$  \\ \hline
\textbf{8D}  &  $(-1,-1,-1)$   &  $(0,-1,-1)$&  $(0,-1,0)$ &  $(0,0,0)$  \\ \hline
\textbf{8E}  & $(-1,-1,-1)$   &  $(0,-1,-1)$&  $(0,0,-1)$ &  $(0,0,0)$  \\ \hline
   
\textbf{8F}  &  $(-1,-1,-1)$   &  $(-1,0,-1)$&  $(0,0,-1)$ &  $(0,0,0)$  \\ \hline

\end{tabular}

  \end{table}

  \end{landscape}